\documentclass[11pt,letterpaper]{article}

\usepackage{amssymb,amsmath,amsfonts}
\usepackage{graphicx,xcolor,enumitem}
\usepackage{epsfig}
\usepackage{amsthm}
\usepackage{bm}
\usepackage{multirow}
\usepackage{subcaption}
\usepackage[round]{natbib}
\usepackage{csquotes}
\renewcommand{\mkbegdispquote}[2]{\itshape}
\usepackage[twoside, hmarginratio=1:1, vmarginratio=1:1, left=1in,top=1in]{geometry}

\RequirePackage[breaklinks=true, hidelinks]{hyperref}
\usepackage{breakcites}

\usepackage{algorithm, algorithmic}

\usepackage{accents}

\newcommand{\ang}[1]{\left\langle  #1 \right\rangle }

\newcommand{\cB}{\mathcal{B}}
\newcommand{\E}{\mathbb{E}}
\newcommand{\R}{\mathbb{R}}
\newcommand{\p}{\mathbb{P}}
\newcommand{\q}{\mathbb{Q}}

\newcommand{\cS}{{\mathcal S}}

\newcommand{\cP}{{\mathcal P}}

\newcommand{\cV}{{\mathcal V}}
\newcommand{\cW}{{\mathcal W}}

\newcommand{\cX}{{\mathcal X}}
\newcommand{\cY}{{\mathcal Y}}
\newcommand{\cT}{{\mathcal T}}

\newcommand{\bN}{{\mathbb N}}

\newtheorem{theorem}{Theorem}

\newtheorem{assumption}[theorem]{Assumption}

\newtheorem{corollary}[theorem]{Corollary}

\newtheorem{definition}[theorem]{Definition}

\newtheorem{lemma}[theorem]{Lemma}

\newtheorem{remark}[theorem]{Remark}

\theoremstyle{definition}

\numberwithin{equation}{section}
\numberwithin{theorem}{section}

\begin{document}

\title{Equilibrium transport with time-inconsistent costs}

\author{
	Erhan Bayraktar \thanks{Department of Mathematics, University of Michigan, Ann Arbor, Email: erhan@umich.edu.}
	\and Bingyan Han \thanks{Thrust of Financial Technology, The Hong Kong University of Science and Technology (Guangzhou), Email: bingyanhan@hkust-gz.edu.cn.}
}

\date{April 9, 2025}
\maketitle

\begin{abstract}
    Given two probability measures on sequential data, we investigate the transport problem with time-inconsistent preferences in a discrete-time setting. Motivating examples are nonlinear objectives, state-dependent costs, and regularized optimal transport with general $f$-divergence. Under the bicausal constraint, we introduce the concept of equilibrium transport. Existence is proved in the semi-discrete Markovian case and the continuous non-Markovian case with strict quasiconvexity, while uniqueness also holds in the second case. We apply our framework to study mean-variance dynamic matching, nonlinear or state-dependent objectives with Gaussian data, and mismatches in job markets. Numerical results indicate a positive relationship between mismatches and state dependence.
	\\[2ex] 
	\\[2ex]
\end{abstract}

\section{Introduction}
Optimal transport (OT) is a method for quantifying the discrepancy between two probability distributions. The origins of OT can be traced back to Monge's mass transfer problem in the 18th century. However, it was not until the advent of linear programming in the mid 20th century that OT was revitalized by Kantorovich. \cite{brenier1991polar} and \cite{gangbo1996geometry} provided new insights into OT through the lens of convex analysis and geometry. Nowadays, OT is recognized as an interdisciplinary tool with a wide range of applications in fields such as mathematics \citep{villani2009optimal}, machine learning \citep{arjovsky2017wasserstein,peyre2019computational}, statistics \citep{torous2021optimal}, distributionally robust optimization \citep{blanchet2019quantifying,gao2016dist,mohajerin2018data}, and quantitative economics \citep{galichon2016optimal}.

While most literature on OT has focused on static data without a chronological structure, temporal data are ubiquitous in finance, statistics, and optimization. To address this, \cite{lassalle2013causal} introduced causal optimal transport (COT), which imposes a causality constraint on the transport plans between discrete stochastic processes. In a less formal sense, COT requires that when the past of one process, $X$, is given, the past of another process, $Y$, should be independent of the future of $X$ under the transport plan. COT has found applications in mathematical finance \citep{backhoff2020adapted}, video modeling \citep{xu2020cot}, mean-field games \citep{acciaio2021cournot,backhoff2023dynamic}, and stochastic optimization \citep{pflug2012distance,pflug2014multistage,acciaio2020causal}.

When making decisions over time, humans often exhibit time inconsistency, where the optimal plan at the present moment is not followed in the future. To address this, consistent planning is suggested as a remedy by  \cite{strotz1955myopia}. Time inconsistency is observed in non-exponential discounting \citep{strotz1955myopia,laibson1997golden}, prospect theory \citep{kahneman1979prospect}, gambling behaviors \citep{barberis2012model}, and mean-variance portfolio selection \citep{basak10}. Mathematical models with time inconsistency include variants of portfolio selection \citep{bjork14FS,bjork14MF,bjork17FS,han2021robust,kovavcova2021time}, optimal stopping \citep{bayraktar2021equilibrium,pichler2022risk}, robust decision-making \citep{epstein2022optimal}, probability distortion \citep{ma2021time}, and risk measures in \citet[Chapter 5]{pflug2014multistage} and \citet[Chapter 11]{follmer2011stochastic}.

In this paper, we generalize the classic Monge-Kantorovich problem to a dynamic setting with time-inconsistent costs. Suppose transport plans are bicausal in the sense that two processes are causal with respect to each other. The cost function is allowed to include time-inconsistent preferences such as nonlinear objectives and state dependence. Moreover, we notice that regularization with general $f$-divergence \citep{bayraktar2022stability,tacskesen2023semi,gonzalez2024quad} can also be a new source of time inconsistency. To address the self-contradictory behavior and the violation of dynamic programming principle (henceforth DPP), we adopt the subgame Nash equilibria solution as in \cite{bjork14FS,basak10}. This concept, given formally in Definition \ref{def:equitrans}, is well-defined and easy to understand in a discrete-time setting.  

Our main technical contribution lies in characterizing equilibrium transport in Definition \ref{def:equitrans} with an extended DPP framework. Lemma \ref{lem:semi_V} and Theorem \ref{thm:lsc} establish the existence of equilibrium transport under a semi-discrete setting with Markovian dynamics. The assumptions imposed in Theorem \ref{thm:lsc} are comparable to the classic assumptions in \citet[Theorem 4.1]{villani2009optimal} and straightforward to verify. Our proof involves refining the Polish topology recursively, drawing on \citet[Theorem 13.11]{kechris2012classical}. Several fundamental claims need to be checked carefully under the refined topology. A crucial topological result is that the weak topology on $\Pi(\mu(dx_{t+1}|x_t), \nu(dy_{t+1} | y_t))$ relies solely on the Borel structure, rather than the topological one. Our proof technique is still new in the well-studied literature on measurable selection and DPP. Theorem \ref{thm:params} proves the existence and uniqueness of parametric couplings in a continuous, non-Markovian setting under a strict quasiconvexity assumption, leveraging Berge's maximum theorem. Section \ref{sec:ex_params} provides a concrete example illustrating Theorem \ref{thm:params}. In contrast to the existence of equilibrium controls \citep{bayraktar2023existence}, the proof techniques in this paper are tailored for OT and rely heavily on the topological and Borel structures of probability spaces.

Our framework introduces a novel methodology with broad interdisciplinary applications. We provide three illustrative examples in Section \ref{sec:illustrate} and analyze mismatches in job markets in Sections \ref{sec:reduced_jobm} and \ref{sec:struct_jobm}:
\begin{enumerate}[label={(\arabic*)}]
	\item Section \ref{sec:mv} explores dynamic matching between supply and demand types with a mean-variance objective. The types are horizontally differentiated, where closer type matches result in lower costs. The probability of mismatches in the next period depends on the current state. Numerically, finding an equilibrium transport is more straightforward than determining a pre-committed OT, as the former is a local optimizer and suitable for parallelization. An example in Section \ref{sec:gaussian} involving Gaussian data has explicit and unique equilibrium parametric transport, while the pre-committed optimal parametric transport is not unique.
	
	\item Roberts' law states that CEO compensation is approximately proportional to $(\text{own firm size})^\kappa$ with $\kappa \simeq 1/3$. However, empirical data sometimes deviate from this relationship, motivating us to investigate why Roberts' law may not always apply. We conjecture that factors such as regulatory constraints, financial performance momentum, and status quo bias contribute to dependence on previous states. Sections \ref{sec:reduced_jobm} and \ref{sec:struct_jobm} develop two job market models incorporating state-dependent terms to capture persistence or inertia in maintaining previous matches. Inspired by reduced-form methods \citep{abowd1999high,borovickova2020high}, Section \ref{sec:reduced_jobm} directly examines the relationship between sales and wages. In contrast, Section \ref{sec:struct_jobm} adopts a structural approach, in the spirit of \cite{gabaix2008has,galichon2016optimal}.
		
	We find a positive relationship between mismatch and inertia as captured by state dependence. This observation is supported by both synthetic and empirical analyses in Sections \ref{sec:exec}, \ref{sec:prof}, and \ref{sec:Cobb_num}. Notably, the postdoctoral job market exhibits even stronger inertia than the executive labor market across industries, a phenomenon often overlooked in the existing literature.

\end{enumerate}

The rest of the paper is organized as follows. We introduce the usual weak topology and bicausal OT in Sections \ref{sec:notation} and \ref{sec:bic}. Section \ref{sec:ti} presents motivations for equilibrium transport and characterization in the discrete case. Section \ref{sec:semi} considers the semi-discrete and Markovian case. Section \ref{sec:params} presents the continuous and non-Markovian case. Section \ref{sec:illustrate} gives three illustrative examples. Section \ref{sec:Roberts} revisits Roberts' law and reviews two main approaches in the literature on wage dispersion. Sections \ref{sec:reduced_jobm} and \ref{sec:struct_jobm} develop two models to analyze state dependence in job markets. Technical proofs are deferred to Section \ref{sec:proofs} in the e-companion. Sections \ref{sec:supp_exec} and \ref{sec:supp_acad} provide details of the numerical implementation. The Python code is available at: \url{https://github.com/hanbingyan/equitrans}.

\subsection{Notation and the usual weak topology}\label{sec:notation} 
Denote the finite number of periods as $T$. For each $t \in \{1, \ldots, T\}$, consider a Polish (complete, separable, and metrizable) space $(\cX_t, \cT_{\cX_t})$, where $\cX_t$ is a closed (but not necessarily bounded) subset of $\R^d$. $\cX_t$ is interpreted as the range of the process at time $t$. Let $C(\cX_t; \cT_{\cX_t})$ be the set of continuous functions $f: (\cX_t, \cT_{\cX_t}) \rightarrow (\R, \cT_\R)$, where $\R$ is always equipped with the usual topology $\cT_\R$.  $C_b(\cX_t; \cT_{\cX_t})$ is the set of continuous and bounded functions. Denote $\cB(\cT_{\cX_t})$ as the Borel $\sigma$-algebra of the topological space $(\cX_t, \cT_{\cX_t})$.

In this paper, we always endow the product space with the product topology. $\cX_{1:T} := \cX_1 \times \ldots \times \cX_T$ is a closed subset of $\R^{T \times d}$ with the product topology $\cT_\cX := \prod^{T}_{t=1} \cT_{\cX_t}$, which is also Polish. With some abuse of notation, we simply write $\cX$ while we mean $\cX_{1:T}$. By \citet[Lemma 1.2]{kallenberg}, the Borel $\sigma$-algebra of the product space satisfies $\cB(\cX) = \cB(\cX_1) \otimes \cdots \otimes \cB(\cX_T)$.  

Denote the set of all Borel probability measures on $\cX$ as $\cP(\cX)$. As in \citet[Section 7.4.2]{bertsekas1978stoch} and \citet[Section II.6]{parthasarathy2005}, the usual weak topology is defined as follows. For a given $\varepsilon > 0$, $\mu \in \cP(\cX)$, and $f \in C_b(\cX; \cT_\cX)$, define the subset of $\cP(\cX)$ as
\begin{equation*}
	V(\mu; f, \varepsilon) := \left\{ \mu' \in \cP(\cX) : \left| \int f d\mu - \int f d\mu' \right| < \varepsilon \right\}. 
\end{equation*}
Consider the collection of these subsets as
\begin{equation*}
	V[C_b(\cX; \cT_\cX)] := \left\{ V(\mu; f, \varepsilon) : \varepsilon >0, \, \mu \in \cP(\cX), \, f \in C_b(\cX; \cT_\cX) \right\}.
\end{equation*}
We endow $\cP(\cX)$ with $\cV[C_b(\cX; \cT_\cX)]$, the weakest topology on $\cP(\cX)$ which contains the collection $V[C_b(\cX; \cT_\cX)]$. By \citet[Theorems 6.2 and 6.5]{parthasarathy2005}, $(\cP(\cX), \cV[C_b(\cX; \cT_\cX)])$ is a Polish topological space if and only if $(\cX, \cT_\cX)$ is so.

In the OT theory, we need to consider another closed set $\cY = \cY_{1:T} = \cY_1 \times ... \times \cY_T$. We also equip each factor space $\cY_t$ with a Polish topology $\cT_{\cY_t}$ and other notations are introduced similarly as in the counterparts for $\cX$.  

For notational convenience, we interpret $\cX_{1:t}$, $\cX_{1:t} \times \cY_{1:t}$, and other similar terms as the product of factor spaces with indices in $1, ..., t$. For a Borel probability measure $\mu$ on temporal data, we denote $\mu(dx_{t+1:T}|x_{1:t})$ as a regular conditional probability kernel of $x_{t+1:T}$ given $x_{1:t}$, which is uniquely determined in a suitable way \cite[Theorem 10.4.14 and Corollary 10.4.17]{bogachev2007measure}. When there is no confusion, we also write the kernel as $\mu^t$ for simplicity. If $\mu$ is a finite discrete measure, we sometimes use $\mu(x_{t+1:T}|x_{1:t})$, which omits the differential, as the conditional probability at the atom $x_{t+1:T}$ given $x_{1:t}$. The initial state at time $0$ is always fixed.

We emphasize that metrics and Wasserstein distances are introduced for the continuous case in Section \ref{sec:params} only. The semi-discrete case in Section \ref{sec:semi} does not rely on a specific choice of the metric.




\subsection{Bicausal optimal transport}\label{sec:bic}

There is a vast literature on the OT theory and its applications \citep{villani2009optimal,kuhn2019tutorial,blanchet2021stat}. However, the previous works have been focusing on data without the time dimension. Motivated by the ubiquitous role of temporal data in OR/MS, finance, and economics, a new notion of adapted Wasserstein distance or causal optimal transport (COT) has been proposed to compare distributions on temporal data \citep{lassalle2013causal,backhoff2017causal}. We present a brief review in this subsection.

Consider two probability measures $\mu \in \cP(\cX)$ and $\nu \in \cP(\cY)$. Denote $\Pi(\mu, \nu)$ as the set of all the couplings that admit $\mu$ and $\nu$ as marginals. Suppose transporting one unit of mass from $x$ to $y$ incurs a cost of $c(x,y)$. A generic OT problem is formulated as 
\begin{equation*}
	\cW(\mu, \nu) := \inf_{ \pi \in \Pi(\mu, \nu)} \int_{\cX \times \cY} c(x, y) \pi(dx, dy).
\end{equation*}

If the data have a temporal structure as $x = (x_1, ..., x_t, ..., x_T)$ and $y = (y_1, ..., y_t, ..., y_T)$, not all couplings $\pi \in \Pi(\mu, \nu)$ will make sense. A natural requirement of the transport plan $\pi(x, y)$ should be the non-anticipative condition. Informally speaking, if the past of $x$ is given, then the past of $y$ should be independent of the future of $x$ under the measure $\pi$. Mathematically, it means a transport plan $\pi$ should satisfy
\begin{equation}\label{eq:causal}
	\pi (dy_t | x_{1:T}) = \pi (dy_t | x_{1:t}), \quad  t = 1, ..., T-1, \quad \pi\text{-a.s.}
\end{equation}
The property \eqref{eq:causal} is known as the causality condition from $x$ to $y$ and the transport plan satisfying \eqref{eq:causal} is called {\it causal} by \cite{lassalle2013causal}. It can be interpreted by the equivalent formulation: $y_t = F_t(x_{1:t}, U_t)$, $\forall \, t \in{1, \ldots, T}$, for some measurable functions $F_t$ and uniform random variables $U_t$. Crucially, $U_t$ is independent of $x_{1:T}$, but the $\{U_t\}_{t=1}^T$ need not to be independent of each other; see \citet[Theorem 8.17]{kallenberg} and \cite{backhoff2017causal}.

If \eqref{eq:causal} holds when we exchange the positions of $x$ and $y$, then the transport plan is called bicausal. Denote $\Pi_{bc}(\mu, \nu)$ as the set of all bicausal transport plans between $\mu$ and $\nu$. The bicausal OT problem considers the optimization over $\Pi_{bc}(\mu, \nu)$ only:
\begin{equation}\label{eq:bc}
	\cW_{bc}(\mu, \nu) := \inf_{ \pi \in \Pi_{bc}(\mu, \nu)} \int_{\cX \times \cY} c(x, y) \pi(dx, dy).
\end{equation}
For applications of (bi)causal OT, see \cite{backhoff2020adapted,xu2020cot,acciaio2021cournot,pflug2012distance,acciaio2020causal} for an incomplete list.

In this paper, we consider two probability measures $\mu$ and $\nu$ with symmetric positions and focus on the bicausal transport plans only. For a given transport plan $\pi \in \Pi(\mu, \nu)$, we can decompose $\pi$ in terms of successive regular kernels:
\begin{align}
	\pi(dx_{1:T}, dy_{1:T}) = & \bar{\pi} (dx_1, dy_1) \prod^{T-1}_{s=1} \pi(dx_{s+1}, dy_{s+1} | x_{1:s}, y_{1:s}). \label{eq:decomp}
\end{align}
By \citet[Proposition 5.1]{backhoff2017causal}, $\pi$ is a bicausal transport plan if and only if
\begin{itemize}
	\item[(1)] $\bar{\pi} \in \Pi(p^1_* \mu, p^1_* \nu)$, and
	\item[(2)] for each $t = 1, ... , T-1$ and $\pi$-almost every path $(x_{1:t}, y_{1:t})$, the following condition holds:
	\begin{equation*}
		\pi(dx_{t+1}, dy_{t+1}| x_{1:t}, y_{1:t}) \in \Pi( \mu(dx_{t+1}| x_{1:t}), \nu(dy_{t+1}|y_{1:t})). 
	\end{equation*} 
\end{itemize}
$p^1_* \mu$ (resp. $p^1_* \nu$) is the pushforward of $\mu$ (resp. $\nu$) by the projection $p^1$ onto the first coordinate.

\section{Equilibrium transport}\label{sec:ti}

Time-consistent problems satisfy the Bellman equation. If a solution is optimal on the time interval $\{t, ..., T-1\}$, then it is also optimal on any subinterval $\{ s, ..., T-1\}$ with $s \geq t$. However, there are various time-inconsistent preferences \citep{strotz1955myopia,kahneman1979prospect,laibson1997golden,basak10}. In this context, we present three instances where time inconsistency arises in the bicausal OT problem, with the regularization case being particularly noteworthy and previously unnoticed.

\subsection{Motivation}
\subsubsection{Regularized bicausal OT} 
In practice, continuous densities are approximated by empirical measures with finite supports. The computational burden of discrete OT can be high, then regularization is adopted in implementation \citep{cuturi2013sinkhorn,pichler2022nested,eckstein2024,gonzalez2024quad}. The discrete bicausal OT with regularization is
\begin{equation}\label{eq:div}
	\inf_{\pi \in \Pi_{bc}(\mu, \nu)} \ang{c(x_{1:T}, y_{1:T}), \pi(x_{1:T}, y_{1:T})}_F + \ang{ f \left( \frac{\pi(x_{1:T}, y_{1:T})}{\mu(x_{1:T}) \otimes \nu(y_{1:T})} \right), \mu(x_{1:T}) \otimes \nu(y_{1:T})}_F
\end{equation}
for a convex function $f$. The last term is known as the $f$-divergence. Here, $$\ang{c(x_{1:T}, y_{1:T}), \pi(x_{1:T}, y_{1:T})}_F = \sum_{x_{1:T} \in \cX, y_{1:T} \in \cY} c(x_{1:T}, y_{1:T})\pi(x_{1:T}, y_{1:T}),$$ 
which sums over all paths. The arguments $(x_{1:T}, y_{1:T})$ in $\ang{c(x_{1:T}, y_{1:T}), \pi(x_{1:T}, y_{1:T})}_F$ should be interpreted as dummy variables similar in the integral \eqref{eq:bc}. $\mu \otimes \nu$ is the independent coupling and $\frac{\pi}{\mu \otimes \nu}$ is interpreted as an element-wise division.

A popular choice of $f(\cdot)$ is the Kullback--Leibler (KL) divergence given by $f(x) = x \ln(x)$. In this case, the regularized objective \eqref{eq:div} is still time-consistent. Indeed, the KL divergence is separable in the sense that 
\begin{align}
	& \text{KL} (\pi(x_{1:T}, y_{1:T}) \| \mu(x_{1:T}) \otimes \nu(y_{1:T}) ) \nonumber \\
	& \qquad = \ang{ \ln \left( \frac{\pi(x_{1:T}, y_{1:T})}{\mu(x_{1:T}) \otimes \nu(y_{1:T})} \right), \pi(x_{1:T}, y_{1:T})}_F \nonumber  \\
	& \qquad = \Big\langle \Big\langle \ln \left( \frac{\pi(x_T, y_T | x_{1:T-1}, y_{1:T-1})}{\mu(x_T | x_{1:T-1}) \otimes \nu(y_T | y_{1:T-1})} \right), \pi(x_T, y_T | x_{1:T-1}, y_{1:T-1}) \Big\rangle_F \label{eq:div-decom} \\
	& \qquad \qquad + \ln \left( \frac{\pi(x_{1:T-1}, y_{1:T-1})}{\mu(x_{1:T-1}) \otimes \nu(y_{1:T-1})} \right), \pi(x_{1:T-1}, y_{1:T-1})  \Big \rangle_F.  \nonumber 
\end{align}
For simplicity, we denote the successive regular kernels of $\pi(x_{1:T}, y_{1:T})$ as $\{ \bar{\pi}, \pi_1, ..., \pi_{T-1} \}$. For example, $\pi_{T-1} = \pi(x_T, y_T | x_{1:T-1}, y_{1:T-1})$. The second equality \eqref{eq:div-decom} shows that the optimization over $\pi_{T-1}$ is only determined by the first term. We can repeat the decomposition on the second term recursively over time. The classic DPP is applicable in this setting.  If $\{ \pi^*_t, ..., \pi^*_{T-1} \}$ is optimal on the interval $\{t, ..., T-1\}$, then $\{ \pi^*_s, ..., \pi^*_{T-1} \}$ is also optimal on the subinterval  $\{s, ..., T-1\}$ with $s > t$.

However, the entropic regularization can be numerically unstable and fails to preserve transport plan sparsity, prompting exploration of alternative choices for $f(\cdot)$ as discussed in \cite{bayraktar2022stability}. Notably, \citet[Proposition 3.5 and Theorem 3.7]{tacskesen2023semi} proved that some $f$-divergence regularization in the primal problem is equivalent to the so-called smooth $c$-transform in the dual problem.

Some non-separable examples with time inconsistency are
\begin{itemize}
	\item[(1)] Squared Hellinger distance: $f(x) = (1 - \sqrt{x})^2$. The $f$-divergence term becomes
	\begin{align*}
		& \ang{f \left( \frac{\pi(x_{1:T}, y_{1:T})}{\mu(x_{1:T}) \otimes \nu(y_{1:T})} \right),  \mu(x_{1:T}) \otimes \nu(y_{1:T})}_F = 2 - 2 \ang{\sqrt{\pi(x_{1:T}, y_{1:T})}, \sqrt{\mu(x_{1:T}) \otimes \nu(y_{1:T})}}_F.
	\end{align*}
	\item[(2)] Le Cam distance: $f(x) = \frac{(x - 1)^2}{2x + 2}$ and
	\begin{align*}
		& \ang{ f \left( \frac{\pi(x_{1:T}, y_{1:T})}{\mu(x_{1:T}) \otimes \nu(y_{1:T})} \right), \mu(x_{1:T}) \otimes \nu(y_{1:T})}_F  = \frac{1}{2} \sum_{x_{1:T}, y_{1:T}} \frac{[\pi(x_{1:T}, y_{1:T}) - \mu(x_{1:T}) \otimes \nu(y_{1:T})]^2}{\pi(x_{1:T}, y_{1:T}) + \mu(x_{1:T}) \otimes \nu(y_{1:T})}.
	\end{align*}
	The summation is over all possible paths $x_{1:T}$ and $y_{1:T}$.
	\item[(3)] Jensen-Shannon divergence: $f(x) = x \ln (\frac{2x}{x+1}) + \ln(\frac{2}{x+1})$, and 
	\begin{equation*}
		\text{JS}(\p, \q) = \text{KL} \left(\p \Big\| \frac{\p + \q}{2} \right) + \text{KL} \left(\q \Big\| \frac{\p + \q}{2} \right),
	\end{equation*}
	with $\p = \pi$, $\q = \mu \otimes \nu$ is the independent coupling.
\end{itemize}
Under these specifications, the objective is no longer separable and the DPP is violated. The cost functional at time $t$ is given by
\begin{align}
	J(x_{1:t}, y_{1:t}; \pi) =& \ang{c(x_{1:T}, y_{1:T}), \pi(x_{t+1:T}, y_{t+1:T}| x_{1:t}, y_{1:t})}_F \label{eq:div-obj} \\
	& + \ang{ f \left( \frac{\pi(x_{t+1:T}, y_{t+1:T}| x_{1:t}, y_{1:t})}{\mu(x_{t+1:T} | x_{1:t}) \otimes \nu(y_{t+1:T} |y_{1:t})} \right), \mu(x_{t+1:T}|x_{1:t}) \otimes \nu(y_{t+1:T}|y_{1:t})}_F. \nonumber
\end{align}
When the iterated expectation relationship (or tower property) fails, a global OT plan on $\{t, ..., T-1\}$ may not be optimal for the subproblems on $\{t+1, ..., T-1\}$. The intuition is that we may sacrifice the optimality on the subinterval $\{t+1, ..., T-1\}$ to attain better solutions on $\{t, ..., T-1\}$. Nevertheless, an agent can still choose to ignore time inconsistency and follow the global optimizer at time $0$ irrevocably. However, he should recognize that the ``today self'' and ``future selves'' have an intertemporal conflict with the ``optimal" plan \citep{strotz1955myopia}.

\begin{remark}
	This paper focuses on existence and uniqueness of solutions introduced in Definition \ref{def:equitrans}. Convergence of numerical algorithms for the transport problem with time-inconsistent costs remains an open problem. Notably, even in the static case with quadratic regularization, this problem has only been examined in a recent work by \cite{gonzalez2024quad}.
\end{remark}

\subsubsection{Nonlinear objectives}
With a nonlinear function $G$, suppose the agent would like to minimize the transport cost given by
\begin{equation}\label{eq:nonlin_G}
	\inf_{\pi \in \Pi_{bc}(\mu, \nu)} G \left(\int h(x_{1:T}, y_{1:T}) \pi(dx_{1:T}, dy_{1:T}) \right).
\end{equation}
An  illustrative example is
\begin{equation*}
	\inf_{\pi \in \Pi_{bc}(\mu, \nu)}  \left(\int h(x_{1:T}, y_{1:T}) \pi(dx_{1:T}, dy_{1:T}) - m \right)^2.
\end{equation*}
The agent wants to match the expected value of $h(\cdot, \cdot)$ to a given level of $m$. Time inconsistency appears since the tower property fails for the objective.


\subsubsection{State-dependent preference} 
For a fixed $t = 0, ... , T-1$, consider an objective functional 
\begin{align} \label{eq:stateobj}
	J(x_{1:t}, y_{1:t}; \pi) =& \int c(x_t, y_t, x_{1:T}, y_{1:T}) \pi(dx_{t+1:T}, dy_{t+1:T} | x_{1:t}, y_{1:t}).
\end{align}
The cost $c$ has four arguments. The first $(x_t, y_t)$ is state-dependent and plays a different role than the counterpart in $x_{1:T}$ and $y_{1:T}$. That is, for a different time $s \neq t$, the integrand becomes $c(x_s, y_s, x_{1:T}, y_{1:T})$. Therefore, the integrand is different when the time changes and it leads to time inconsistency.  

A classic example is non-exponential discounting \citep{strotz1955myopia,laibson1997golden}. At time $t$, consider
\begin{equation}\label{eq:nonexp}
	\inf_{\pi \in \Pi_{bc}(\mu(dx_{t+1:T}|x_{1:t}), \nu(dy_{t+1:T}|y_{1:t}))} \int  \sum^T_{s = t+1} \varphi(s - t) c(x_s, y_s) \pi(dx_{t+1:T}, dy_{t+1:T} | x_{1:t}, y_{1:t}),
\end{equation}
where $\varphi(\cdot)$ is a discounting function besides the power function. Note that the form \eqref{eq:nonexp} can be reformulated into \eqref{eq:stateobj} by redefining a new $\tilde x_t = (t, x_t)$. Since transitions of time are deterministic, the reformulation does not change the conditional probability kernels.


\subsection{Definition of equilibrium transport}
To handle time inconsistency and define a reasonable concept of ``optimal" solutions, we should consider consistent plans that the agent can follow, instead of pre-committed solutions that are optimal only at the fixed initial point. The agent who is aware of time inconsistency should keep in mind that the objectives in the future are different. Inspired by the concept of subgame perfect Nash equilibrium, several works \citep{strotz1955myopia,bjork14FS,bjork17FS} propose to reformulate the problem as a non-cooperative game between agents at time $t+1, ..., T-1$ that are incarnations of the agent at time $t$, that is, a game between the current self and the future selves. An equilibrium solution, denoted as $\{ \pi^*_t, ..., \pi^*_{T-1} \}$ in a discrete-time setting, should satisfy the following property: If ourselves in the future time $t+1, ... , T-1$ stick with $\{ \pi^*_{t+1}, ..., \pi^*_{T-1} \}$, then it is optimal for us at the current time $t$ to adopt $\pi^*_t$. This backward recursive definition guarantees that the agent at each time $t$ will not deviate from $\pi^*_t$.    

Inspired by the above discussion, we define the equilibrium transport as follows. We refer to a transport plan $\pi$ by its successive regular kernels $\{ \bar{\pi}, \pi_1, ..., \pi_{T-1} \}$ in \eqref{eq:decomp}. At time $t \in \{0, ..., T-1\}$, denote a generic cost functional as $J(x_{1:t}, y_{1:t}; \pi)$, which implies that the current path is $(x_{1:t}, y_{1:t})$ and the agent uses the transport plan $\pi$ in time $t, ..., T-1$.
\begin{definition}\label{def:equitrans}
	Consider a given bicausal transport plan $\pi^*$ with successive regular kernels $\{ \bar{\pi}^*, \pi^*_1, $ $..., \pi^*_{T-1} \}$.
	\begin{itemize}
		\item[1.] For any $t = 1, ..., T-1$ and $\pi^*$-almost every path $(x_{1:t}, y_{1:t})$, consider any transport plan $\gamma \in \Pi( \mu(dx_{t+1}| x_{1:t}), \nu(dy_{t+1}|y_{1:t}) )$. When $t=0$, we adopt the convention that $\cX_0$ and $\cY_0$ are singletons and consider $\bar{\gamma} \in \Pi(p^1_* \mu, p^1_* \nu)$.
		\item[2.] Define a perturbed transport plan $\pi^{t, \gamma}$ as
		\begin{align*}
			\pi^{t, \gamma} (dx_{t+1:T}, dy_{t+1:T} | x_{1:t}, y_{1:t}) =& \gamma( dx_{t+1}, dy_{t+1} | x_{1:t}, y_{1:t}) \prod^{T-1}_{s=t+1} \pi^*(dx_{s+1}, dy_{s+1} | x_{1:s}, y_{1:s}).
		\end{align*}
	\end{itemize}
	Then $\pi^*$ is a subgame perfect Nash equilibrium bicausal transport plan if for every $t = 0, ..., T-1$ and $\pi^*$-almost every $(x_{1:t}, y_{1:t})$, we have
	\begin{equation*}
		\inf_{\gamma \in \Pi(\mu(dx_{t+1}|x_{1:t}), \nu(dy_{t+1} | y_{1:t}) )} J(x_{1:t}, y_{1:t}; \pi^{t, \gamma}) = J(x_{1:t}, y_{1:t}; \pi^{t, *}),
	\end{equation*}
	where
	\begin{align*}
		\pi^{t, *} (dx_{t+1:T}, dy_{t+1:T} | x_{1:t}, y_{1:t}) =& \prod^{T-1}_{s=t} \pi^*(dx_{s+1}, dy_{s+1} | x_{1:s}, y_{1:s}).
	\end{align*}
	We call $\pi^*$ an equilibrium transport for simplicity. If $\pi^*$ exists, define the equilibrium value function $V$ as
	\begin{equation*}
		V_t (x_{1:t}, y_{1:t}) = J(x_{1:t}, y_{1:t}; \pi^{t, *}).
	\end{equation*}
\end{definition}

Next, we compare equilibrium transport with pre-committed OT. Since $\cX_0$ and $\cY_0$ are singletons, the initial state at time $0$ is fixed. The cost functional at time $0$ is denoted as $J_0(\pi)$, where the initial state is omitted. 
	\begin{definition}\label{def:pre-OT}
		Suppose there is a bicausal optimizer, denoted as $\pi^{pre}$, for the following problem with a fixed initial point:
		\begin{equation}\label{prob:precommitted}
			\inf_{\pi \in \Pi_{bc} (\mu, \nu)} J_0(\pi).
		\end{equation}
		Then $\pi^{pre}$ is called a pre-committed OT solution to \eqref{prob:precommitted} at time $0$.
	\end{definition}
	
	Consider the nonlinear objective in \eqref{eq:nonlin_G} as an example. Then, $$J_0(\pi) = G(\int h(x_{1:T}, y_{1:T}) \pi(dx_{1:T}, dy_{1:T})).$$ A pre-committed OT solution, if it exists, minimizes $J_0(\pi)$ over all $\pi \in \Pi_{bc}(\mu, \nu)$. While pre-committed OT solutions are globally optimal, they overlook the issue of time inconsistency. In contrast, equilibrium transport plans may not minimize $J_0(\pi)$, but they are subgame optimal and time-consistent.

\subsection{The discrete case}
To understand Definition \ref{def:equitrans} with an example, we revisit the regularized objective \eqref{eq:div-obj} and characterize an equilibrium transport as follows. If time moves to $t+1$, the cost functional is given by
\begin{align*}
	& J(x_{1:t+1}, y_{1:t+1}; \pi) = \ang{c(x_{1:T}, y_{1:T}), \pi(x_{t+2:T}, y_{t+2:T}| x_{1:t+1}, y_{1:t+1})}_F \\
	& \qquad  + \ang{ f \left( \frac{\pi(x_{t+2:T}, y_{t+2:T}| x_{1:t+1}, y_{1:t+1})}{\mu(x_{t+2:T} | x_{1:t+1}) \otimes \nu(y_{t+2:T} |y_{1:t+1})} \right), \mu(x_{t+2:T}|x_{1:t+1}) \otimes \nu(y_{t+2:T}|y_{1:t+1})}_F.
\end{align*}
If all agents at time $t+1, ..., T-1$ adopt the equilibrium transport $\pi^*$, then
\begin{equation*}
	J(x_{1:t+1}, y_{1:t+1}; \pi^*) = V_{t+1} (x_{1:t+1}, y_{1:t+1}).
\end{equation*} 
Next, we argue by backward induction and consider the subproblem at time $t$ with a perturbed transport plan $\pi^{t, \gamma} = \pi^*(x_{t+2:T}, y_{t+2:T}| x_{1:t+1}, y_{1:t+1}) \gamma(x_{t+1}, y_{t+1})$. We have
\begin{align*}
	& J(x_{1:t}, y_{1:t}; \pi^{t, \gamma})  = \ang{c(x_{1:T}, y_{1:T}), \pi^*(x_{t+2:T}, y_{t+2:T}| x_{1:t+1}, y_{1:t+1}) \gamma(x_{t+1}, y_{t+1})}_F \\
	& \qquad + \ang{ f \left( \frac{\pi^*(x_{t+2:T}, y_{t+2:T}| x_{1:t+1}, y_{1:t+1}) \gamma(x_{t+1}, y_{t+1})}{\mu(x_{t+1:T} | x_{1:t}) \otimes \nu(y_{t+1:T} |y_{1:t})} \right), \mu(x_{t+1:T}|x_{1:t}) \otimes \nu(y_{t+1:T}|y_{1:t})}_F.
\end{align*}

Therefore, combining the two equations, we obtain the following recursive relationship between $J(x_{1:t}, y_{1:t}; \pi^{t, \gamma})$ and $J(x_{1:t+1}, y_{1:t+1}; \pi^*)$:
\begin{align*}
	&J (x_{1:t}, y_{1:t}; \pi^{t, \gamma}) = \ang{J(x_{1:t+1}, y_{1:t+1}; \pi^*), \gamma(x_{t+1}, y_{t+1})}_F \\
	& \quad -  \ang{ f \left( \frac{\pi^*(x_{t+2:T}, y_{t+2:T}| x_{1:t+1}, y_{1:t+1})}{\mu(x_{t+2:T} | x_{1:t+1}) \otimes \nu(y_{t+2:T} |y_{1:t+1})} \right), \mu(x_{t+2:T}|x_{1:t+1}) \otimes \nu(y_{t+2:T}|y_{1:t+1}) \gamma(x_{t+1}, y_{t+1})}_F   \\
	& \quad + \ang{ f \left( \frac{\pi^*(x_{t+2:T}, y_{t+2:T}| x_{1:t+1}, y_{1:t+1}) \gamma(x_{t+1}, y_{t+1})}{\mu(x_{t+1:T} | x_{1:t}) \otimes \nu(y_{t+1:T} |y_{1:t})} \right), \mu(x_{t+1:T}|x_{1:t}) \otimes \nu(y_{t+1:T}|y_{1:t})}_F.
\end{align*}

By Definition \ref{def:equitrans}, if all agents at time $t+1, ..., T-1$ adopt the equilibrium transport $\{\pi^*_{t+1}, ..., $ $\pi^*_{T-1} \}$, then it is also optimal for the agent at time $t$ to follow $\pi^*_t$. Therefore, the recursive relationship leads to the following equation for the value function:
\begin{align*}
	& V_t (x_{1:t}, y_{1:t}) = \inf_{\gamma \in \Pi(\mu(x_{t+1}|x_{1:t}), \nu(y_{t+1} | y_{1:t}) )} \Big[ \ang{V_{t+1} (x_{1:t+1}, y_{1:t+1}), \gamma(x_{t+1}, y_{t+1})}_F \\ 
	&\quad  - \ang{f \left( \frac{\pi^*(x_{t+2:T}, y_{t+2:T} | x_{1:t+1}, y_{1:t+1})}{\mu(x_{t+2:T} | x_{1:t+1}) \otimes \nu(y_{t+2:T}| y_{1:t+1})} \right), \mu(x_{t+2:T} | x_{1:t+1}) \otimes \nu(y_{t+2:T}| y_{1:t+1}) \gamma(x_{t+1}, y_{t+1})}_F \\
	&\quad  + \ang{f \left( \frac{\pi^*(x_{t+2:T}, y_{t+2:T} | x_{1:t+1}, y_{1:t+1}) \gamma(x_{t+1}, y_{t+1})}{\mu(x_{t+1:T} | x_{1:t}) \otimes \nu(y_{t+1:T}| y_{1:t})} \right), \mu(x_{t+1:T} | x_{1:t}) \otimes \nu(y_{t+1:T}| y_{1:t})}_F  \Big]. 
\end{align*}
We follow \cite{bjork14FS,bjork17FS} and call the recursive relationship of the value function $V_t$ above the extended dynamic programming (DP) equation. 
Since the problem is finite and discrete, the infimum over $\gamma$ is attained and induces the conditional kernel $\pi^*_t$ of an equilibrium transport. Moreover, there can be multiple equilibrium transport plans, since the extended DP equation may have more than one optimizer. The value function is continuous automatically under the discrete topology. For the discrete case with general nonlinear and state-dependent objectives, see Remark \ref{rem:discrete} in Section \ref{sec:para}.

As a side remark, if we consider the regularized objective \eqref{eq:div-obj} with continuous measures $\mu$ and/or $\nu$ instead, the usual weak topology in Section \ref{sec:notation} is not strong enough to prove the continuity of value function. Even convergence of probabilities under Wasserstein distance does not necessarily lead to convergence under $f$-divergence, including Jensen-Shannon, KL, reverse KL, total variation, and other divergences. See \citet[Example 1]{arjovsky2017wasserstein}. Since the continuity of value function is needed for the recursion, we can only show the existence of equilibrium transport in the discrete case for the regularized objective \eqref{eq:div-obj} under the current weak topology.

\section{The semi-discrete and Markovian case}\label{sec:semi}
In this section, we work with topological spaces without fixing any particular metrics. Lemma \ref{lem:semi_V} and Theorem \ref{thm:lsc} provide the existence of equilibrium transport when $\cX$ is a finite discrete set with the discrete topology and $\cY$ can be a general closed set in $\R^{T \times d}$. The current formulation is already computationally challenging. Indeed, \citet[Theorem 2.2]{tacskesen2023semi} showed that the computational complexity of the Wasserstein distance between a discrete probability measure supported on two points and the Lebesgue measure on the standard hypercube is already $\#$P-hard. 

Suppose the nonlinear and state-dependent objective is given by  
\begin{align}
	J(x_{t}, y_{t}; \pi) :=& \int \sum^T_{k=t+1} c_k(x_t, y_t, x_k, y_k) \pi(dx_{t+1:T}, dy_{t+1:T} | x_{t}, y_{t}) \nonumber \\
	& + G \left(x_t, y_t, \int h(x_{T}, y_{T}) \pi(dx_{T}, dy_{T} | x_{t}, y_{t}) \right). \label{eq:semi_obj}
\end{align} 
The function $h$ only depends on $(x_T, y_T)$. The cost function $c$ is separable and given by
\begin{equation*}
	c(w, v, x_{1:T}, y_{1:T}) := \sum^T_{k=1} c_k(w, v, x_k, y_k),
\end{equation*}
where $w$ and $v$ capture the state dependence on $x$ and $y$. The first two terms in $G(x_t, y_t, \cdot )$ are also state-dependent. 

By convention, we interpret $\cX_0$ and $\cY_0$ as singletons.
\begin{assumption}\label{a:semi_obj}
	\begin{itemize}
		\item[(1)] For each $t \in \{1, ..., T\}$, $\cX_t$ is a non-empty finite discrete set equipped with the discrete topology $\cT_{\cX_t}$. And $\cX$ is endowed with the product discrete topology.
		\item[(2)] For each $k \in \{1, ..., T\}$ and $i \in \{0, ..., T\}$, $c_k(x_i, y_i, x_k, y_k):  (\cX_i \times \cY_i \times \cX_k \times \cY_k, \cT_{\cX_i} \times \cT_{\cY_i} \times \cT_{\cX_k} \times \cT_{\cY_k}) \rightarrow (\R, \cT_\R)$ is continuous and bounded. Hence, when $i=k$, $c_k(x_k, y_k, x_k, y_k)$ is $\cT_{\cX_k} \times \cT_{\cY_k}$-continuous and bounded.
		\item[(3)] $h(x_T, y_T): (\cX_T \times \cY_T, \cT_{\cX_T} \times \cT_{\cY_T}) \rightarrow (\R, \cT_\R)$ is continuous and bounded.
		\item[(4)] For each $t \in \{0, ..., T\}$, $G(x_t, y_t, g): (\cX_t \times \cY_t \times \R, \cT_{\cX_t} \times \cT_{\cY_t} \times \cT_\R ) \rightarrow (\R, \cT_\R)$ is continuous and bounded. 
	\end{itemize}
\end{assumption}

\begin{assumption}\label{a:semi_meas}
	For each $t \in \{1, ..., T-1\}$,
	\begin{itemize}
		\item[(1)] the regular conditional kernels of $\mu$ and $\nu$ are Markovian and denoted as $\mu(dx_{t+1}| x_{t})$ and $\nu(dy_{t+1}| y_{t})$.
		\item[(2)] $\nu(dy_{t+1}|y_t): (\cY_t, \cT_{\cY_t}) \rightarrow (\cP(\cY_{t+1}),  \cV[C_b(\cY_{t+1}; \cT_{\cY_{t+1}})])$ is Borel measurable.
	\end{itemize}
\end{assumption}
Here, the Markov property means that the kernels depend on the current $x_t$ (or $y_t$) instead of the whole path $x_{1:t}$ (or $y_{1:t}$). We note that $\mu(dx_{t+1}| x_t)$ is continuous with respect to the discrete topology.

Under the current formulation, the extended DP equation is also Markovian:
\begin{align}
	V_t (x_{t}, y_{t}) =   \inf_{\gamma \in \Pi(\mu(dx_{t+1}|x_{t}), \nu(dy_{t+1} | y_{t}) )} & \Big[ \int \Big(c_{t+1} (x_t, y_t, x_{t+1}, y_{t+1}) +  V_{t+1} (x_{t+1}, y_{t+1}) \Big) \gamma(dx_{t+1}, dy_{t+1}) \nonumber \\ 
	& - \int G(x_{t+1}, y_{t+1}, g_{t+1} (x_{t+1}, y_{t+1})) \gamma(dx_{t+1}, dy_{t+1})  \label{eq:semi_dpp} \\
	& + G\Big(x_t, y_t, \int g_{t+1} (x_{t+1}, y_{t+1}) \gamma(dx_{t+1}, dy_{t+1}) \Big) \nonumber \\
	& - \int \sum^T_{k=t+2} b_k (x_{t+1}, y_{t+1}, x_{t+1}, y_{t+1})  \gamma(dx_{t+1}, dy_{t+1})  \nonumber \\
	& +  \int \sum^T_{k=t+2} b_k (x_t, y_t, x_{t+1}, y_{t+1})  \gamma(dx_{t+1}, dy_{t+1})  \Big]. \nonumber 
\end{align}

In this equation,
\begin{itemize}
	\item[(a)] the boundary condition for $V$ is
	\begin{equation*}
		V_T (x_{T}, y_{T}) =  G \left(x_T, y_T, h(x_{T}, y_{T}) \right);
	\end{equation*}
	\item[(b)] let $\pi^*(dx_{t+1}, dy_{t+1} | x_{t}, y_{t})$ denote an optimizer for the equation \eqref{eq:semi_dpp} at time $t$. Concatenate these kernels to obtain
	\begin{align}
		\pi^*(dx_{t+1:T}, dy_{t+1:T} | x_{t}, y_{t}) = & \prod^{T-1}_{s=t} \pi^*(dx_{s+1}, dy_{s+1} | x_{s}, y_{s}); \label{eq:semipi*}
	\end{align}	
	\item[(c)] the function sequences $g_{t+1}$ and $b_{k}$ in \eqref{eq:semi_dpp} are given by
	\begin{align*}
		g_{t+1} (x_{t+1}, y_{t+1}) & := \int h(x_{T}, y_{T}) \pi^*(dx_{T}, dy_{T} | x_{t+1}, y_{t+1}), \\
		b_{k} (w, v, x_{t+1}, y_{t+1}) & := \int c_k(w, v, x_k, y_k) \pi^*(dx_{k}, dy_{k} | x_{t+1}, y_{t+1}),
	\end{align*}
	with $\pi^*(dx_k, dy_k | x_{t+1}, y_{t+1})$ obtained by backward induction using \eqref{eq:semipi*}. Here, the time $t \in \{0, \ldots, T-2\}$ and the index $k \in \{t+2, ..., T\}$. The state-dependent terms $w \in \cX_i$ and $v \in \cY_i$, with $i \in \{0, ..., t+1\}$.
\end{itemize}

There is an essential difficulty to prove the existence of equilibrium transport. Recall the Borel space framework in \citet[Chapters 7 and 8]{bertsekas1978stoch} for stochastic optimal control. The model treated there has Borel state, control, and disturbance spaces; see \citet[Definition 8.1]{bertsekas1978stoch}. The existence of optimal policies relies on measurable selection theorems. The optimizer usually has weaker properties compared with the objective function. For example, if the objective is lower semicontinuous (l.s.c.), then an optimizer, if it exists, is Borel measurable \cite[Proposition 7.33]{bertsekas1978stoch}. If the objective is lower semianalytic, then an $\varepsilon$-optimizer is universally measurable \cite[Proposition 7.50]{bertsekas1978stoch}. However, the extended DP equation introduces auxiliary functions $g_{t+1}$ and $b_k$, tied to the optimizer, and thus inheriting the weaker properties. Consequently, demonstrating the existence of an equilibrium via a recursive application of measurable selection becomes challenging.

Our main idea to prove Lemma \ref{lem:semi_V} is the Borel measurable selection theorem \citep{brown1973} but with a careful treatment of the Polish topology \citep[Theorem 13.11]{kechris2012classical}. The proof is given in Section \ref{sec:proof_semi} of the e-companion. \citet[Theorem 13.11]{kechris2012classical} shows that there is a finer Polish topology with the same Borel sets, such that a Borel function is continuous under the finer topology. We apply this technique recursively. However, one has to be careful to verify that several fundamental claims are still valid under the new topology:
	\begin{enumerate}[label={(\arabic*)}]
		\item Tightness of transport plans \cite[Lemma 4.4]{villani2009optimal} uses the fact that the topology on $\cX_{t+1} \times \cY_{t+1}$ is a product topology. But \citet[Theorem 13.11]{kechris2012classical} does not guarantee that the refined topology on $\cX_{t+1} \times \cY_{t+1}$ is still a product topology in general. Therefore, we imposed that $\cX_{t+1}$ is equipped with the discrete topology and refine the topology on $\cY_{t+1}$. The Markov property is imposed for a similar reason, to ensure that $\cY_{t+1:T}$ has a product topology after applying the refining technique.
		\item Recall that compact Hausdorff topologies $\cT$ are both minimal Hausdorff and maximal compact \cite[Section I.3, p.25]{steen1978}. In other words, no topology strictly smaller than $\cT$ can be Hausdorff and no topology strictly larger than $\cT$ can be compact. Consequently, it is crucial to emphasize that $\Pi(\mu(dx_{t+1}|x_t), \nu(dy_{t+1} | y_t))$ remains compact after refining the Polish topology. Indeed, the weak topology on $\Pi(\mu(dx_{t+1}|x_t), \nu(dy_{t+1} | y_t))$ only depends on the Borel structure of $\cX_{t+1}$ and $\cY_{t+1}$, instead of the topological one. It means that if we replace the topologies on $\cX_{t+1}$ and $\cY_{t+1}$ with new ones having the same Borel sets, then the weak topology on $\Pi(\mu(dx_{t+1}|x_t), \nu(dy_{t+1} | y_t))$ is unchanged. This technical topological result is also observed in \citet[Lemma 2.3 and Remark 2.4]{beiglbock2012duality}.
	\end{enumerate}	
	To the best of our knowledge, the refining topology technique is still new in the established literature on measurable selection and DP, dating back to \cite{brown1973,bertsekas1978stoch}.

\begin{lemma}\label{lem:semi_V}
	Suppose Assumptions \ref{a:semi_obj} and \ref{a:semi_meas} hold. Then for $t \in \{0, \ldots, T-1\}$, there exists a Borel measurable optimizer $\pi^*(dx_{t+1}, dy_{t+1} | x_{t}, y_{t})$  for \eqref{eq:semi_dpp}. Moreover, $\pi^*$ defined recursively by \eqref{eq:semipi*} is an equilibrium transport in the sense of Definition \ref{def:equitrans}. The value function $V_t(x_{t}, y_{t})$ in \eqref{eq:semi_dpp} is Borel measurable in $(x_t, y_t)$. Moreover, there exists a finer Polish topology such that $\pi^*(dx_{t+1}, dy_{t+1} | x_{t}, y_{t})$ and $V_t(x_{t}, y_{t})$ are continuous in $(x_t, y_t)$.
\end{lemma}

Note that for $t \in \{0, \ldots, T-2\}$,
\begin{align*}
	V_{t+1} (x_{t+1}, y_{t+1}) = & \sum^T_{k=t+2} b_k (x_{t+1}, y_{t+1}, x_{t+1}, y_{t+1}) + G(x_{t+1}, y_{t+1}, g_{t+1} (x_{t+1}, y_{t+1})). 		
\end{align*}
Then the extended DP equation \eqref{eq:semi_dpp} reduces to
\begin{align}
	V_t (x_{t}, y_{t}) =   \inf_{\gamma \in \Pi(\mu(dx_{t+1}|x_{t}), \nu(dy_{t+1} | y_{t}) )} & \Big[ \int c_{t+1} (x_t, y_t, x_{t+1}, y_{t+1})  \gamma(dx_{t+1}, dy_{t+1}) \label{eq:semi_dpp_simple} \\
	& + G\Big(x_t, y_t, \int g_{t+1} (x_{t+1}, y_{t+1}) \gamma(dx_{t+1}, dy_{t+1}) \Big) \nonumber \\
	& +  \int \sum^T_{k=t+2} b_k (x_t, y_t, x_{t+1}, y_{t+1})  \gamma(dx_{t+1}, dy_{t+1})  \Big]. \nonumber 
\end{align}
The original expression \eqref{eq:semi_dpp} demonstrates the derivation of the value function $V_t$ from $V_{t+1}$, establishing a recursive connection that aligns closely with \cite{bjork14FS}. Under the continuity assumption, the proof of Lemma \ref{lem:semi_V} can rely on \eqref{eq:semi_dpp} and the reduction \eqref{eq:semi_dpp_simple} is not needed. However, when the costs are l.s.c. only, it becomes unclear whether \eqref{eq:semi_dpp} can maintain the l.s.c. property and the reduction \eqref{eq:semi_dpp_simple} is needed. Both \eqref{eq:semi_dpp} and  \eqref{eq:semi_dpp_simple} are used to give more intuitive proofs.
\begin{assumption}\label{a:semi_lsc}
	\begin{itemize}
		\item[(1)] Assumption \ref{a:semi_obj} (1) holds.
		\item[(2)] For each $k \in \{1, ..., T\}$ and $i \in \{0, ..., T\}$, $c_k(x_i, y_i, x_k, y_k):  (\cX_i \times \cY_i \times \cX_k \times \cY_k, \cT_{\cX_i} \times \cT_{\cY_i} \times \cT_{\cX_k} \times \cT_{\cY_k}) \rightarrow (\R, \cT_\R)$ is l.s.c. and bounded from below. Hence, when $i=k$, $c_k(x_k, y_k, x_k, y_k)$ is $\cT_{\cX_k} \times \cT_{\cY_k}$-l.s.c. and bounded from below.
		\item[(3)] $h(x_T, y_T): (\cX_T \times \cY_T, \cT_{\cX_T} \times \cT_{\cY_T}) \rightarrow (\R, \cT_\R)$ is l.s.c. and bounded from below.
		\item[(4)] For each $t \in \{0, ..., T\}$, $G(x_t, y_t, g): (\cX_t \times \cY_t \times \R, \cT_{\cX_t} \times \cT_{\cY_t} \times \cT_\R ) \rightarrow (\R, \cT_\R)$ is l.s.c. and $G(x_t, y_t, \cdot)$ is nondecreasing for each $(x_t, y_t)$. 
	\end{itemize}
\end{assumption}
Theorem \ref{thm:lsc} extends Lemma \ref{lem:semi_V} to l.s.c. costs. Hence, Theorem \ref{thm:lsc} is comparable to the classic existence of an optimal coupling \cite[Theorem 4.1]{villani2009optimal}. Moreover, it is direct to verify the additional assumptions such as $G(x_t, y_t, \cdot)$ is nondecreasing for each $(x_t, y_t)$.
\begin{theorem}\label{thm:lsc}
	Suppose Assumptions \ref{a:semi_meas} and \ref{a:semi_lsc} hold. Then for $t \in \{0, \ldots, T-1\}$, there exists a Borel measurable optimizer $\pi^*(dx_{t+1}, dy_{t+1} | x_{t}, y_{t})$  for \eqref{eq:semi_dpp} or equivalently \eqref{eq:semi_dpp_simple}. Moreover, $\pi^*$ defined recursively by \eqref{eq:semipi*} is an equilibrium transport in the sense of Definition \ref{def:equitrans}. The value function $V_t(x_{t}, y_{t})$ in \eqref{eq:semi_dpp} is Borel measurable in $(x_t, y_t)$. Moreover, there exists a finer Polish topology such that $\pi^*(dx_{t+1}, dy_{t+1} | x_{t}, y_{t})$ is continuous in $(x_t, y_t)$ and $V_t(x_{t}, y_{t})$ is l.s.c. in $(x_t, y_t)$.
\end{theorem}

\section{Continuous and non-Markovian cases with parametric couplings}\label{sec:params}
\subsection{Metric spaces and Wasserstein distances}\label{sec:metrics}
For generic probability measures and unbounded cost functions, metric spaces are needed for growth rate conditions; see Lemma \ref{lem:metric_cont} in the e-companion. Consider $\cX$ and $\cY$ with the Euclidean topology. Given a fixed $p \in [1, \infty)$, introduce the metric as $d_\cX(x_{1:T}, x'_{1:T}) = \left[ \sum^T_{t=1} d_{\cX_t}(x_t, x'_t)^p \right]^{1/p}$ for $x, x' \in \cX$, where $d_{\cX_t}(x_t, x'_t) := |x_t - x'_t|$. The Wasserstein space of order $p$ is given by
\begin{equation*}
	\cP_p(\cX) := \left\{ \mu \in \cP(\cX) \Big| \int_{\cX} d_\cX(x_{1:T}, \bar{x}_{1:T})^p \mu(dx) < \infty \right\}
\end{equation*}
for some fixed $\bar{x}_{1:T} \in \cX$. We always equip $\cP_p(\cX)$ with the Wasserstein distance of order $p$:
\begin{equation}\label{eq:Wp}
	W_p(\mu, \mu') = \left( \inf_{\gamma \in \Pi(\mu, \mu')} \int_{\cX \times \cX} d_\cX(x, x')^p \gamma (dx, dx') \right)^{1/p}.
\end{equation}
Let
\begin{equation*}
	C_p(\cX) := \left\{ f \in C(\cX) \Big|  \exists\, \text{ constant } C > 0, \; |f(x)| \leq C(1 + d_\cX(x, \bar{x})^p)\right\}
\end{equation*}
be the set of continuous functions from $\cX$ to $\R$ with a growth rate of order $p$ in $d_\cX(x, \bar{x})$. Consider a sequence of probability measures $\{ \mu^k \}_{k \in \bN}$ in $\cP_p(\cX)$ and $\mu$ another probability measure in $\cP_p(\cX)$. Recall that $\mu^k$ converges weakly in $\cP_p(\cX)$ to $\mu$ means that $\mu^k$ converges to $\mu$ in the usual weak convergence and $\int d_\cX(x, \bar{x})^p \mu^k(dx) \rightarrow \int d_\cX(x, \bar{x})^p \mu(dx)$; see \citet[Definition 6.7]{villani2009optimal}. By \citet[Theorem 6.8]{villani2009optimal}, $W_p(\mu^k, \mu) \rightarrow 0$ is equivalent to $\mu^k$ converges weakly in $\cP_p(\cX)$ to $\mu$.

Similarly, we can replace $(\cX, d_\cX)$ with any metric space $(\cS, d_\cS)$ that is Polish and define $\cP_p(\cS)$ and $C_p(\cS)$ accordingly. For example, $\cS = \cX_{t+1}$ or $\cS = \cX_{t+1} \times \cY_{t+1}$. In particular, we denote the metric for $\cX \times \cY$ as
\begin{align*}
	d((x, y), (\bar{x}, \bar{y})) = [d_{\cX}(x, \bar{x})^p + d_{\cY}(y, \bar{y})^p]^{1/p}.
\end{align*}



\subsection{Parametric couplings}\label{sec:para}

For each time $t \in \{0, ..., T-1\}$, we use $ \gamma(dx_{t+1}, dy_{t+1}|\theta_{t+1})$ to model elements in a subset of $\cP_p(\cX_{t+1} \times \cY_{t+1})$. Suppose the parameter $\theta_{t+1}$ is in a Polish topological vector space $(\Theta_{t+1}, \cT_{\Theta_{t+1}})$ with a complete compatible metric $d_{\Theta_{t+1}}$. The parametric model $\gamma$ may not enumerate every element in $\cP_p(\cX_{t+1} \times \cY_{t+1})$ and can be regarded as a form of dimension reduction. For example, the agent may want to narrow down the candidate couplings and consider specific parametric distributions only. The dual potential functions are commonly modeled by neural networks, which provide a parametric form for the couplings if regularization is also considered; see \citet[Equation (8)]{seguy2018large}. Another method in \cite{delon2020wasserstein} restricts to Gaussian mixtures models as the couplings. Parametric methods can reduce the computational burden while sacrifice the accuracy. If we take $\Theta_{t+1} = \cP_p(\cX_{t+1} \times \cY_{t+1})$, it collapses to the classic formulation. We term the parametric kernel as the parametric transport (coupling) and the concatenation is denoted as
\begin{align}\label{eq:params_pi}
	& \gamma(dx_{t+1:T}, dy_{t+1:T} | \theta_{t+1:T})  := \prod^{T-1}_{s=t} \gamma(dx_{s+1}, dy_{s+1} | \theta_{s+1}).
\end{align} 

Consider a general nonlinear and state-dependent objective which is not necessarily separable:
\begin{align}
	J(x_{1:t}, y_{1:t}; \theta_{t+1:T} ) :=& \int c(x_t, y_t, x_{1:T}, y_{1:T}) \gamma(dx_{t+1:T}, dy_{t+1:T} | \theta_{t+1:T}) \nonumber \\
	& + G \left(x_t, y_t, \int h(x_{1:T}, y_{1:T}) \gamma(dx_{t+1:T}, dy_{t+1:T} | \theta_{t+1:T}) \right). \label{eq:nonMark_Obj}
\end{align} 
The minimization is over $\theta_{t+1:T}$, which should satisfy the bicausal constraints $\gamma(dx_{t+1:T}, dy_{t+1:T} | $ $\theta_{t+1:T}) \in \Pi_{bc}(\mu(dx_{t+1:T}|x_{1:t}), \nu(dy_{t+1:T}|y_{1:t}))$. 

We impose the following assumptions on the objective and probability measures.
\begin{assumption}\label{a:params_obj}
	\begin{itemize}
		\item[(1)]  For each $t \in \{ 0, ..., T\}$, $c(x_t, y_t, x_{1:T}, y_{1:T}) \in C_p(\cX_{1:T} \times \cY_{1:T})$.
		\item[(2)] $h(\cdot, \cdot)$ is continuous and
		$$|h(x_{1:T}, y_{1:T})| \leq C(1 + d((x_{1:T}, y_{1:T}), (\bar{x}_{1:T}, \bar{y}_{1:T}))^{1/r})$$
		for some constant $r>0$ satisfying $1/r \leq p$.
		\item[(3)] For each $t \in \{ 0, ..., T\}$, $G(x_t, y_t, g):$ $(\cX_t \times \cY_t \times \R, \cT_{\cX_t} \times \cT_{\cY_t} \times \cT_\R) \rightarrow (\R, \cT_\R)$ is continuous and $|G(x_t, y_t, g)| \leq C [1 + |g|^{pr} + d((x_t, y_t), (\bar{x}_t, \bar{y}_t))^p]$.
	\end{itemize}
\end{assumption}
\begin{assumption}\label{a:measure}
	\begin{itemize}
		\item[(1)] 	For each $t \in \{0, ..., T-1\}$, $x_{1:t} \mapsto \mu(dx_{t+1}|x_{1:t})$ is continuous with respect to metric $d_{\cX}$ on the domain $\cX_{1:t}$ and the $W_p$ metric on the range $\cP_p (\cX_{t+1})$.
		\item[(2)] For each $t \in \{0, ..., T-1\}$ and constant $a \in \{1/r,\, p\}$, where the constant $r>0$ is the same as in Assumption \ref{a:params_obj},
		\begin{equation*}
			\int_{\cX_{t+1:T}} d_\cX( x_{1:T}, \bar{x}_{1:T})^a \mu(dx_{t+1:T}|x_{1:t}) \leq C[1 +  d_\cX(x_{1:t}, \bar{x}_{1:t})^a].
		\end{equation*} 
	\end{itemize}
	Same assumptions on $\nu$ hold with the metric $d_\cY$.
\end{assumption}
Assumption \ref{a:measure}(2) is used in \eqref{eq:a} for the growth rate.

Introduce the correspondence $D_t: \cX_{1:t} \times \cY_{1:t} \twoheadrightarrow \Theta_{t+1}$ as
\begin{equation} \label{eq:paramD}
	(x_{1:t}, y_{1:t}) \mapsto \{ \theta_{t+1} \in \Theta_{t+1} |  \gamma(dx_{t+1}, dy_{t+1}|\theta_{t+1}) \in  \Pi(\mu(dx_{t+1}|x_{1:t}), \nu(dy_{t+1} | y_{1:t})) \}.	
\end{equation}
We use a notation $\twoheadrightarrow$ to highlight that the correspondence maps a point to a subset.

The equilibrium parametric transport is characterized by:
\begin{align}
	V_t (x_{1:t}, y_{1:t}) =   \inf_{\theta_{t+1} \in D_t(x_{1:t}, y_{1:t})} f(x_{1:t}, y_{1:t}, \theta_{t+1}), \label{eq:params_dpp}  
\end{align}
where
\begin{align}
	f(x_{1:t}, y_{1:t}, \theta_{t+1})  :=  & G\Big(x_t, y_t, \int g_{t+1} (x_{1:t+1}, y_{1:t+1}) \gamma(dx_{t+1}, dy_{t+1}| \theta_{t+1}) \Big) \label{eq:para_t_obj} \\
	& + \int b_{t+1} (x_t, y_t, x_{1:t+1}, y_{1:t+1}) \gamma(dx_{t+1}, dy_{t+1} | \theta_{t+1}). \nonumber 
\end{align}

In this equation,
\begin{itemize}
	\item[(a)] the boundary condition for $V$ is
	\begin{equation*}
		V_T (x_{1:T}, y_{1:T}) = c(x_T, y_T, x_{1:T}, y_{1:T}) + G \left(x_T, y_T, h(x_{1:T}, y_{1:T}) \right);
	\end{equation*}
	\item[(b)] denote an optimizer for the equation \eqref{eq:params_dpp} at time $t$ as  $\theta^*_{t+1}(x_{1:t}, y_{1:t})$. It yields the equilibrium parametric kernels in the sense of Definition \ref{def:equitrans} when kernels are restricted in the parametric spaces: 
	\begin{align}
		& \gamma(dx_{t+1:T}, dy_{t+1:T} | \theta^*_{t+1:T}(x_{1:t}, y_{1:t})) = \prod^{T-1}_{s=t} \gamma(dx_{s+1}, dy_{s+1} | \theta^*_{s+1}(x_{1:s}, y_{1:s}));  \label{eq:pi*params} 
	\end{align}	
	\item[(c)] the function sequences $g_{t+1}$ and $b_{t+1}$ in \eqref{eq:params_dpp}  are given by
	\begin{align*}
		g_{t+1} (x_{1:t+1}, y_{1:t+1}) & := \int h(x_{1:T}, y_{1:T}) \gamma(dx_{t+2:T}, dy_{t+2:T} | \theta^*_{t+2:T}(x_{1:t+1}, y_{1:t+1})), \\
		b_{t+1} (w, v, x_{1:t+1}, y_{1:t+1}) & := \int c(w, v, x_{1:T}, y_{1:T}) \gamma(dx_{t+2:T}, dy_{t+2:T} | \theta^*_{t+2:T}(x_{1:t+1}, y_{1:t+1})),
	\end{align*}
	with $\gamma(dx_{t+2:T}, dy_{t+2:T} | \theta^*_{t+2:T}(x_{1:t+1}, y_{1:t+1}))$ in \eqref{eq:pi*params} and $\theta^*_{i+1}(x_{1:i}, y_{1:i})$, $i = t+1, ..., T-1$, are obtained by backward induction using \eqref{eq:params_dpp}.
\end{itemize} 

Recall that a function $f: \cS \rightarrow \R$ defined on a convex subset $\cS$ of a real vector space is quasiconvex if for all $s, s' \in \cS$ and $\lambda \in [0, 1]$, we have
$f(\lambda s + (1- \lambda) s') \leq \max\{f(s), f(s')\}$. If furthermore $f(\lambda s + (1 - \lambda) s') < \max\{f(s), f(s')\}$ for all $s \neq s'$ and $\lambda \in (0, 1)$, then $f$ is strictly quasiconvex. Every convex function is quasiconvex. Besides, we recall the definition of upper and lower hemicontinuity:
	\begin{definition}{\cite[Definition 17.2]{charalambos2013infinite}}
		A correspondence $\varphi: X \twoheadrightarrow Y$ between topological spaces is:
		\begin{itemize}
			\item upper hemicontinuous at the point $x$ if for every neighborhood $U$ of $\varphi(x)$, there is a neighborhood $V$ of $x$ such that $z \in V$ implies $\varphi(z) \subset U$.
			\item lower hemicontinuous at $x$ if for every open set $U$ that meets $\varphi(x)$ (i.e. $\varphi(x) \cap U \neq \emptyset$), there is a neighborhood $V$ of $x$ such that $z \in V$ implies $\varphi(z) \cap U \neq \emptyset$.
			\item continuous at $x$ if it is both upper and lower hemicontinuous at $x$.
		\end{itemize}
	\end{definition}

\begin{assumption}\label{a:params}
	For each $t\in\{0, ..., T-1\}$,
	\begin{itemize}
		\item[(1)] $\gamma(dx_{t+1}, dy_{t+1}|\theta_{t+1}): (\Theta_{t+1}, d_{\Theta_{t+1}}) \rightarrow \cP_p(\cX_{t+1} \times \cY_{t+1})$ is continuous;
		\item[(2)] $D_t$ in \eqref{eq:paramD} is a continuous correspondence and $D_t(x_{1:t}, y_{1:t})$ is non-empty, convex, and compact, under the product topology $\cT_{\cX_{1:t}} \times \cT_{\cY_{1:t}} \times \cT_{\Theta_{t+1}}$;
		\item[(3)] for an arbitrary path $(x_{1:t}, y_{1:t}) \in \cX_{1:t} \times \cY_{1:t}$, we have $f(x_{1:t}, y_{1:t}, \theta_{t+1})$ in \eqref{eq:para_t_obj} strictly quasiconvex in $\theta_{t+1} \in \Theta_{t+1}$.
	\end{itemize}
\end{assumption}

We apply a version of Berge's maximum theorem with strict quasiconcavity, see \citet[Theorem 9.14 and Corollary 9.20]{sundaram1996first} or \citet[Theorem 17.31]{charalambos2013infinite}, to prove existence and uniqueness of the equilibrium parametric transport. The proof is given in Section \ref{sec:proof_params}.

\begin{theorem}\label{thm:params}
	Suppose Assumptions \ref{a:params_obj}, \ref{a:measure}, and \ref{a:params} hold. Then for each $t \in \{0, ..., T-1\}$,
	\begin{itemize}
		\item [(a)] there is a continuous and unique optimizer, $\theta^*_{t+1}( x_{1:t}, y_{1:t}): (\cX_{1:t} \times \cY_{1:t}, d) \rightarrow (\Theta_{t+1}, d_{\Theta_{t+1}})$, for the extended DP equation \eqref{eq:params_dpp};
		\item[(b)] the equilibrium parametric kernel $$\gamma(dx_{t+1:T}, dy_{t+1:T} | \theta^*_{t+1:T}(x_{1:t}, y_{1:t})): (\cX_{1:t} \times \cY_{1:t}, d) \rightarrow \cP_p(\cX_{t+1:T} \times \cY_{t+1:T})$$ in \eqref{eq:pi*params} is continuous and unique;
		\item[(c)] the corresponding value function $V_t$ satisfies \eqref{eq:params_dpp} and $V_t \in C_p(\cX_{1:t} \times \cY_{1:t})$.
	\end{itemize}
\end{theorem}
\begin{remark}
	$\theta^*_{t+1}$ is unique among all correspondences, including functions, with graphs that are subsets of the graph of $D_t$. If the objective $f(x_{1:t}, y_{1:t}, \theta_{t+1})$ satisfies the quasiconvexity but not strictly in Assumption \ref{a:params} (3), we could also include a regularization term on $\theta_{t+1}$ directly, such as $|\theta_{t+1}|^2$. If the regularization term is $\cT_{\Theta_{t+1}}$-continuous, then it can be treated similarly by Theorem \ref{thm:params}.
\end{remark} 

It is direct to recover the classic formulation when the correspondence $$D_t(x_{1:t}, y_{1:t}) = \Pi (\mu(d x_{t+1} |x_{1:t}), \nu(d y_{t+1}|y_{1:t}))$$ by setting $\Theta_{t+1} = \cP_p(\cX_{t+1} \times \cY_{t+1})$, $t \in \{0, ..., T-1\}$.

\begin{corollary}\label{cor:classical_nonMark}
	Suppose 
	\begin{itemize}
		\item[(1)] Assumptions \ref{a:params_obj} and \ref{a:measure} hold;
		\item[(2)] $\Theta_{t+1} = \cP_p(\cX_{t+1} \times \cY_{t+1})$, $t \in \{0, ..., T-1\}$;
		\item[(3)] for each given $t \in \{0, ..., T-1\}$ and an arbitrary path $(x_{1:t}, y_{1:t}) \in \cX_{1:t} \times \cY_{1:t}$, suppose $f(x_{1:t}, y_{1:t}, \gamma)$ in \eqref{eq:para_t_obj}  is strictly quasiconvex in any $$\gamma(d x_{t+1}, d y_{t+1}) \in \Pi (\mu(d x_{t+1} |x_{1:t}), \nu(d y_{t+1}|y_{1:t})).$$
	\end{itemize}
	Then for each $t \in \{0, ..., T-1\}$,
	\begin{itemize}
		\item [(a)] there is a continuous and unique optimizer, $\pi^*(dx_{t+1}, dy_{t+1} | x_{1:t}, y_{1:t}): (\cX_{1:t} \times \cY_{1:t}, d) \rightarrow \cP_p(\cX_{t+1} \times \cY_{t+1})$, for the extended DP equation \eqref{eq:params_dpp};
		\item[(b)] the equilibrium transport $\pi^*(dx_{t+1:T}, dy_{t+1:T} |x_{1:t}, y_{1:t}): (\cX_{1:t} \times \cY_{1:t}, d) \rightarrow \cP_p(\cX_{t+1:T} \times \cY_{t+1:T})$ in \eqref{eq:pi*params}, is continuous and unique;
		\item[(c)] the corresponding value function $V_t$ satisfies \eqref{eq:params_dpp} and $V_t \in C_p(\cX_{1:t} \times \cY_{1:t})$.
	\end{itemize}
\end{corollary}
\begin{remark}\label{rem:discrete}
	For the special case that $\mu$ and $\nu$ are discrete measures with finite supports, the extended DPP is straightforward. When $\cX$ and $\cY$ are equipped with discrete topology, the value function $V_t (x_{1:t}, y_{1:t})$ in \eqref{eq:params_dpp} and the correspondence $D_t$ in \eqref{eq:paramD} are continuous automatically. There exists an equilibrium transport $\pi^*$ which may not be unique. The strict quasiconvexity is not needed for the existence. However, we can regard probability masses as parameters and obtain the uniqueness if strict quasiconvexity holds. 
\end{remark}

\section{Illustrative examples}\label{sec:illustrate}
		\subsection{Dynamic matching under a mean-variance objective}\label{sec:mv}
	
	Inspired by \cite{hu2022dynamic}, we examine an example from operations management involving the matching of supply $x$ with demand $y$. Consider two supply types and two demand types, with type labels $\cX_t = \{0, 1\}$ and $\cY_t= \{0, 1\}$, respectively. In each period, supply and demand of different types arrive randomly: $x_{1:T} \sim \mu$ and $y_{1:T} \sim \nu$. The supply and demand types exhibit idiosyncratic taste towards one another, where matching closer types incurs lower costs. They are known as horizontally differentiated demand and supply types \cite[Section 4]{hu2022dynamic}. A practical example is a ride-hailing platform, where types are determined by the locations of riders and drivers, and closer matches result in reduced costs.
	
	The one-period cost matrix $f(x_t, y_t)$ is presented in Table \ref{tab:match}.
	\begin{table}
		\centering
		\begin{tabular}{ccc}
			\hline
			& $x = 0$ & $x = 1$ \\
			\hline
			$y=0$ & $1.0$ & $2.0$  \\
			$y=1$ & $2.0$ & $0.0$  \\
			\hline 
		\end{tabular}
		\caption{One-period cost matrix of supply/demand matching}\label{tab:match}
	\end{table}
	The total cost over $T$ periods is given by $c(x_{1:T}, y_{1:T}) = \sum^T_{t=1} \beta^t f(x_t, y_t)$, where $\beta \in [0, 1]$ acts as a discount factor. In the multi-period matching, bicausality is a natural requirement, meaning that no future information can be utilized. The agent aims to find a bicausal matching plan $\pi$ that minimizes both cost and variance:
	\begin{equation}\label{eq:mv}
		\inf_{\pi \in \Pi_{bc}(\mu, \nu) } \E_\pi \left[c(x_{1:T}, y_{1:T}) \right] + \gamma \text{Var}_\pi \left[c(x_{1:T}, y_{1:T}) \right].
	\end{equation} 
	Here, constant $\gamma > 0$ represents the tolerance level for variance.
	
	As an illustration, we consider a scenario where the number of periods $T = 2$, the discount factor $\beta = 1$, and the marginals $\mu$ and $\nu$ are given as follows: 
	\begin{equation}
		\begin{aligned}
			\mu(x_1) = 0.1 \delta_0 + 0.9 \delta_1, \quad \mu(x_2| x_1 = 0) & = 0.8 \delta_0 + 0.2 \delta_1, \quad \mu(x_2| x_1 = 1) = 0.2 \delta_0 + 0.8 \delta_1, \\
			\nu(y_1) = 0.5 \delta_0 + 0.5 \delta_1, \quad \nu(y_2 | y_1 = 0) & = 0.9 \delta_0 + 0.1 \delta_1, \quad \nu(y_2 | y_1 = 1) = 0.1 \delta_0 + 0.9 \delta_1.
		\end{aligned}
	\end{equation}
	Here, $\delta_i$ represents the Dirac measure at type $i$. The transition kernels are designed in a way that the next type is more likely to be the same as the previous one.
	
	The pre-committed solution to \eqref{eq:mv}, denoted as $\pi^{pre}$, is given by
	\begin{equation}\label{eq:pre}
		\begin{aligned}
			\pi^{pre}(x_1, y_1) &= 0.1 \delta_{(0, 0)} + 0.4 \delta_{(1, 0)} + 0.5 \delta_{(1, 1)}, \\
			\pi^{pre}(x_2, y_2| x_1 = 0, y_1 = 0) & = 0.8 \delta_{(0, 0)} + 0.1 \delta_{(1, 0)} + 0.1 \delta_{(1, 1)}, \\
			\pi^{pre}(x_2, y_2| x_1 = 1, y_1 = 0) & = 0.2 \delta_{(0, 0)} + 0.7 \delta_{(1, 0)} + 0.1 \delta_{(1, 1)}, \\
			\pi^{pre}(x_2, y_2| x_1 = 1, y_1 = 1) & = 0.2 \delta_{(0, 1)} + 0.1 \delta_{(1, 0)} + 0.7 \delta_{(1, 1)}.
		\end{aligned}
	\end{equation}
	Similarly, $\delta_{(i, j)}$ is the Dirac measure at $(x_2, y_2) = (i, j)$. The corresponding mean and variance of the cost are calculated as:
	\begin{equation}
		\E_{\pi^{pre}} \left[c(x_{1:2}, y_{1:2}) \right] = 1.94 \quad \text{and} \quad \text{Var}_{\pi^{pre}} \left[c(x_{1:2}, y_{1:2}) \right] = 2.6164.
	\end{equation}
	Consequently, the optimal value of the objective in \eqref{eq:mv} amounts to $4.5564$.
	
	In contrast, the equilibrium transport $\pi^*$ is given as follows:
	\begin{equation}\label{eq:equi}
		\begin{aligned}
			\pi^*(x_1, y_1) &= 0.1 \delta_{(0, 0)} + 0.4 \delta_{(1, 0)} + 0.5 \delta_{(1, 1)}, \\
			\pi^*(x_2, y_2| x_1 = 0, y_1 = 0) & = 0.8 \delta_{(0, 0)} + 0.1 \delta_{(1, 0)} + 0.1 \delta_{(1, 1)}, \\
			\pi^*(x_2, y_2| x_1 = 1, y_1 = 0) & = 0.1 \delta_{(0, 0)} + 0.1 \delta_{(0, 1)} + 0.8 \delta_{(1, 0)}, \\
			\pi^*(x_2, y_2| x_1 = 1, y_1 = 1) & = 0.1 \delta_{(0, 0)} + 0.1 \delta_{(0, 1)} + 0.8 \delta_{(1, 1)}.
		\end{aligned}
	\end{equation}
	The corresponding value of the objective in \eqref{eq:mv} is $5.0519$, with the mean and variance given by
	\begin{equation}
		\E_{\pi^*} \left[c(x_{1:2}, y_{1:2}) \right] = 1.91 \quad \text{and} \quad \text{Var}_{\pi^*} \left[c(x_{1:2}, y_{1:2}) \right] = 3.1419.
	\end{equation}
	
	We have the following comments on this example:
	\begin{enumerate}[label={(\arabic*)}]
		\item The pairs $(0, 1)$ and $(1, 0)$ are referred to as mismatches. Compared with the pre-committed OT,  \eqref{eq:equi} can increase or decrease the probabilities of mismatches in $(x_2, y_2)$, depending on the previous state $(x_1, y_1)$. Specifically, when $(x_1, y_1) = (1, 0)$, the total probability for $(x_2, y_2) = (0, 1)$ and $(x_2, y_2) = (1, 0)$ increases from $0.7$ to $0.9$. Conversely, when $(x_1, y_1) = (1, 1)$, this sum decreases from $0.3$ to $0.1$. Structural properties of equilibrium transport become more complex due to the dependence on the temporal structures of the marginals $\mu$ and $\nu$, a challenge not encountered in the static formulation \citep{boerma2023composite} or the infinite horizon setting \citep{shimer2000assortative}.
		
		\item Both the pre-committed OT \eqref{eq:pre} and the equilibrium transport \eqref{eq:equi} are not Monge maps. In the discrete setting, it is common to split probability masses, see also \citet[Section 3.1]{boerma2023composite}. Exploring the denseness of specific Monge maps for the equilibrium transport presents another intriguing and demanding task, extending the previous research by \cite{schrott2023denseness}.
		
		\item Regarding computational algorithms, the optimization problem \eqref{eq:mv} is not a linear program due to the presence of variance operator. While solvers like Gurobi can efficiently address this specific example, the general problem can be challenging because of the nonlinear $G$. In contrast to \cite{eckstein2024,gonzalez2024quad}, different methodologies are needed for algorithm design and convergence analysis. Despite this, identifying an equilibrium transport is easier than obtaining a pre-committed OT, given its local optimization nature and the possibility of solving sub-problems \eqref{eq:params_dpp} across states in parallel.
	\end{enumerate}

\subsection{Gaussian data}\label{sec:gaussian}
It is rare for continuous OT problems to have explicit solutions, even in the single-period case. One exception is the Gaussian distribution \citep{givens1984class,gunasingam2024adapted}. For simplicity, suppose $x_t \in \R$ and $y_t \in \R$ are one-dimensional. Consider $\mu$ and $\nu$ as Gaussian distributions with linear dynamics: 
\begin{equation}\label{eq:GaussianState}
	\begin{aligned}
		x_{t+1} & = x_t + \lambda_t, \quad \lambda_t \sim N(0, 1), \\
		y_{t+1} & = y_t + \eta_t, \quad \eta_t \sim N(0, 1).
	\end{aligned}
\end{equation}
The white Gaussian noise process $\{\lambda_t\}$ consists of standard normal random variables that are independent of each other. Impose the same condition on $\{\eta_t\}$.

\subsubsection{Normal distribution as parametric couplings}\label{sec:ex_params}
In this example for Theorem \ref{thm:params}, we consider a two-period problem with a nonlinear and state-dependent objective:
\begin{equation}\label{eq:normal}
	\left( x_0 y_0 - \int x_2 y_2 \pi(dx_{2}, dy_{2} | x_0, y_0) \right)^2, 
\end{equation}
with the state process given by \eqref{eq:GaussianState}. The initial states $x_0$ and $y_0$ are fixed as constants. 

We assume that the agent considers only normal distribution as couplings. At time $t=1$, the agent seeks an optimal coupling which is a bivariate normal distribution $\gamma(dx_2, dy_2| \theta_2)$ with marginals  $x_2 \sim N(x_1, 1)$ and $y_2 \sim N(y_1, 1)$ and $\theta_2 \in [-1, 1]$ as the correlation between $x_2$ and $y_2$. Since the first term in the square is state-dependent, the objective at time $1$ is
\begin{align*}
	\left( x_1 y_1 - \int x_2 y_2 \gamma(dx_2, dy_2| \theta_2) \right)^2.
\end{align*}
Recalling $f$ in \eqref{eq:para_t_obj}, we have $f(x_{1}, y_{1}, \theta_{2}) = (x_1y_1  - \theta_2 \times 1 \times 1 - x_1y_1)^2 = \theta^2_{2}$, which is strictly quasiconvex in $\theta_2 \in [-1, 1]$. For other conditions in Assumptions \ref{a:params_obj}, \ref{a:measure}, and \ref{a:params}, we take $1/r =2$ and $p=4$ as the growth rate. Assumptions \ref{a:params_obj} and \ref{a:measure} are satisfied. $\gamma(dx_2, dy_2| \theta_2)$ is continuous in $\theta_2$ with the Wasserstein distance of order 4, by properties of the normal distribution and \citet[Theorem 6.8]{villani2009optimal}. Clearly, the correspondence $D_1: (x_1, y_1) \twoheadrightarrow [-1, 1]$ satisfies Assumption \ref{a:params} (2). Theorem \ref{thm:params} shows that $\theta^*_2(x_1, y_1) = 0$ is the unique optimizer at time $1$. Similarly, $\theta^*_1(x_0, y_0) = 0$, which is also unique. In this case, the restriction on normal couplings is not too restrictive since it gives a value function equal to zero.

Consider the agent who ignores time inconsistency and minimizes the objective \eqref{eq:normal} at time 0 only. Any constant $(\theta_1, \theta_2)$ with $\theta_2 + \theta_1 = 0$ is an optimizer for $ \left( x_0 y_0 - \int x_2 y_2 \pi(dx_{2}, dy_{2} | x_0, y_0) \right)^2 = (\theta_2 + \theta_1)^2$. However, for nonzero $\theta_1$, the agent at time $1$ will find that it is optimal to deviate from $\theta_2 = - \theta_1$. Therefore, it is possible to have a unique equilibrium parametric transport but multiple globally optimal parametric transports.  


\subsubsection{State dependence with alternating signs} To provide another explicit example, we set $\varphi(k) = (-1)^{k+1}$ and consider the following artificial objective, which may lack any economic motivation: 
\begin{equation}\label{eq:alt}
	\inf_{\pi \in \Pi_{bc}(\mu(dx_{t+1:T}|x_{1:t}), \nu(dy_{t+1:T}|y_{1:t}))} \int  \sum^T_{s = t+1} \varphi(s - t) (x_s - y_s)^2 \pi(dx_{t+1:T}, dy_{t+1:T} | x_{1:t}, y_{1:t}),
\end{equation}
where the state process is still given by \eqref{eq:GaussianState}. Since $\varphi(\cdot)$ can be negative, it is no longer interpreted as a discounting function, as in \eqref{eq:nonexp}. The term $\varphi(s - t)$ reflects the state dependence on time $t$.
	
Unlike the previous example, we allow all bicausal couplings. Suppose $T = 2$. A pre-committed solution is any optimizer of the following problem at time $t=0$: 
\begin{equation}\label{eq:pre_alt}
	\inf_{\pi \in \Pi_{bc}(\mu, \nu)} \int  \left[ (-1)^{1-0+1}(x_1 - y_1)^2 + (-1)^{2-0+1} (x_2 - y_2)^2 \right] \pi(dx_{1:2}, dy_{1:2}).
\end{equation}
After expressing $\pi$ in terms of successive regular kernels, the problem in \eqref{eq:pre_alt} can be solved as follows:
	\begin{align*}
		& \inf_{\bar{\pi} \in \Pi(p^1_*\mu, p^1_*\nu)} \int  \Big[ (x_1 - y_1)^2 - \sup_{\gamma \in \Pi(\mu(dx_2|x_1), \nu(dy_2|y_1))} \int (x_2 - y_2)^2 \gamma(dx_{2}, dy_{2}) \Big]\bar{\pi}(dx_{1}, dy_{1}) \\
		& = \inf_{\bar{\pi} \in \Pi(p^1_*\mu, p^1_*\nu)} \int  \left[ (x_1 - y_1)^2 - 4 - (x_1 - y_1)^2 \right]\bar{\pi}(dx_{1}, dy_{1}) \\
		& = - 4.
	\end{align*}  
The inner maximization is derived from \citet[Lemma 1]{han2025IEEE}. For $t=1$, an optimizer in $\Pi(\mu(dx_2|x_1), \nu(dy_2|y_1))$ is the coupling under which $X_2$ and $Y_2$ follow the bivariate normal distribution with correlation $-1$ and the given marginals. Any coupling $\bar{\pi} \in \Pi(p^1_*\mu, p^1_*\nu)$ minimizes the objective at $t=0$.

By Definition \ref{def:equitrans}, to find an equilibrium transport, we first solve the problem at time $t=1$, given by
	\begin{align*}
		\inf_{\gamma \in \Pi(\mu(dx_2|x_1), \nu(dy_2|y_1))} \int \varphi(2-1) (x_2 - y_2)^2 \gamma(dx_{2}, dy_{2}) = (x_1 - y_1)^2.
	\end{align*}
An optimizer of this problem, denoted as $\pi^*(dx_2, dy_2 | x_1, y_1)$, is the coupling under which $X_2$ and $Y_2$ follow the bivariate normal distribution with correlation $1$ and the given marginals. Then the problem at time $t=0$ becomes
	\begin{align*}
		& \inf_{\bar{\pi} \in \Pi(p^1_*\mu, p^1_*\nu)} \int  \Big[ \varphi(1-0)(x_1 - y_1)^2 + \varphi(2-0) \int (x_2 - y_2)^2 \pi^*(dx_2, dy_2 | x_1, y_1) \Big]\bar{\pi}(dx_{1}, dy_{1}) \\
		& = \inf_{\bar{\pi} \in \Pi(p^1_*\mu, p^1_*\nu)} \int  \left[ (x_1 - y_1)^2 - (x_1 - y_1)^2 \right]\bar{\pi}(dx_{1}, dy_{1}) \\
		& = 0,
	\end{align*}  
where any coupling $\bar{\pi} \in \Pi(p^1_*\mu, p^1_*\nu)$ is a minimizer.
	
The state dependence and alternating signs of $\varphi(\cdot)$ result in equilibrium transport that differs from the pre-committed solutions. This example demonstrates the subgame formulation of equilibrium transport.

\section{Roberts' law: Literature and methodologies}\label{sec:Roberts}	
	This section briefly reviews empirical findings in labor markets and relevant literature, providing motivation for the models developed in the following sections.
	
\begin{figure}
		\centering
		\includegraphics[width=0.4\textwidth]{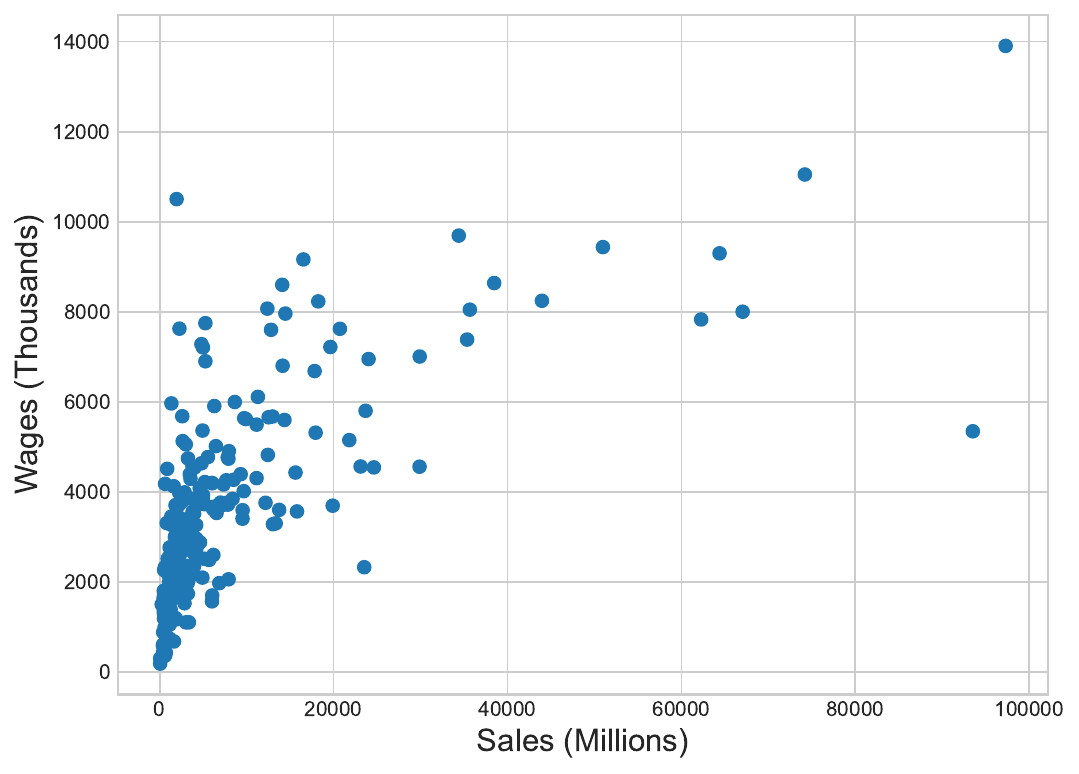}
		\caption{Sale-wage relationship in the Industrials sector, fiscal year 2021. \label{fig:roberts}}
	\end{figure}

	Roberts' law \cite[Section II]{gabaix2008has} states that CEO compensation is proportional to $(\text{own firm size})^{\kappa}$, with a typical empirical exponent $\kappa \simeq 1/3$. However, data from executive labor markets do not always show a perfect correlation between firm size and wages. Figure \ref{fig:roberts} plots wage and net sales data for the Industrials sector in the 2021 fiscal year. A power-law relationship holds approximately, but some large firms deviate from Roberts' law and underpay their managers.
	
	\cite{bonhomme2019} summarize two influential approaches to studying wage dispersion and worker-firm sorting:
	\begin{itemize}
		\item \cite{abowd1999high} (hereafter AKM) propose a decomposition of wages into worker and firm components, along with an error term. Specifically, suppose worker $i$ is employed at firm $J(i, t)$ at time $t$ and earns a wage $w_{i, J(i, t)}$. The AKM specification assumes:
		\begin{equation}\label{AKM} 
			\log[w_{i, J(i, t)}] = r_i + \psi_{J(i, t)} + \varepsilon_{it},
		\end{equation}
		where $r_i$ is the worker fixed effect, $\psi_{J(i, t)}$ is the firm fixed effect, and $\varepsilon_{it}$ is a mean-zero residual. Estimating \eqref{AKM} requires matched employer-employee data tracking both workers and firms over time.
		
		The AKM approach is a type of reduced-form, or non-structural, analysis. This econometric method models the relationship between a dependent variable and explanatory variables without specifying the underlying economic mechanisms. Its strength lies in tractability, particularly in dealing with unobserved heterogeneity. However, it may rely on restrictive assumptions, such as the exclusion of interactions between worker and firm attributes. Related studies include \cite{card2013workplace,song2019firming,borovickova2020high}, among others.
		
		\item The second approach follows a structural methodology, which specifies a theoretical model of sorting to explain observed labor market outcomes. Notable examples include \cite{becker1973theory,shimer2000assortative,postel2002equilibrium,gabaix2008has,hagedorn2017identifying}. Structural models aim to identify the mechanisms driving outcomes but are typically harder to estimate. Concerns often arise regarding the extent to which results depend on functional form assumptions.
	\end{itemize}

	In the following two sections, we propose two models to examine state dependence in job markets from different perspectives. Inspired by the reduced-form methods, Section \ref{sec:reduced_jobm} directly analyzes the relationship between sales and wages by approximating observed pairs with equilibrium transport plans. In contrast, Section \ref{sec:struct_jobm} takes a structural approach, following the methodologies of \cite{gabaix2008has,galichon2016optimal}.

\section{State dependence between sales/rankings and wages}\label{sec:reduced_jobm}
	
\subsection{Formulation}
	
	Consider a time horizon of $T$ years. The agent observes two time series, $x = x_{1:T} \sim \mu$ and $y = y_{1:T} \sim \nu$, where $x_t$ and $y_t$ take values in the same range. In the context of the CEO job market discussed in Section \ref{sec:exec}, $x_t$ denotes the rank of net sales and $y_t$ the rank of wages at time $t$, both after clustering. We choose to compare sales and wages directly, as both are observable and their transition matrices are easier to estimate. Similarly, \cite{borovickova2020high} define a worker's type as the average log wage received over a long horizon.
	
	The goal is to quantify the statistical association between $x$ and $y$. In practice, their joint (empirical) distribution $(x, y) \sim \pi_r$ is typically available. For instance, the empirical data record wages paid by firms. A standard approach is to compute Spearman or Kendall rank correlations between $x_t$ and $y_t$ at each time $t$. Alternatively, one may adopt the AKM methodology and perform regression analysis as in \eqref{AKM}. However, both approaches are static, in the sense that they do not capture the dependence of the future pair $(x_{t+1}, y_{t+1})$ on the current pair $(x_t, y_t)$.
	
	
	Bicausal OT in \eqref{eq:bc} can also quantify the discrepancy between $\mu$ and $\nu$ for two time series. From this perspective, the agent observes the transport plan $\pi_r$, but the cost function $c(x, y)$ is unknown. A similar question is considered in the OT literature by \cite{stuart2020inverse}, which introduces inverse OT to infer unknown cost functions from noisy observations of OT plans. These problems are generally ill-posed, necessitating suitable simplifications.
	
	In the bicausal OT framework, the joint distribution with marginals $\mu$ and $\nu$ is determined as an optimizer of \eqref{eq:bc}. With an appropriately chosen cost function, the optimizer should be close to the observed plan $\pi_r$. We begin with the simple case in which Roberts' law holds exactly: the firm ranked $n$ pays wages ranked $n$. In this case, the $L_p$ norm is a natural choice for the cost function, making $\pi_r$ an optimizer. This choice is not unique, and alternatives such as $-x_t^a y_t^b$ may also be considered. Notably, the specific case $-x_t y_t$ yields the same primal solution as the $L_2$ norm $(x_t - y_t)^2$.
	
	When $y_t$ represents wage ranks, the dual problem does not admit a wage-equation interpretation as in \cite{boerma2023composite}. In this section, we focus on the statistical association between sales and wages, rather than interpreting the problem as a {\it decentralized} matching process between risk-neutral firms and managers. This latter interpretation from the dual perspective motivates the model in Section \ref{sec:struct_jobm}.

	For simplicity, we adopt the $L_1$ norm and consider \eqref{eq:PAM} when Roberts' law holds exactly:
	\begin{align}\label{eq:PAM}
		\inf_{\pi \in \Pi_{bc}(\mu, \nu)} \int  \sum^T_{t=1} \beta^t |x_t - y_t| \, \pi(dx_{1:T}, dy_{1:T}),
	\end{align}
	where $\beta \in [0, 1]$ is a discount factor.
	
	When empirical data deviate from Roberts' law, as shown in Figure \ref{fig:roberts}, can the cost in \eqref{eq:PAM} be modified to produce a solution closer to the observed plan $\pi_r$? We conjecture that there exist incentives to preserve the previous pair $(x_t, y_t)$, even if it does not align with Roberts' law. To capture this, we introduce a new term that reduces the cost if we keep similar matching:
	\begin{align}\label{eq:incon_match}
		\inf_{\pi \in \Pi_{bc}(\mu^t, \nu^t)} \int \left[ - \alpha e^{- \frac{|x_{t+1} - x_t| + |y_{t+1} - y_t|}{\tau}} + \sum^T_{s=t+1} \beta^{s-t} |x_s - y_s| \right] \pi(dx_{t+1:T}, dy_{t+1:T} \mid x_{1:t}, y_{1:t}).
	\end{align}
	Here, we denote $\mu^t := \mu(dx_{t+1:T} \mid x_{1:t})$ and $\nu^t := \nu(dy_{t+1:T} \mid y_{1:t})$. The parameter $\alpha$ is a constant to be calibrated, and $\tau > 0$ is a scaling factor. We refer to the first term as the state-dependent preference function, which decays rapidly when $(x_{t+1}, y_{t+1})$ deviates from $(x_t, y_t)$. Interpreting $0/0 = 0$ and setting $\tau = 0$, this function reduces to the indicator function $\mathbf{1}_{\{ x_{t+1} = x_t,\, y_{t+1} = y_t \}}$ as a special case.

	The motivation for the state-dependent term can be further explained as follows:
	\begin{enumerate}[label={(\arabic*)}]
		\item From a regulatory perspective, public companies are generally required to disclose executive compensation for the past {\it three} fiscal years [\footnote{\url{https://www.sec.gov/answers/execcomp.htm}}]. This disclosure helps stakeholders assess year-over-year changes in compensation.
		
		\item From a financial perspective, executive compensation typically consists of a base salary, cash incentives, and long-term equity awards. While the base salary tends to be stable, incentives and equity awards are often tied to stock performance, which may exhibit momentum and depend on past outcomes.
		
		\item From a psychological and behavioral perspective, the specification aligns with the \textit{status quo bias} when $\alpha > 0$, reflecting a preference for maintaining the current state. Employers may preserve compensation rankings among industry peers, and employees may remain in underpaying positions due to inertia or relocation costs.
	\end{enumerate}
	
	Overall, a state-dependent term with $\alpha > 0$ introduces persistence or inertia, favoring matches close to the previous pair $(x_t, y_t)$. It is important to note that the state dependence here is under sales (or university rankings) $x_t$ and wages $y_t$, while Section \ref{sec:struct_jobm} focuses on the state dependence under firm size and worker talent.
	
	The problem \eqref{eq:incon_match} becomes state-dependent and therefore time-inconsistent. Since $x_t$ and $y_t$ represent ranks, the extended DP equation in the discrete case suffices. Let $\pi(\alpha)$ denote the equilibrium transport plan corresponding to a given $\alpha$. We calculate the classic Wasserstein distance $\cW(\pi(\alpha), \pi_r )$ between $\pi(\alpha)$ and the observed actual transport plan (matching) $\pi_r$:
	\begin{equation}\label{eq:wass_alpha}
		\cW(\pi(\alpha), \pi_r ) := \inf_{\gamma \in \Pi(\pi(\alpha), \pi_r)} \int \sum^T_{t=1} (|x_t - x'_t| + |y_t - y'_t| ) \, d\gamma.
	\end{equation}
	By varying $\alpha$, we identify the value that minimizes $\cW(\pi(\alpha), \pi_r)$. A negative optimal $\alpha$ suggests that the data favor deviation from the previous matching, whereas a positive value indicates a tendency to preserve it.

	Compared with the literature discussed in Section \ref{sec:Roberts}, reduced-form methods examine the association using regression models such as \eqref{AKM}, whereas we approximate the empirical plan $\pi_r$ with a distribution $\pi(\alpha)$ parameterized by $\alpha$.
	
	In the next two subsections, we investigate the inertia, measured by state dependence, of two job markets: top-ranking executives and academia (faculty and postdocs). These markets are selected primarily due to data availability. For simplicity and tractability, we impose Assumption \ref{a:Markov}. These simplifications address data limitations and improve calibration robustness, particularly in the academic job market setting.
	\begin{assumption}\label{a:Markov}
		The conditional kernels $\mu(dx_{t+1}|x_{1:t})$ and $\nu(dy_{t+1}|y_{1:t})$ are Markov and time-homogeneous. That is, $\mu(dx_{t+1}|x_{1:t}) = \mu(dx_{t+1}|x_t)$ and is the same for any $t = 1, ..., T-1$. The same condition holds for $\nu(dy_{t+1}|y_{1:t})$ as well.
	\end{assumption}

	\subsection{Executive job market}\label{sec:exec}
	
	Consider a five-year time horizon from 2017 to 2021. Firms are categorized into industries based on the Global Industry Classification Standard (GICS). To identify which sectors align more closely with Roberts' law, we compute the Spearman and Kendall rank correlations between net sales and compensation, as shown in Table \ref{tab:corr}. The Consumer Discretionary (GICS Code 25), Real Estate (GICS Code 60), and Information Technology (GICS Code 45) sectors exhibit the weakest correlations between net sales and wages. In contrast, the Utilities sector (GICS Code 55) shows the strongest correlation. For clarity, we define an efficient job market as follows.
	
	\begin{definition}
		A job market is said to be more efficient if the sale-wage correlation is higher. We refer to this correlation as the efficiency coefficient.
	\end{definition}
	
	As a validation, we first apply model \eqref{eq:incon_match} to perfectly matched synthetic data, where the firm ranked $n$ always pays the wage ranked $n$. The first row of Table \ref{tab:6} shows that the corresponding optimal benchmark $\alpha$ is close to zero across most sectors. We then apply the model to bootstrap samples from the real data. The row labeled ``raw $\alpha$" in Table \ref{tab:6} shows that the optimal $\alpha$ shifts in a positive direction. Sector 50 yields a smaller value due to minor fluctuations in the Wasserstein distance. The final row of Table \ref{tab:6} reports the difference between the raw and benchmark $\alpha$. Overall, the results suggest the presence of state dependence or inertia in most industries, regardless of whether the raw or adjusted $\alpha$ is considered.

	\begin{table}
		\centering
		\begin{tabular}{cccccccccccc}
			\hline
			Sector & 10 & 15 & 20 & 25 & 30 & 35 & 40 & 45 & 50 & 55 & 60 \\
			\hline
			Spearman & 0.72 & 0.793 & 0.752 & 0.515 & 0.755 & 0.768 & 0.745 & 0.646 & 0.75 &
			0.871 & 0.559 \\
			\hline
			Kendall & 0.543 & 0.598 & 0.564 & 0.362 & 0.571 & 0.587 & 0.579 & 0.475 & 0.564 & 0.692 & 0.388 \\
			\hline
		\end{tabular}
		\caption{Spearman and Kendall rank correlation between net sales and compensations in five years.} \label{tab:corr}
	\end{table}

	\begin{table}
		\centering
		{\scriptsize
			\begin{tabular}{cccccccccccc}
				\hline
				Sector & 10 & 15 & 20 & 25 & 30 & 35 & 40 & 45 & 50 & 55 & 60 \\
				\hline
				Benchmark $\alpha$ & $-0.084$ & $-0.06$  & $-0.06$ & $-0.084$ & $-0.042$ & $-0.114$ & $-0.288$ & $-0.096$ & $-0.258$ & $-0.12$ & $-0.06$ \\
				Raw $\alpha$ &  0.438 & 0.006 &  0.276 &  0.522 &  0.054 & 0.03 & $-0.048$ &  0.66 & $-0.384$ &  0.018 &  0.702 \\
				Adjusted $\alpha$ & 0.522 &  0.066 &  0.336 &  0.606 &  0.096 &  0.144 &  0.24 &  0.756 & $-0.126$ &  0.138 &  0.762 \\
				\hline 
			\end{tabular}
		}
		\caption{Mean values of the optimal $\alpha$ in ten simulations. The number of clusters is set as 6. If there are multiple optimal $\alpha$, we choose the one that is closest to zero.}\label{tab:6}
	\end{table}

\begin{figure}
		\centering
			\begin{minipage}{0.45\textwidth}
				\centering
				\includegraphics[width=0.95\textwidth]{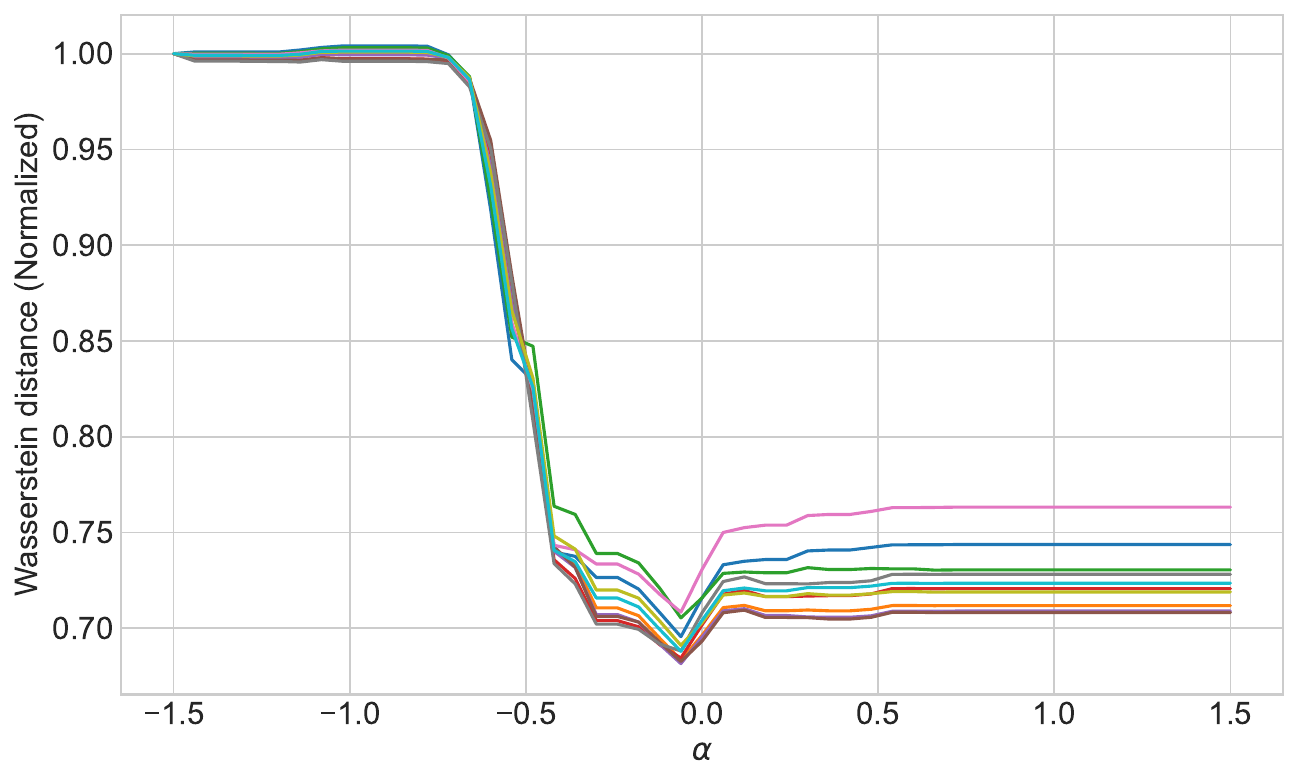}
				\subcaption{Synthetic data with perfect matching}\label{fig:MT_zero}
			\end{minipage}
			\begin{minipage}{0.45\textwidth}
				\centering
				\includegraphics[width=0.95\textwidth]{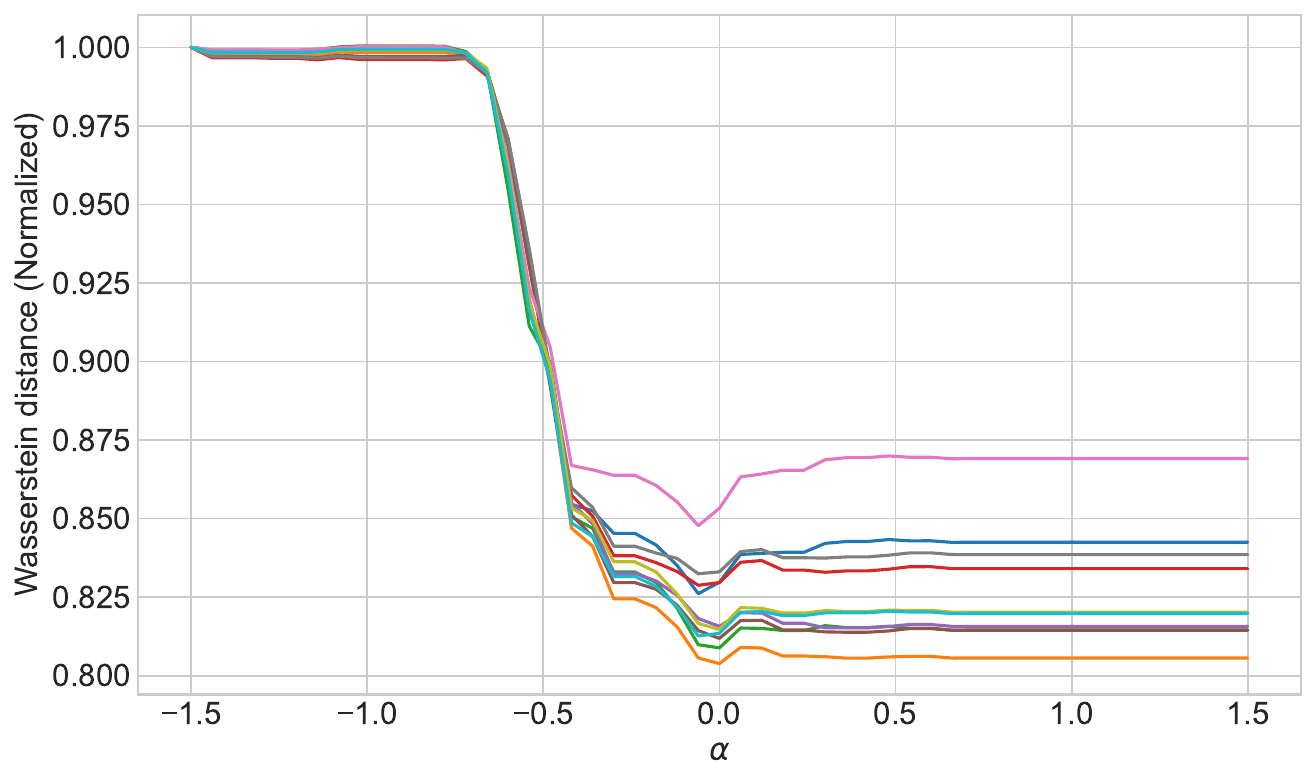}
				\subcaption{Resampled real data}\label{fig:MT_real}
			\end{minipage}%
		\caption{Calibration curves for Materials sector (GICS Code 15). Ten curves in each subplot represent ten independent simulations. \label{fig:MTRL}}
	\end{figure}

\begin{figure}
		\centering
			\begin{minipage}{0.45\textwidth}
				\centering
				\includegraphics[width=0.95\textwidth]{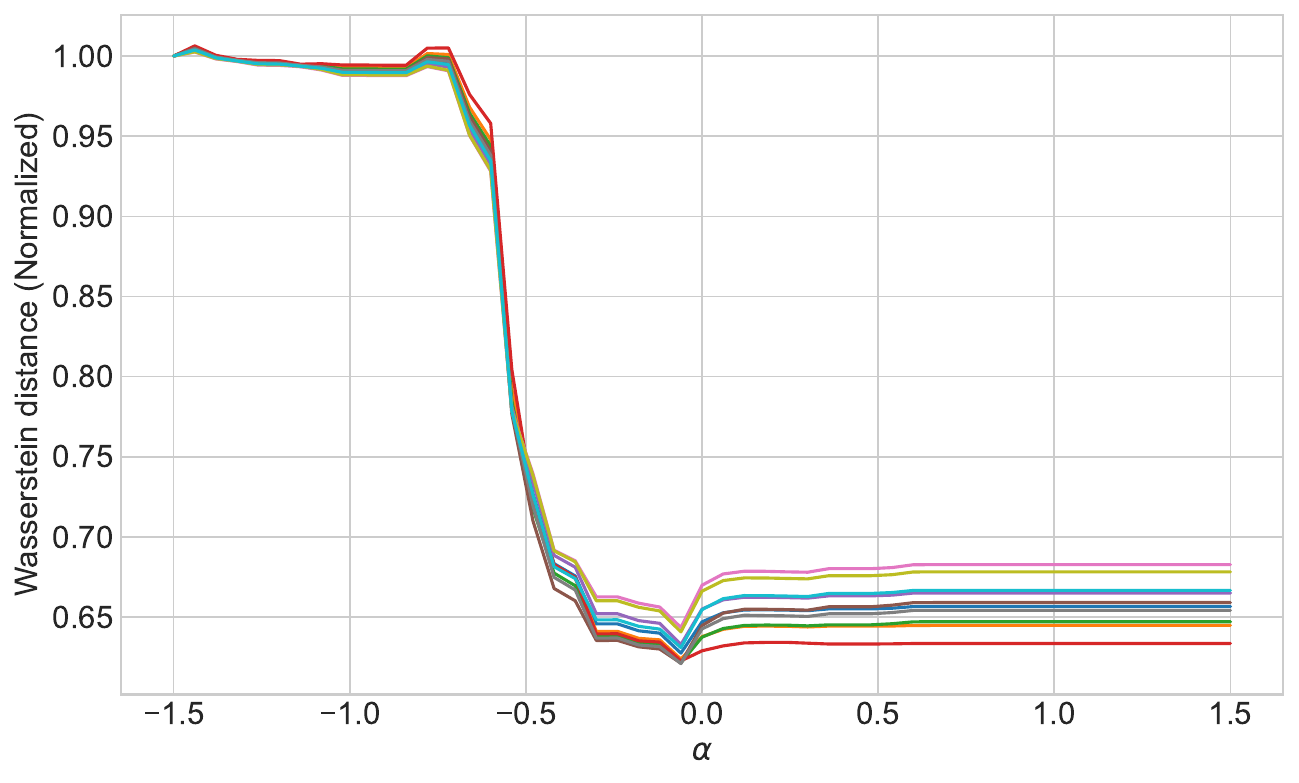}
				\subcaption{Synthetic data with perfect matching}\label{fig:RE_zero}
			\end{minipage}
			\begin{minipage}{0.45\textwidth}
				\centering
				\includegraphics[width=0.95\textwidth]{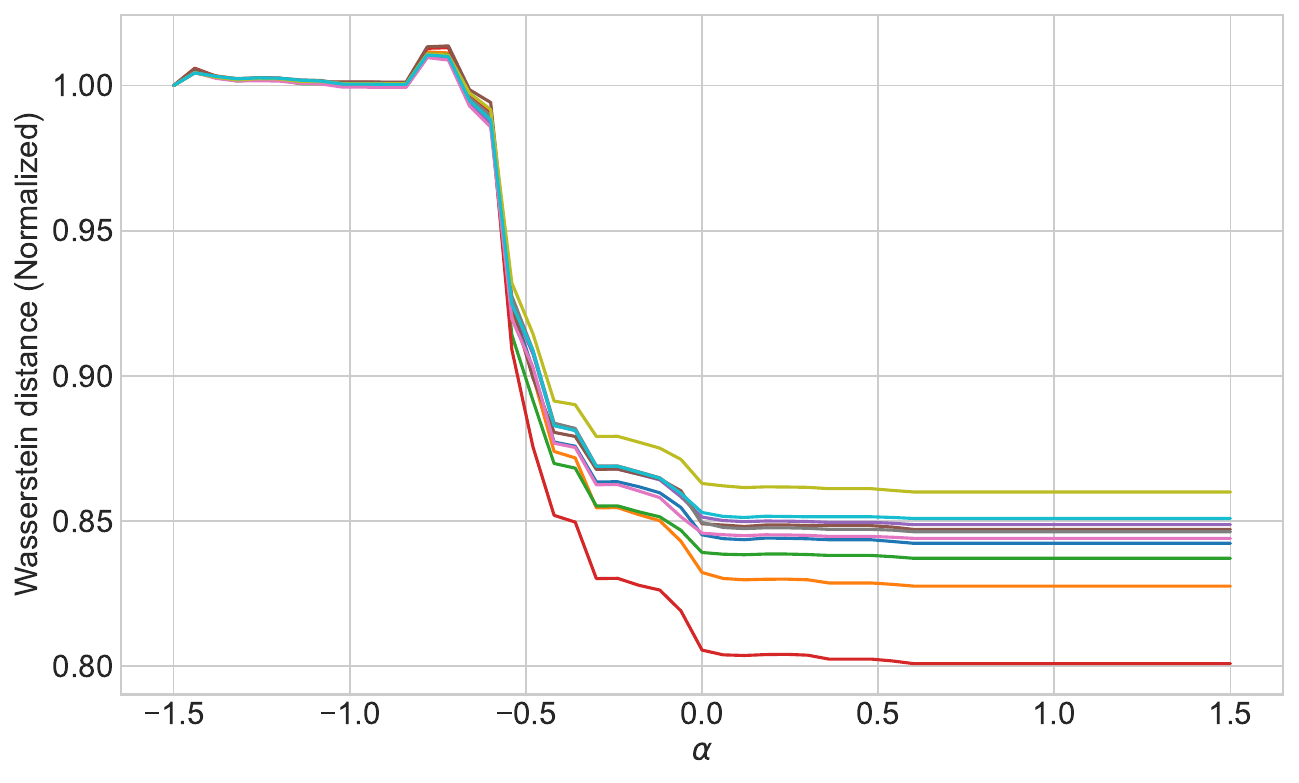}
				\subcaption{Resampled real data}\label{fig:RE_real}
			\end{minipage}%
		\caption{Calibration curves for Real Estate sector (GICS Code 60). Ten curves in each subplot represent ten independent simulations. \label{fig:RLEST}}
	\end{figure}

	For the Materials and Real Estate sectors, Figures \ref{fig:MTRL} and \ref{fig:RLEST} present the calibration curves, where the vertical axis represents the normalized Wasserstein distance $\cW(\pi(\alpha), \pi_r)$ from \eqref{eq:wass_alpha}. The distance is normalized by $\cW(\pi(-1.5), \pi_r)$ so that all curves begin at $1.0$ when $\alpha = -1.5$. Each subplot displays ten curves, corresponding to ten independent bootstrap samples. Comparing Figures \ref{fig:MT_zero} and \ref{fig:MT_real}, we observe that the shape of the calibration curves remains largely unchanged when perfectly matched data are replaced with resampled real data. This indicates that the Materials sector does not exhibit significant state dependence or inertia in sale-wage pairs. In contrast, the Real Estate sector shows inertia: the optimal $\alpha$ shifts toward positive values in Figure \ref{fig:RE_real}, in contrast to the synthetic perfectly matched data in Figure \ref{fig:RE_zero}.

	\subsubsection*{Main observation}
	A negative relationship emerges between job market efficiency and inertia. This is supported by the negative correlation between the optimal $\alpha$ values in Table \ref{tab:6} and the sale-wage efficiency coefficients in Table \ref{tab:corr}. Table \ref{tab:exec_test} further shows that these correlation values are generally below $-0.6$ and statistically significant at the 5\% level.
	
	\begin{table}
		\centering
		\begin{tabular}{ccc}
			\hline
			Correlation & Spearman ($p$-value) & Kendall ($p$-value) \\
			\hline 
			Raw $\alpha$ & $-0.645$ $(0.032)$ & $-0.600$ $(0.010)$ \\
			\hline
			Adjusted $\alpha$ & $- 0.773$ $(0.005)$ & $-0.600$ $(0.010)$ \\
			\hline
		\end{tabular}
		\caption{The relation between job market efficiency and the inertia effect. The correlations are between the optimal $\alpha$ in Table \ref{tab:6} and the sale-wage efficiency coefficient in Table \ref{tab:corr}. The number of clusters is 6.}\label{tab:exec_test}
	\end{table}

	\subsection{Academic job market}\label{sec:prof}
	
	Compared to the executive job market, the academic job market has received less attention in the literature. Our analysis fills this gap using University of California (UC) compensation data from 2017 to 2021.
	
	Table \ref{tab:uc_corr} reports Spearman and Kendall rank correlations between university rankings and wages. Professor-level positions are drawn from Business, Economics, and Engineering (B/E/E) departments, while postdocs are from all departments. Salaries for B/E/E assistant professors align more closely with university rankings. In contrast, postdoc salaries exhibit the weakest correlation with university rankings among all job positions and business sectors, motivating an investigation into whether postdoc wages exhibit stronger inertia.
	
	\begin{table}[H]
		\centering
		\begin{tabular}{ccccc}
			\hline
			Position & Professor & Associate Professor & Assistant Professor & Postdoc \\
			\hline
			Spearman    & 0.78 &  0.747 &  0.868 &  0.403 \\
			Kendall & 0.616 & 0.585 & 0.707 &  0.298 \\
			\hline 
		\end{tabular}
		\caption{Correlations between university rankings and wages in 2017--2021.}\label{tab:uc_corr}
	\end{table}
	
	\begin{figure}[H]
		\centering
		\begin{minipage}{0.48\textwidth}
			\centering
			\includegraphics[width=0.9\textwidth]{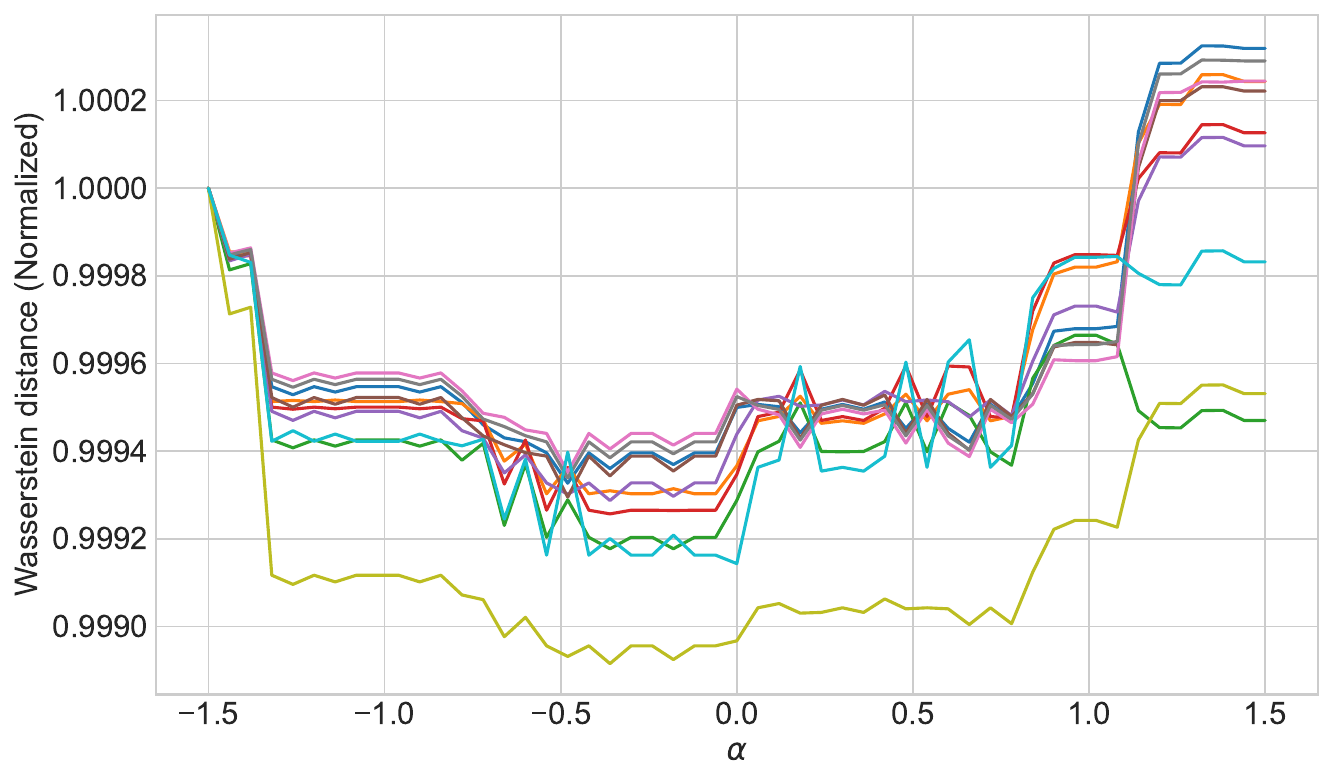}
			\subcaption{Synthetic data with perfect matching}\label{fig:postdoc_zero}
		\end{minipage}
		\begin{minipage}{0.48\textwidth}
			\centering
			\includegraphics[width=0.9\textwidth]{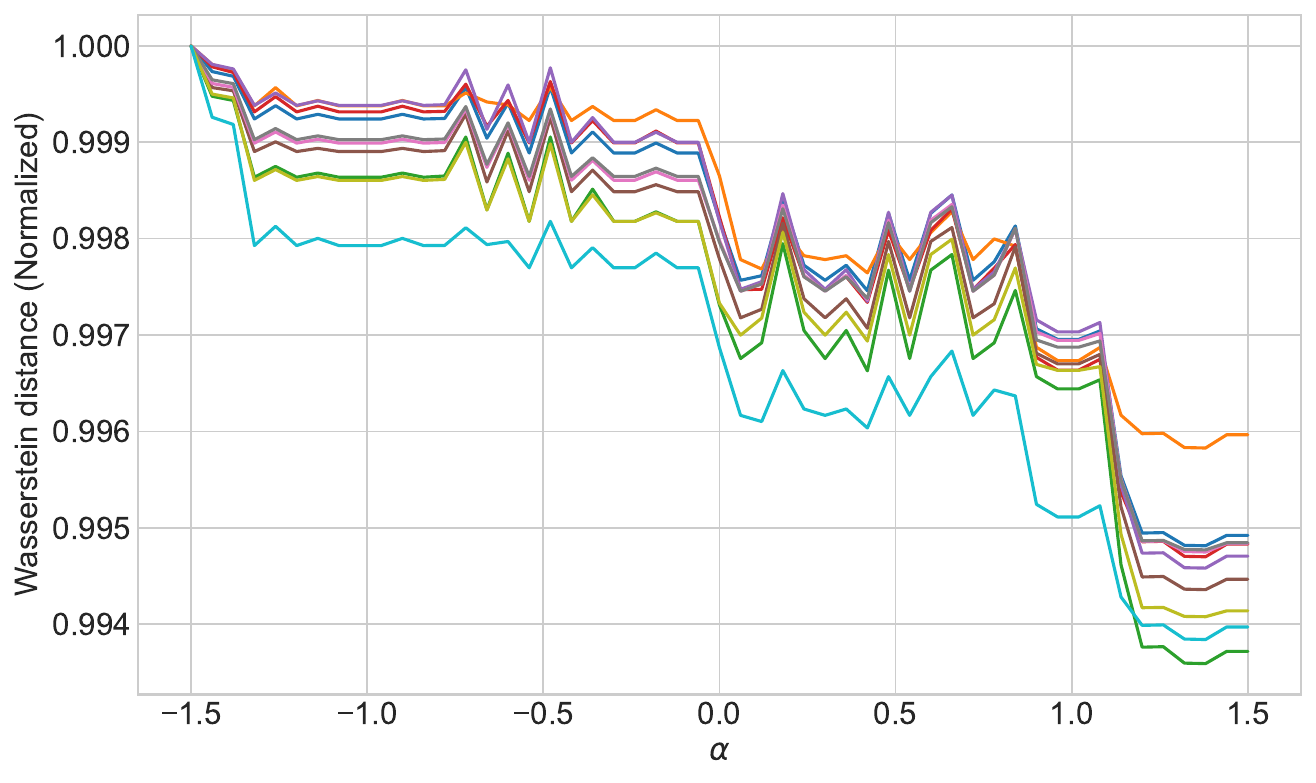}
			\subcaption{Resampled real data}\label{fig:postdoc_real}
		\end{minipage}%
		\caption{Calibration curves for postdocs. \label{fig:postdoc}}
	\end{figure}
	
	As a validation, Figure \ref{fig:postdoc_zero} shows that the optimal $\alpha$ is near zero when using synthetic perfectly matched data. However, the trend shifts notably when bootstrap samples from the postdoc data are used. As shown in Figure \ref{fig:postdoc_real}, the optimal $\alpha$ exceeds $1.0$, indicating strong inertia. This trend is also more pronounced than in the executive data. In contrast, the calibration curves for the three professor-level positions do not indicate significant inertia. Further details are provided in the e-companion.

\section{State dependence between sizes and talents: A structural approach}\label{sec:struct_jobm}
	
	\subsection{Formulation}
	In this section, we follow \cite{gabaix2008has} and introduce talents as latent variables for workers. Unlike Section \ref{sec:reduced_jobm}, the structural model developed here is primarily motivated by a dual formulation; see \citet[Chapter 4]{galichon2016optimal}.
	
	Consider an economy with an equal number of firms and workers. A firm of type $n$ has size $S(n)$, and a worker of type $m$ possesses talent $Q(m)$. A lower value of $n$ corresponds to a larger firm, while a lower $m$ indicates a more talented worker. Note that $n$ and $m$ are not necessarily ranks. For simplicity, assume a finite number of types for both firms and workers. Over a time horizon of length $T$, the marginal distributions of types are given by $(n_1, \ldots, n_t, \ldots, n_T) \sim \mu$ and $(m_1, \ldots, m_t, \ldots, m_T) \sim \nu$, where $n_t$ and $m_t$ denote the firm and worker types at time $t$, respectively. Suppose Assumption \ref{a:Markov} holds for both $\mu$ and $\nu$.
	
	
	Following \cite{gabaix2008has}, we assume that the earnings generated when a worker of type $m$ is matched with a firm of type $n$ are given by $C_1 S^a(n) Q(m)$, where $C_1 > 0$ is a constant. In line with the minimization framework, we consider the cost as the negative of earnings: $c(n, m) = - C_1 S^a(n) Q(m)$.
	
	At each time $t = 1, \ldots, T-1$, a central planner minimizes the following cost functional, which incorporates state dependence:
	\begin{align}\label{eq:Cobb}
		\inf_{\pi \in \Pi_{bc}(\mu^t, \nu^t)} \int \Big[ - \alpha e^{- \frac{(n_{t+1} - n_t)^2 + (m_{t+1} - m_t)^2}{\tau}} + \sum^T_{k=t+1} \beta^{k - t} c(n_k, m_k) \Big] \pi(dn_{t+1:T}, d m_{t+1:T} |n_{t}, m_{t}).
	\end{align}
	
	Compared with \eqref{eq:incon_match}, we adopt the quadratic form $(n_{t+1} - n_t)^2 + (m_{t+1} - m_t)^2$ in \eqref{eq:Cobb} to facilitate the computation of derivatives. The model without the state-dependent term has been considered in \cite{gabaix2008has,galichon2016optimal}. 
	
	For later use, given an equilibrium transport plan $\pi^*(dn_{t+2:T}, dm_{t+2:T} | n_{t+1}, m_{t+1})$, define the expected cost excluding the state-dependent term as
	\begin{align*}
		L_{t+1}(n_{t+1}, m_{t+1}) := \int  \sum^T_{k=t+2} \beta^{k - (t+1)} c(n_k, m_k) \pi^*(dn_{t+2:T}, dm_{t+2:T} | n_{t+1}, m_{t+1}).
	\end{align*} 
	With this definition, the equilibrium transport problem \eqref{eq:semi_dpp_simple} at time $t$ becomes
	\begin{align}
		\inf_{\gamma \in \Pi(\mu(dn_{t+1}|n_t), \nu(dm_{t+1}|m_t))} \int \Big[ & - \alpha e^{- \frac{(n_{t+1} - n_t)^2 + (m_{t+1} - m_t)^2}{\tau}} \label{eq:CobbEqui} \\
		& + \beta c(n_{t+1}, m_{t+1}) + \beta L_{t+1}(n_{t+1}, m_{t+1}) \Big] \gamma(dn_{t+1}, dm_{t+1}). \nonumber
	\end{align}
	The effect of state dependence appears in two ways. It contributes directly to the cost via the first term, and indirectly by influencing the equilibrium transport plan $\pi^*$ after $t+1$, thereby affecting $L_{t+1}$. 
	
	To understand the state-dependent effect from the first term of \eqref{eq:CobbEqui}, we examine the (local) dual problem of \eqref{eq:CobbEqui}, even though strong duality for the full problem over $t=1, \ldots, T-1$ remains an open question. Denote the wage function as $w(m_{t+1}; n_t, m_t)$, which depends on the current state $(n_t, m_t)$. In the dual problem of \eqref{eq:CobbEqui}, the firm $n_{t+1}$ selects the optimal worker that minimizes the sum of costs and wages:
	\begin{align*}
		\inf_{m_{t+1}} \Big\{ w(m_{t+1}; n_t, m_t) - \alpha e^{- \frac{(n_{t+1} - n_t)^2 + (m_{t+1} - m_t)^2}{\tau}} + \beta c(n_{t+1}, m_{t+1}) + \beta L_{t+1}(n_{t+1}, m_{t+1}) \Big\}.
	\end{align*}
	Informally, regarding $m_{t+1}$ as a continuous variable and assuming sufficient smoothness, the first-order condition yields
	\begin{equation}\label{eq:FOC}
		\begin{aligned}
			w'(m_{t+1}; n_t, m_t) = & - \frac{2\alpha(m_{t+1} - m_t)}{\tau} e^{- \frac{(n_{t+1} - n_t)^2 + (m_{t+1} - m_t)^2}{\tau}} \\
			& + \beta C_1 S^a(n_{t+1}) Q'(m_{t+1}) - \beta \partial_{m_{t+1}}L_{t+1}(n_{t+1}, m_{t+1}).
		\end{aligned}
	\end{equation}
	Note that the second term in \eqref{eq:FOC} is negative. If $n_{t+1} = n_t$ and the first term on the right-hand side dominates, it influences the wage as follows:
	\begin{itemize}
		\item When $m_{t+1} > m_t$, meaning the worker at time $t+1$ is less talented than $m_t$, the first term becomes negative. This causes wages to decrease more rapidly as $m_{t+1}$ increases.
		\item When $m_{t+1} < m_t$, meaning the worker at time $t+1$ is more talented than $m_t$, the first term is positive. Combined with the negative second term, this leads to slower wage increases, or even decreases, as $m_{t+1}$ declines.
	\end{itemize}
	Overall, due to the state-dependent component, the wage function $w(\cdot; n_t, m_t)$ may attain a local maximum at $m_{t+1} = m_t$. This creates an incentive for workers whose talent levels are close to $m_t$ to work in firm $n_{t+1} = n_t$.
	
	The above analysis is informal, as it neglects the dependence of the equilibrium transport plan $\pi^*(dn_{t+2:T}, dm_{t+2:T} | n_{t+1}, m_{t+1})$ on the parameter $\alpha$, and thus overlooks its impact on $L_{t+1}(n_{t+1}, m_{t+1})$. We investigate this effect numerically in the next subsection.
	
	\subsection{Numerical analysis}\label{sec:Cobb_num}
	In practice, worker talents typically need to be estimated, and estimating the corresponding transition matrix poses additional challenges, particularly when data are limited. These statistical issues are beyond the scope of this paper and merit separate study; see \cite{demerjian2012quantifying} for a related approach. As a result, we rely on synthetic data rather than empirical observations. 
	
	The model setup includes five types of firms and five types of workers, each labeled $\{1, 2, 3, 4, 5\}$. The firm transition matrix is specified as:
	\begin{equation}
		\begin{aligned}
			\mu(n_{t+1} =1 | n_t = 1) &= 1.0, \\
			\mu(n_{t+1} =2 | n_t = 2) &= 0.7, \quad \mu(n_{t+1} =3 | n_t = 2) = 0.3, \\
			\mu(n_{t+1} =2 | n_t = 3) &= 0.3, \quad \mu(n_{t+1} =3 | n_t = 3) = 0.7, \\
			\mu(n_{t+1} =4 | n_t = 4) &= 0.9, \quad \mu(n_{t+1} =5 | n_t = 4) = 0.1, \\
			\mu(n_{t+1} =4 | n_t = 5) &= 0.1, \quad \mu(n_{t+1} =5 | n_t = 5) = 0.9.
		\end{aligned}
	\end{equation}
	The worker transition matrix is given by: 
	\begin{equation}
		\begin{aligned}
			\nu(m_{t+1} = j | m_t = j) &= 0.4, \quad \text{ for } j = 1, \ldots, 5, \\
			\nu(m_{t+1} =k | m_t = j) &= 0.15, \quad \text{ for } k \neq j \text{ and } j = 1, \ldots, 5.
		\end{aligned}
	\end{equation}
	Hence, firm types are more stable than worker types. In particular, firms at the top and bottom are less likely to transition to other types. Workers retain their type with probability $0.4$; otherwise, they switch to one of the other types with equal probability. The initial distributions are defined as follows:
	\begin{equation}
		p^1_*\mu(n_1 = 1) = p^1_*\mu(n_1 = 5) = 0.125, \quad  p^1_*\mu(n_1 = 2) = p^1_*\mu(n_1 = 3) = p^1_*\mu(n_1 = 4) = 0.25,
	\end{equation}
	and $p^1_*\nu(m_1 = k) = 0.2, \; k = 1, \ldots, 5$.

	Following \cite{gabaix2008has}, we adopt the functional forms $S(n) = C_2/n^\xi$ and $Q(m) = C_3 - C_4 m^\eta$, where $C_2, C_3, C_4$ are positive constants. Empirical findings suggest $ \xi \simeq 1$, $\eta \simeq 2/3$, and $a \simeq 1$. Based on this, we define the cost function as
	\begin{equation}
		c(n, m) = - \frac{1}{n} (5^{0.6} - m^{0.6}).
	\end{equation}
	In addition, we set the discount factor to $\beta = 0.9$, the scaling constant to $\tau = 2$, and the time horizon to $T = 5$.
	
	To study the impact of $\alpha$ on $L_{t+1}(n_{t+1}, m_{t+1})$, we introduce the OT problem at time $t = 0$ without state dependence on $(n_0, m_0)$:
	\begin{equation}\label{eq:prob0}
		\inf_{\gamma \in \Pi(p^1_*\mu, p^1_*\nu) } \int \left[ \beta c(n_1, m_1) + \beta L_1(n_1, m_1) \right] \gamma(dn_1, dm_1). 
	\end{equation}
	In this formulation, it is equivalent to treat $L_1(n_1, m_1)$ as the expected costs adjusted by state dependence. The problem does not include a state-dependent term on $(n_0, m_0)$, as the focus is on examining the impact of $\alpha$ on $L_1(n_1, m_1)$ through the primal and dual solutions to \eqref{eq:prob0}.
	
\begin{figure}
		\centering
		\includegraphics[width=0.4\textwidth]{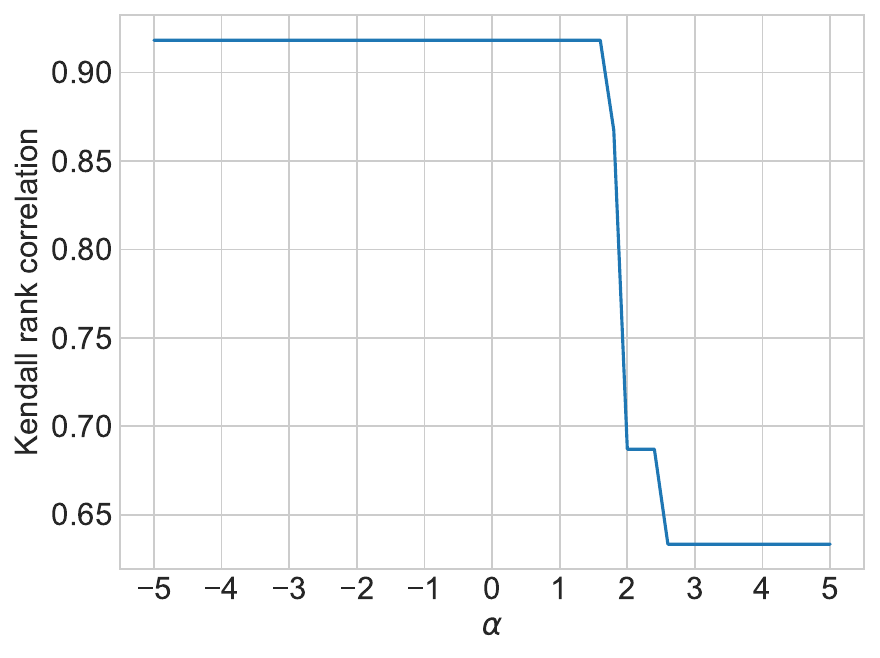}
		\caption{Kendall rank correlation decreases when state dependence is stronger in later periods. \label{fig:Cobb_corr} }
	\end{figure}

	We analyze mismatches in the optimal coupling $\pi^*(n_1, m_1)$ for the problem \eqref{eq:prob0} at time 0. Specifically, we compute the Kendall rank correlation between the types of matched firms and workers, weighted by the probability measure $\pi^*(n_1, m_1)$. As the state-dependent coefficient $\alpha$ varies, the loss function $L_1(n_1, m_1)$ changes accordingly. As shown in Figure \ref{fig:Cobb_corr}, the Kendall rank correlation decreases as $\alpha$ increases. It indicates that mismatches in $\pi^*(n_1, m_1)$ become more pronounced when state dependence is stronger in later time periods ($t=1, \ldots, T-1$).  In this example, the correlation remains constant at $0.918$ for $\alpha \in [-5.0, 1.6]$, primarily due to marginal constraints. Even in the absence of state dependence ($\alpha = 0.0$), minor mismatches persist, so the correlation does not reach one.
	
	Next, to test the model in \eqref{eq:Cobb} and \eqref{eq:prob0} on data, we perform calibrations similar to those in Sections \ref{sec:exec} and \ref{sec:prof}. Since the dual of the equilibrium transport problem in \eqref{eq:Cobb} over the entire horizon is unclear, we calibrate $\alpha$ using the method in \eqref{eq:wass_alpha}, which adopts a primal perspective.
	
\begin{figure}
		\centering
			\begin{minipage}{0.45\textwidth}
				\centering
				\includegraphics[width=0.95\textwidth]{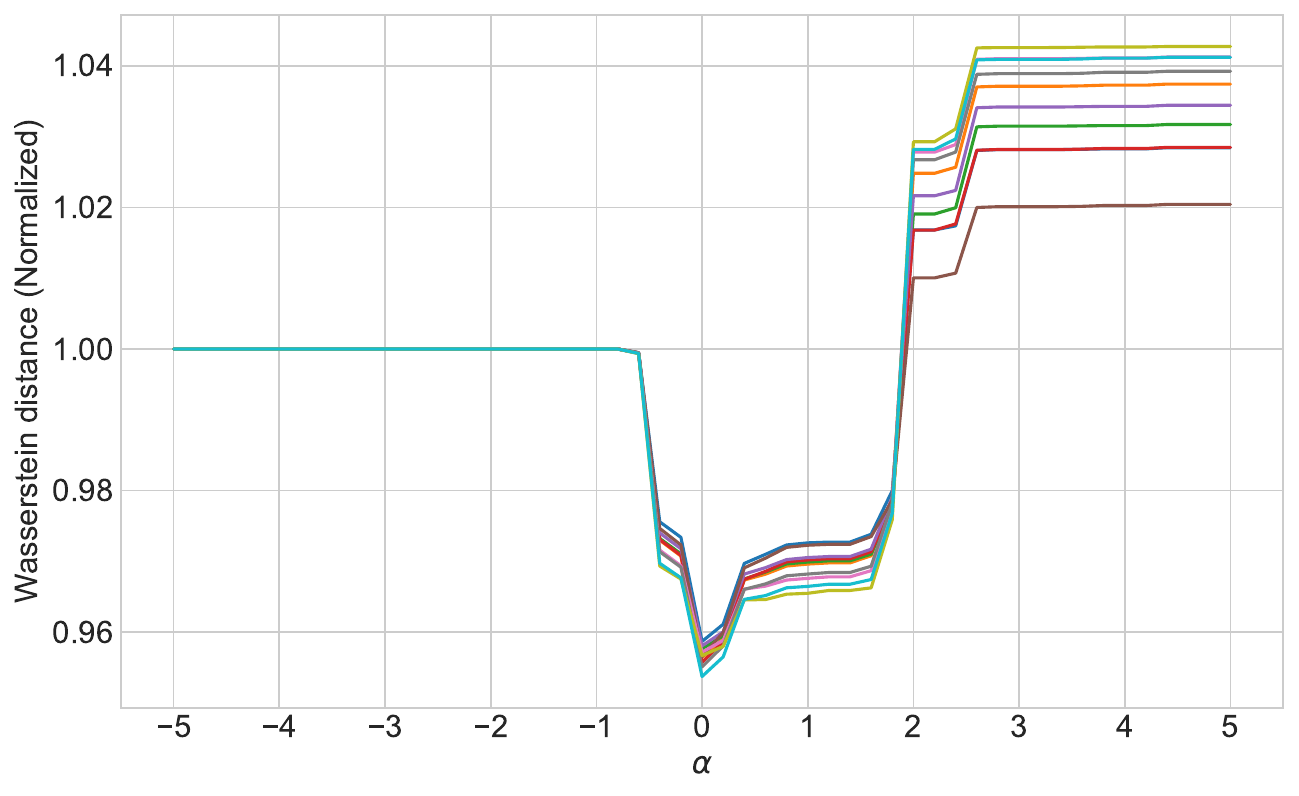}
				\subcaption{Perfectly matched synthetic data}\label{fig:matched_corr}
			\end{minipage}
			\begin{minipage}{0.45\textwidth}
				\centering
				\includegraphics[width=0.95\textwidth]{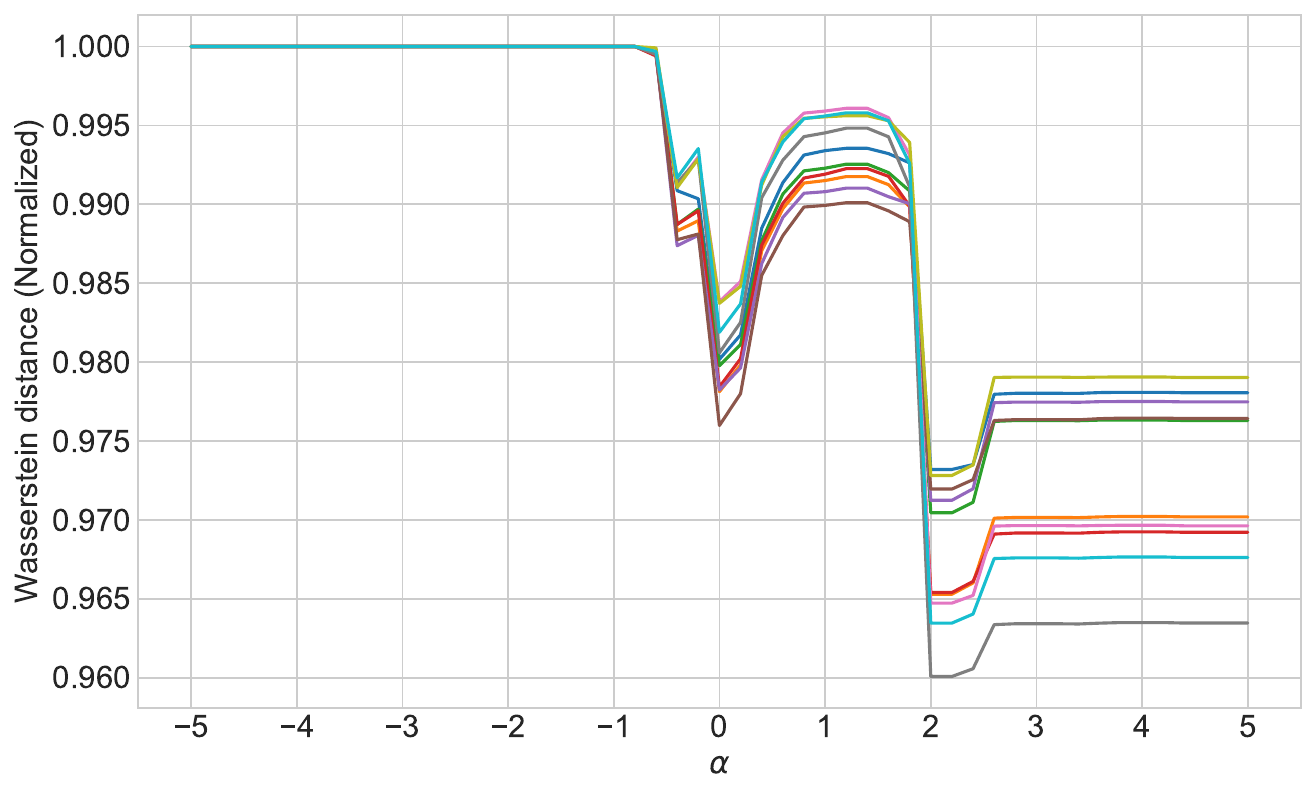}
				\subcaption{Mismatched synthetic data}\label{fig:mismatched_corr}
			\end{minipage}%
		\caption{Calibration curves with matched and mismatched synthetic data. \label{fig:CaliCobb}}
	\end{figure}

	We generate synthetic data with perfectly matched pairs and calibrate $\alpha$ over ten simulation runs. Figure \ref{fig:matched_corr} shows that the optimal value of $\alpha$ is zero. In contrast, Figure \ref{fig:mismatched_corr} adopts synthetic data where mismatches mainly occur among medium types. The results in Figure \ref{fig:mismatched_corr} indicate that the optimal value of $\alpha$ is approximately $2.0$, suggesting the presence of state-dependent effects.

	Overall, Figures \ref{fig:Cobb_corr} and \ref{fig:CaliCobb} both indicate a positive relationship between mismatches and state dependence, consistent with the findings in Section \ref{sec:reduced_jobm}.
	
\begin{figure}
		\centering
			\begin{minipage}{0.45\textwidth}
				\centering
				\includegraphics[width=0.95\textwidth]{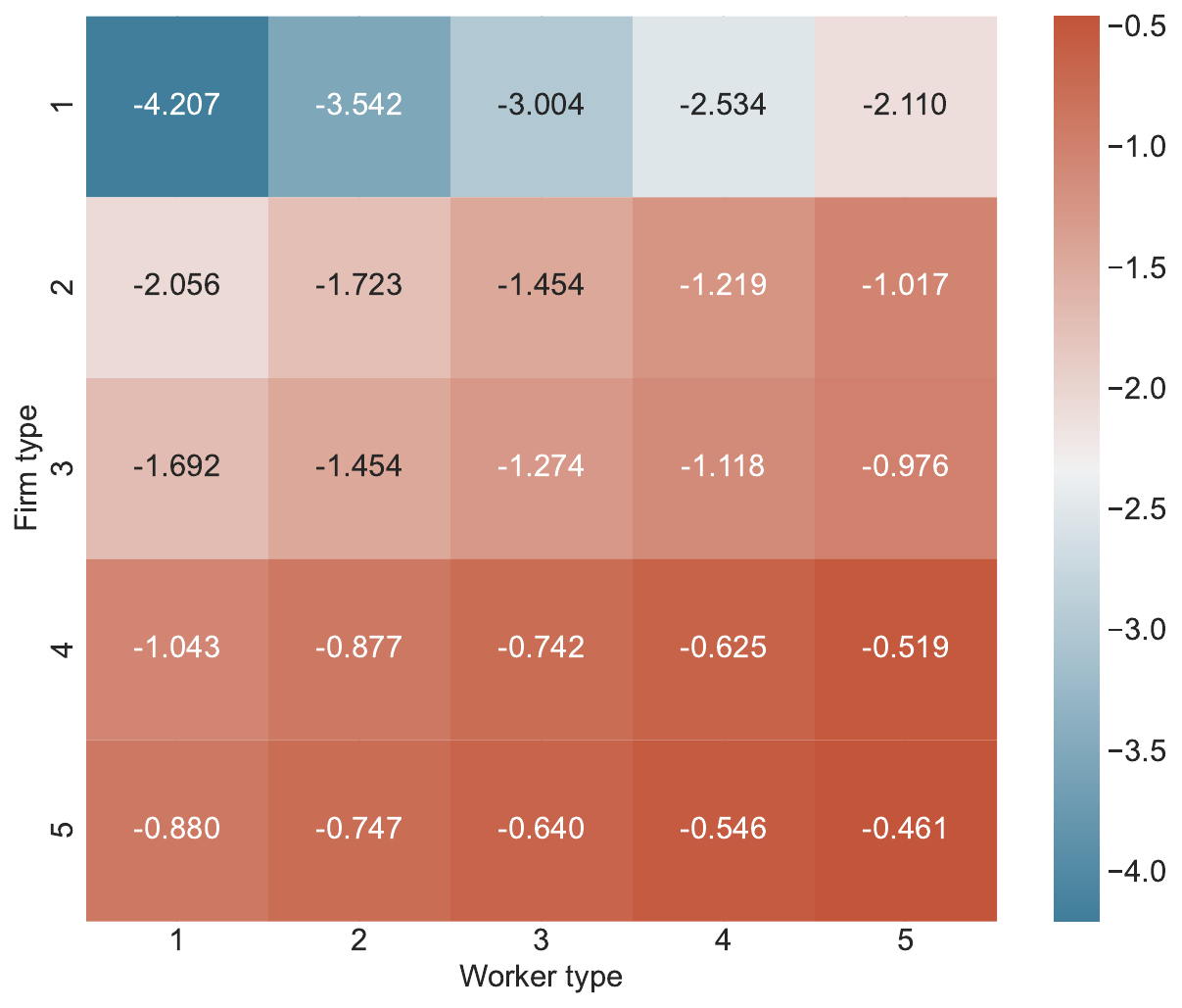}
				\subcaption{$\alpha=0.0$}\label{fig:matched_cost}
			\end{minipage}
			\begin{minipage}{0.45\textwidth}
				\centering
				\includegraphics[width=0.95\textwidth]{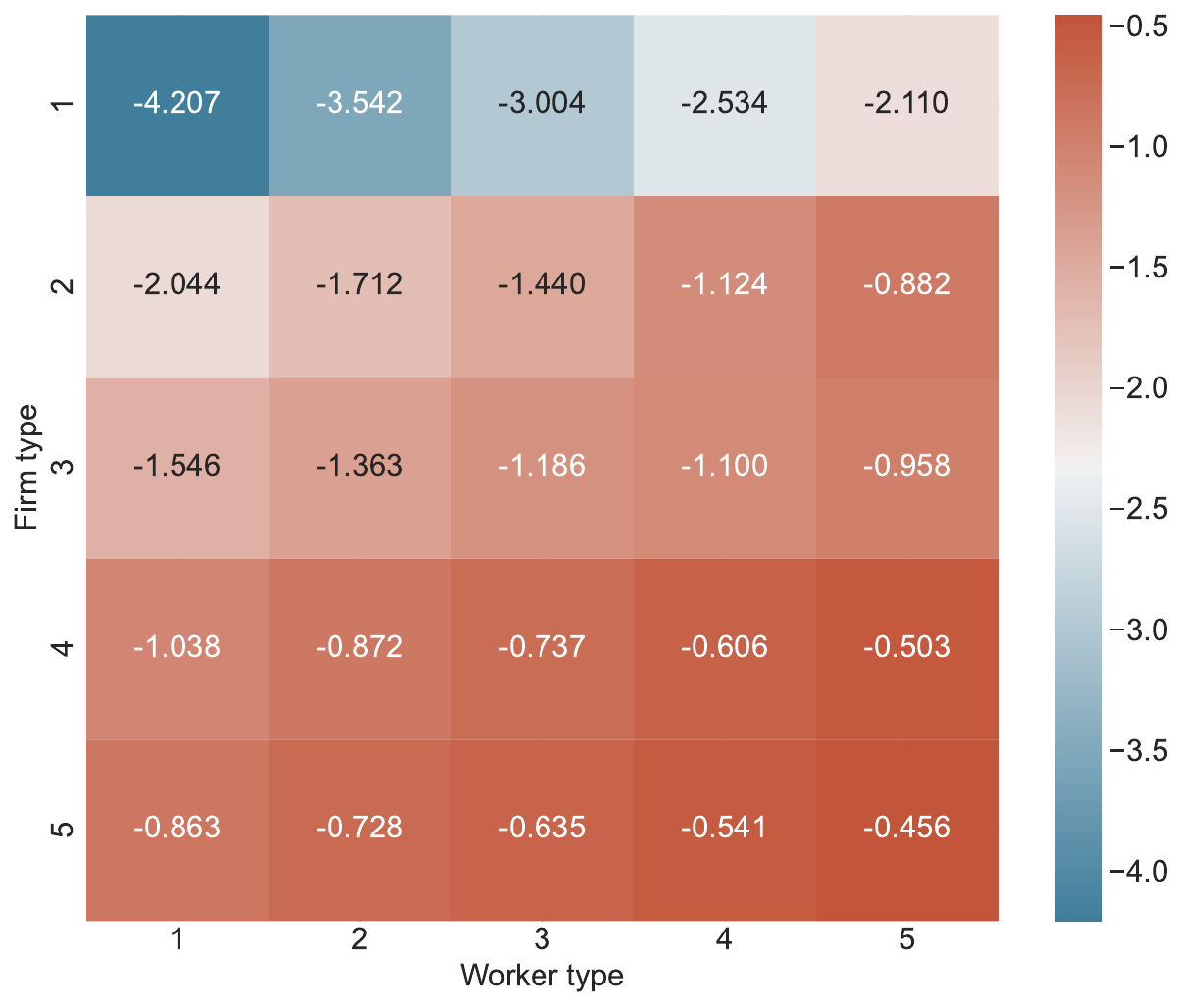}
				\subcaption{$\alpha=2.0$}\label{fig:mismatched_cost}
			\end{minipage}%
		\caption{Cost matrices $c(n_1, m_1) + L_1(n_1, m_1)$. \label{fig:cost_mat}}
	\end{figure}

	To provide a deeper analysis, following Figure \ref{fig:CaliCobb}, we further examine the problem \eqref{eq:prob0} at time $0$, using $L_1(n_1, m_1)$ implied by $\alpha = 0.0$ and $\alpha = 2.0$, respectively. Figure \ref{fig:matched_cost} displays the undiscounted cost matrix $c(n_1, m_1) + L_1(n_1, m_1)$ when $\alpha = 0.0$, while Figure \ref{fig:mismatched_cost} corresponds to $\alpha = 2.0$. The costs associated with firm type $1$ remain unchanged, whereas those for other types differ. A numerical check reveals that the cost matrix for $\alpha = 0.0$ is submodular, favoring to match larger firms with more talented workers. In contrast, the cost matrix for $\alpha = 2.0$ is not always submodular. For instance, submodularity is violated at pairs involving firm types $\{3, 4\}$ and worker types $\{3, 4\}$. Moreover, the cost matrix for $\alpha = 2.0$ is not always supermodular either.
	
\begin{figure}
		\centering
			\begin{minipage}{0.42\textwidth}
				\centering
				\includegraphics[width=0.95\textwidth]{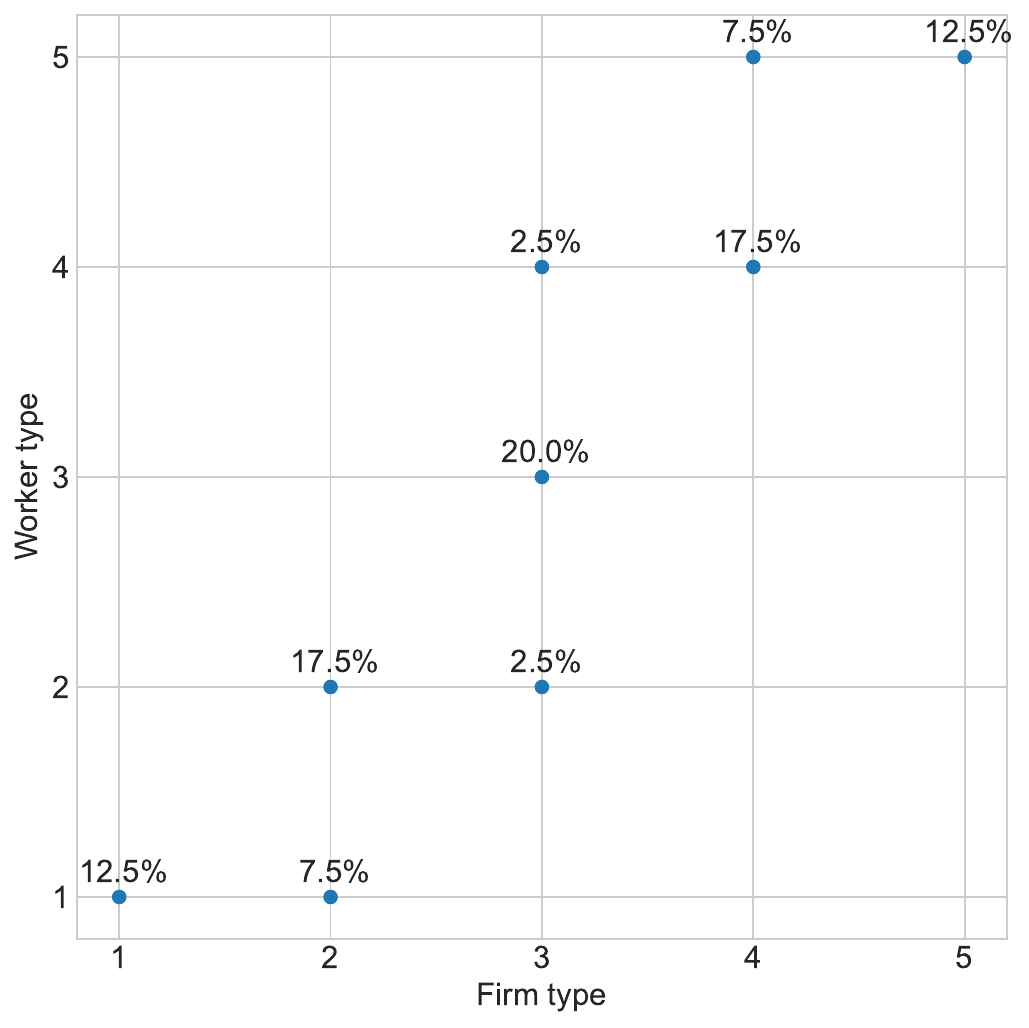}
				\subcaption{$\alpha=0.0$}\label{fig:matched_fw}
			\end{minipage}
			\begin{minipage}{0.42\textwidth}
				\centering
				\includegraphics[width=0.95\textwidth]{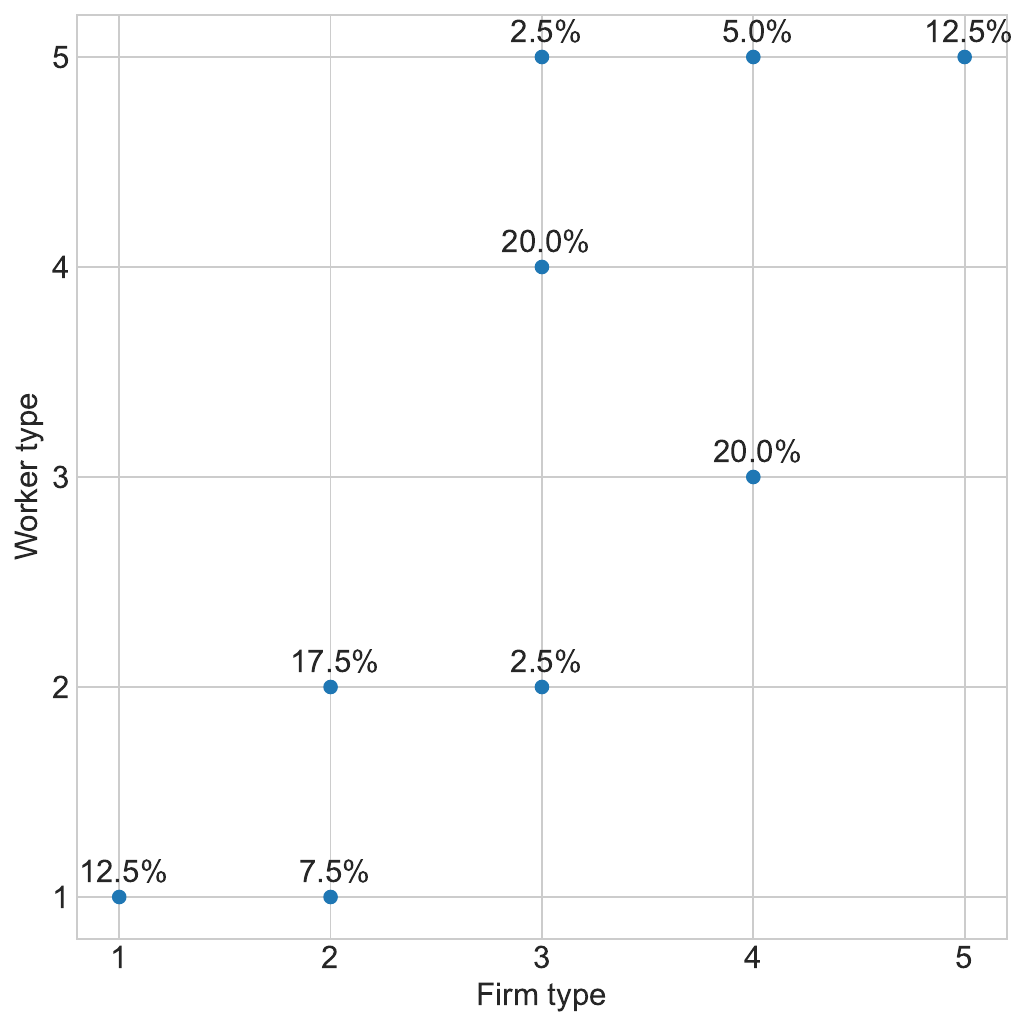}
				\subcaption{$\alpha=2.0$}\label{fig:mismatched_fw}
			\end{minipage}%
		\caption{Optimal coupling $\pi^*(n_1, m_1)$ for problem \eqref{eq:prob0}. \label{fig:Cobb_coupling}}
	\end{figure}

	Using the cost matrices in Figure \ref{fig:cost_mat}, Figure \ref{fig:Cobb_coupling} presents the optimal couplings, where the percentages indicate the probability assigned to each pair. Due to state dependence, the mismatch rate increases from $20\%$ in Figure \ref{fig:matched_fw} to $57.5\%$ in Figure \ref{fig:mismatched_fw}. Most mismatches occur between types 3 and 4, as the cost matrix violates submodularity in this region. In contrast, top and bottom types are more positively matched.

	\begin{table}
		\centering
		\begin{tabular}{cccccc}
			\hline
			Worker type  & 1 & 2 & 3 & 4 & 5 \\
			\hline
			Wages ($\alpha = 0.0$) & 0.697 & 0.398 & 0.236 & 0.095 & 0.0 \\ 
			\hline
			Wages ($\alpha = 2.0$) & 0.663 & 0.364 & 0.211 & 0.127 & 0.0  \\
			\hline
		\end{tabular}
		\caption{Wages in problem \eqref{eq:prob0} with $L_1(n_1, m_1)$ implied by different $\alpha$. The lowest wage is normalized to be zero.} \label{tab:Cobb_wages}
	\end{table}
	
Table \ref{tab:Cobb_wages} reports the wages derived from the dual problem of \eqref{eq:prob0}, with the lowest wage normalized to zero. In both cases, wages remain increasing in talent. Therefore, the mismatches between firms and wage levels are driven by the mismatches between firms and talents, as shown in Figure \ref{fig:Cobb_coupling}. A firm may employ workers with varying talent levels and thus offer different salaries based solely on talent. It mainly happens among firms of medium types $\{2, 3, 4\}$.

\section*{Acknowledgment}
The authors express gratitude to the anonymous referees and editors for their valuable comments and suggestions that have greatly improved this manuscript. Erhan Bayraktar is partially supported by the National Science Foundation under grant DMS-2106556 and by the Susan M. Smith chair. Bingyan Han is partially supported by The Hong Kong University of Science and Technology (Guangzhou) Start-up Fund G0101000197, the Guangzhou-HKUST(GZ) Joint Funding Program (No. 2024A03J0630), and the National Natural Science Foundation of China (Grant No. 12401621).  This work was partially conducted when Bingyan Han was a postdoctoral researcher in the Department of Mathematics at the University of Michigan. He expresses gratitude to the University of Michigan for providing support and an atmosphere conducive to this work. 


\begin{thebibliography}{}
	
	\bibitem[Abowd et~al., 1999]{abowd1999high}
	Abowd, J.~M., Kramarz, F., and Margolis, D.~N. (1999).
	\newblock High wage workers and high wage firms.
	\newblock {\em Econometrica}, 67(2):251--333.
	
	\bibitem[Acciaio et~al., 2021]{acciaio2021cournot}
	Acciaio, B., Backhoff-Veraguas, J., and Jia, J. (2021).
	\newblock {Cournot--Nash} equilibrium and optimal transport in a dynamic
	setting.
	\newblock {\em SIAM Journal on Control and Optimization}, 59(3):2273--2300.
	
	\bibitem[Acciaio et~al., 2020]{acciaio2020causal}
	Acciaio, B., Backhoff-Veraguas, J., and Zalashko, A. (2020).
	\newblock Causal optimal transport and its links to enlargement of filtrations
	and continuous-time stochastic optimization.
	\newblock {\em Stochastic Processes and their Applications}, 130(5):2918--2953.
	
	\bibitem[Arjovsky et~al., 2017]{arjovsky2017wasserstein}
	Arjovsky, M., Chintala, S., and Bottou, L. (2017).
	\newblock Wasserstein generative adversarial networks.
	\newblock In {\em International Conference on Machine Learning}, pages
	214--223. PMLR.
	
	\bibitem[Backhoff-Veraguas et~al., 2020]{backhoff2020adapted}
	Backhoff-Veraguas, J., Bartl, D., Beiglb{\"o}ck, M., and Eder, M. (2020).
	\newblock Adapted {Wasserstein} distances and stability in mathematical
	finance.
	\newblock {\em Finance and Stochastics}, 24(3):601--632.
	
	\bibitem[Backhoff-Veraguas et~al., 2022]{backhoff2020estimating}
	Backhoff-Veraguas, J., Bartl, D., Beiglb{\"o}ck, M., and Wiesel, J. (2022).
	\newblock Estimating processes in adapted {Wasserstein} distance.
	\newblock {\em The Annals of Applied Probability}, 32(1):529--550.
	
	\bibitem[Backhoff-Veraguas et~al., 2017]{backhoff2017causal}
	Backhoff-Veraguas, J., Beiglbock, M., Lin, Y., and Zalashko, A. (2017).
	\newblock Causal transport in discrete time and applications.
	\newblock {\em SIAM Journal on Optimization}, 27(4):2528--2562.
	
	\bibitem[Backhoff-Veraguas and Zhang, 2023]{backhoff2023dynamic}
	Backhoff-Veraguas, J. and Zhang, X. (2023).
	\newblock Dynamic {Cournot-Nash} equilibrium: The non-potential case.
	\newblock {\em Mathematics and Financial Economics}, 17(2):153--174.
	
	\bibitem[Barberis, 2012]{barberis2012model}
	Barberis, N. (2012).
	\newblock A model of casino gambling.
	\newblock {\em Management Science}, 58(1):35--51.
	
	\bibitem[Basak and Chabakauri, 2010]{basak10}
	Basak, S. and Chabakauri, G. (2010).
	\newblock Dynamic mean-variance asset allocation.
	\newblock {\em The Review of Financial Studies}, 23(8):2970--3016.
	
	\bibitem[Bayraktar et~al., 2025]{bayraktar2022stability}
	Bayraktar, E., Eckstein, S., and Zhang, X. (2025).
	\newblock Stability and sample complexity of divergence regularized optimal
	transport.
	\newblock {\em Bernoulli}, 31(1):213--239.
	
	\bibitem[Bayraktar and Han, 2023]{bayraktar2023existence}
	Bayraktar, E. and Han, B. (2023).
	\newblock Existence of {Markov} equilibrium control in discrete time.
	\newblock {\em SIAM Journal on Financial Mathematics}, 14(4):SC60--SC71.
	
	\bibitem[Bayraktar et~al., 2021]{bayraktar2021equilibrium}
	Bayraktar, E., Zhang, J., and Zhou, Z. (2021).
	\newblock Equilibrium concepts for time-inconsistent stopping problems in
	continuous time.
	\newblock {\em Mathematical Finance}, 31(1):508--530.
	
	\bibitem[Becker, 1973]{becker1973theory}
	Becker, G.~S. (1973).
	\newblock A theory of marriage: Part {I}.
	\newblock {\em Journal of Political Economy}, 81(4):813--846.
	
	\bibitem[Beiglb{\"o}ck et~al., 2022]{beiglbock2022approx}
	Beiglb{\"o}ck, M., Jourdain, B., Margheriti, W., and Pammer, G. (2022).
	\newblock Approximation of martingale couplings on the line in the adapted weak
	topology.
	\newblock {\em Probability Theory and Related Fields}, 183(1):359--413.
	
	\bibitem[Beiglb{\"o}ck and Pratelli, 2012]{beiglbock2012duality}
	Beiglb{\"o}ck, M. and Pratelli, A. (2012).
	\newblock Duality for rectified cost functions.
	\newblock {\em Calculus of Variations and Partial Differential Equations},
	45:27--41.
	
	\bibitem[Bertsekas and Shreve, 1978]{bertsekas1978stoch}
	Bertsekas, D. and Shreve, S.~E. (1978).
	\newblock {\em Stochastic Optimal Control: The Discrete-time Case}.
	\newblock Academic Press.
	
	\bibitem[Bj{\"o}rk et~al., 2017]{bjork17FS}
	Bj{\"o}rk, T., Khapko, M., and Murgoci, A. (2017).
	\newblock On time-inconsistent stochastic control in continuous time.
	\newblock {\em Finance and Stochastics}, 21(2):331--360.
	
	\bibitem[Bj{\"o}rk and Murgoci, 2014]{bjork14FS}
	Bj{\"o}rk, T. and Murgoci, A. (2014).
	\newblock A theory of {Markovian} time-inconsistent stochastic control in
	discrete time.
	\newblock {\em Finance and Stochastics}, 18(3):545--592.
	
	\bibitem[Bj{\"o}rk et~al., 2014]{bjork14MF}
	Bj{\"o}rk, T., Murgoci, A., and Zhou, X.~Y. (2014).
	\newblock Mean--variance portfolio optimization with state-dependent risk
	aversion.
	\newblock {\em Mathematical Finance}, 24(1):1--24.
	
	\bibitem[Blanchet and Murthy, 2019]{blanchet2019quantifying}
	Blanchet, J. and Murthy, K. (2019).
	\newblock Quantifying distributional model risk via optimal transport.
	\newblock {\em Mathematics of Operations Research}, 44(2):565--600.
	
	\bibitem[Blanchet et~al., 2021]{blanchet2021stat}
	Blanchet, J., Murthy, K., and Nguyen, V.~A. (2021).
	\newblock Statistical analysis of {Wasserstein} distributionally robust
	estimators.
	\newblock In {\em Tutorials in Operations Research: Emerging Optimization
		Methods and Modeling Techniques with Applications}, pages 227--254. INFORMS.
	
	\bibitem[Boerma et~al., 2023]{boerma2023composite}
	Boerma, J., Tsyvinski, A., Wang, R., and Zhang, Z. (2023).
	\newblock Composite sorting.
	\newblock {\em arXiv preprint arXiv:2303.06701}.
	
	\bibitem[Bogachev, 2007]{bogachev2007measure}
	Bogachev, V.~I. (2007).
	\newblock {\em Measure Theory}, volume~II.
	\newblock Springer Science \& Business Media.
	
	\bibitem[Bonhomme et~al., 2019]{bonhomme2019}
	Bonhomme, S., Lamadon, T., and Manresa, E. (2019).
	\newblock A distributional framework for matched employer employee data.
	\newblock {\em Econometrica}, 87(3):699--739.
	
	\bibitem[Borovickov{\'a} and Shimer, 2020]{borovickova2020high}
	Borovickov{\'a}, K. and Shimer, R. (2020).
	\newblock High wage workers work for high wage firms.
	
	\bibitem[Brenier, 1991]{brenier1991polar}
	Brenier, Y. (1991).
	\newblock Polar factorization and monotone rearrangement of vector-valued
	functions.
	\newblock {\em Communications on Pure and Applied Mathematics}, 44(4):375--417.
	
	\bibitem[Brown and Purves, 1973]{brown1973}
	Brown, L.~D. and Purves, R. (1973).
	\newblock Measurable selections of extrema.
	\newblock {\em The Annals of Statistics}, pages 902--912.
	
	\bibitem[Card et~al., 2013]{card2013workplace}
	Card, D., Heining, J., and Kline, P. (2013).
	\newblock Workplace heterogeneity and the rise of {West German} wage
	inequality.
	\newblock {\em The Quarterly Journal of Economics}, 128(3):967--1015.
	
	\bibitem[Charalambos and Aliprantis, 2013]{charalambos2013infinite}
	Charalambos, D. and Aliprantis, B. (2013).
	\newblock {\em Infinite Dimensional Analysis: A Hitchhiker's Guide}.
	\newblock Springer.
	
	\bibitem[Condat, 2013]{condat2013direct}
	Condat, L. (2013).
	\newblock A direct algorithm for {1-D} total variation denoising.
	\newblock {\em IEEE Signal Processing Letters}, 20(11):1054--1057.
	
	\bibitem[Cuturi, 2013]{cuturi2013sinkhorn}
	Cuturi, M. (2013).
	\newblock Sinkhorn distances: Lightspeed computation of optimal transport.
	\newblock {\em Advances in Neural Information Processing Systems}, 26.
	
	\bibitem[Delon and Desolneux, 2020]{delon2020wasserstein}
	Delon, J. and Desolneux, A. (2020).
	\newblock A {Wasserstein}-type distance in the space of {Gaussian} mixture
	models.
	\newblock {\em SIAM Journal on Imaging Sciences}, 13(2):936--970.
	
	\bibitem[Demerjian et~al., 2012]{demerjian2012quantifying}
	Demerjian, P., Lev, B., and McVay, S. (2012).
	\newblock Quantifying managerial ability: A new measure and validity tests.
	\newblock {\em Management Science}, 58(7):1229--1248.
	
	\bibitem[Eckstein and Pammer, 2024]{eckstein2024}
	Eckstein, S. and Pammer, G. (2024).
	\newblock Computational methods for adapted optimal transport.
	\newblock {\em The Annals of Applied Probability}, 34(1A):675--713.
	
	\bibitem[Epstein and Ji, 2022]{epstein2022optimal}
	Epstein, L.~G. and Ji, S. (2022).
	\newblock Optimal learning under robustness and time-consistency.
	\newblock {\em Operations Research}, 70(3):1317--1329.
	
	\bibitem[F{\"o}llmer and Schied, 2011]{follmer2011stochastic}
	F{\"o}llmer, H. and Schied, A. (2011).
	\newblock {\em Stochastic Finance: An Introduction in Discrete Time}.
	\newblock Walter de Gruyter.
	
	\bibitem[Gabaix and Landier, 2008]{gabaix2008has}
	Gabaix, X. and Landier, A. (2008).
	\newblock Why has {CEO} pay increased so much?
	\newblock {\em The Quarterly Journal of Economics}, 123(1):49--100.
	
	\bibitem[Galichon, 2016]{galichon2016optimal}
	Galichon, A. (2016).
	\newblock {\em Optimal Transport Methods in Economics}.
	\newblock Princeton University Press.
	
	\bibitem[Gangbo and McCann, 1996]{gangbo1996geometry}
	Gangbo, W. and McCann, R.~J. (1996).
	\newblock The geometry of optimal transportation.
	\newblock {\em Acta Mathematica}, 177(2):113--161.
	
	\bibitem[Gao and Kleywegt, 2022]{gao2016dist}
	Gao, R. and Kleywegt, A. (2022).
	\newblock Distributionally robust stochastic optimization with {Wasserstein}
	distance.
	\newblock {\em Mathematics of Operations Research}.
	
	\bibitem[Givens and Shortt, 1984]{givens1984class}
	Givens, C.~R. and Shortt, R.~M. (1984).
	\newblock A class of {Wasserstein} metrics for probability distributions.
	\newblock {\em Michigan Mathematical Journal}, 31(2):231--240.
	
	\bibitem[Gonz{\'a}lez-Sanz and Nutz, 2024]{gonzalez2024quad}
	Gonz{\'a}lez-Sanz, A. and Nutz, M. (2024).
	\newblock Quantitative convergence of quadratically regularized linear
	programs.
	\newblock {\em arXiv preprint arXiv:2408.04088}.
	
	\bibitem[Gunasingam and Wong, 2024]{gunasingam2024adapted}
	Gunasingam, M. and Wong, T.-K.~L. (2024).
	\newblock Adapted optimal transport between {Gaussian} processes in discrete
	time.
	\newblock {\em arXiv preprint arXiv:2404.06625}.
	
	\bibitem[Hagedorn et~al., 2017]{hagedorn2017identifying}
	Hagedorn, M., Law, T.~H., and Manovskii, I. (2017).
	\newblock Identifying equilibrium models of labor market sorting.
	\newblock {\em Econometrica}, 85(1):29--65.
	
	\bibitem[Han, 2025]{han2025IEEE}
	Han, B. (2025).
	\newblock Distributionally robust {Kalman} filtering with volatility
	uncertainty.
	\newblock {\em IEEE Transactions on Automatic Control}.
	
	\bibitem[Han et~al., 2021]{han2021robust}
	Han, B., Pun, C.~S., and Wong, H.~Y. (2021).
	\newblock Robust state-dependent mean--variance portfolio selection: A
	closed-loop approach.
	\newblock {\em Finance and Stochastics}, 25(3):529--561.
	
	\bibitem[Hu and Zhou, 2022]{hu2022dynamic}
	Hu, M. and Zhou, Y. (2022).
	\newblock Dynamic type matching.
	\newblock {\em Manufacturing \& Service Operations Management}, 24(1):125--142.
	
	\bibitem[Jenks, 1967]{jenks}
	Jenks, G.~F. (1967).
	\newblock The data model concept in statistical mapping.
	\newblock {\em International Yearbook of Cartography}, 7:186--190.
	
	\bibitem[Kahneman and Tversky, 1979]{kahneman1979prospect}
	Kahneman, D. and Tversky, A. (1979).
	\newblock Prospect theory: An analysis of decision under risk.
	\newblock {\em Econometrica}, 47(2):263--292.
	
	\bibitem[Kallenberg, 2021]{kallenberg}
	Kallenberg, O. (2021).
	\newblock {\em Foundations of Modern Probability}.
	\newblock Springer Science \& Business Media.
	\newblock The third edition.
	
	\bibitem[Kechris, 2012]{kechris2012classical}
	Kechris, A. (2012).
	\newblock {\em Classical Descriptive Set Theory}, volume 156.
	\newblock Springer Science \& Business Media.
	
	\bibitem[Kov{\'a}{\v{c}}ov{\'a} and Rudloff, 2021]{kovavcova2021time}
	Kov{\'a}{\v{c}}ov{\'a}, G. and Rudloff, B. (2021).
	\newblock Time consistency of the mean-risk problem.
	\newblock {\em Operations research}, 69(4):1100--1117.
	
	\bibitem[Kuhn et~al., 2019]{kuhn2019tutorial}
	Kuhn, D., Esfahani, P.~M., Nguyen, V.~A., and Shafieezadeh-Abadeh, S. (2019).
	\newblock Wasserstein distributionally robust optimization: Theory and
	applications in machine learning.
	\newblock In {\em Operations research \& management science in the age of
		analytics}, pages 130--166. INFORMS.
	
	\bibitem[Laibson, 1997]{laibson1997golden}
	Laibson, D. (1997).
	\newblock Golden eggs and hyperbolic discounting.
	\newblock {\em The Quarterly Journal of Economics}, 112(2):443--478.
	
	\bibitem[Lassalle, 2013]{lassalle2013causal}
	Lassalle, R. (2013).
	\newblock Causal transference plans and their {Monge-Kantorovich} problems.
	\newblock {\em arXiv preprint arXiv:1303.6925}.
	
	\bibitem[Ma et~al., 2021]{ma2021time}
	Ma, J., Wong, T.-K.~L., and Zhang, J. (2021).
	\newblock Time-consistent conditional expectation under probability distortion.
	\newblock {\em Mathematics of Operations Research}, 46(3):1149--1180.
	
	\bibitem[Mohajerin~Esfahani and Kuhn, 2018]{mohajerin2018data}
	Mohajerin~Esfahani, P. and Kuhn, D. (2018).
	\newblock Data-driven distributionally robust optimization using the
	{Wasserstein} metric: Performance guarantees and tractable reformulations.
	\newblock {\em Mathematical Programming}, 171(1):115--166.
	
	\bibitem[Neufeld and Sester, 2021]{neufeld2021stability}
	Neufeld, A. and Sester, J. (2021).
	\newblock On the stability of the martingale optimal transport problem: A
	set-valued map approach.
	\newblock {\em Statistics \& Probability Letters}, 176:109131.
	
	\bibitem[Parthasarathy, 2005]{parthasarathy2005}
	Parthasarathy, K.~R. (2005).
	\newblock {\em Probability Measures on Metric Spaces}, volume 352.
	\newblock American Mathematical Soc.
	
	\bibitem[Peyr{\'e} and Cuturi, 2019]{peyre2019computational}
	Peyr{\'e}, G. and Cuturi, M. (2019).
	\newblock Computational optimal transport: With applications to data science.
	\newblock {\em Foundations and Trends{\textregistered} in Machine Learning},
	11(5-6):355--607.
	
	\bibitem[Pflug and Pichler, 2012]{pflug2012distance}
	Pflug, G.~C. and Pichler, A. (2012).
	\newblock A distance for multistage stochastic optimization models.
	\newblock {\em SIAM Journal on Optimization}, 22(1):1--23.
	
	\bibitem[Pflug and Pichler, 2014]{pflug2014multistage}
	Pflug, G.~C. and Pichler, A. (2014).
	\newblock {\em Multistage Stochastic Optimization}, volume 1104.
	\newblock Springer.
	
	\bibitem[Pichler et~al., 2022]{pichler2022risk}
	Pichler, A., Liu, R.~P., and Shapiro, A. (2022).
	\newblock Risk-averse stochastic programming: Time consistency and optimal
	stopping.
	\newblock {\em Operations Research}, 70(4):2439--2455.
	
	\bibitem[Pichler and Weinhardt, 2022]{pichler2022nested}
	Pichler, A. and Weinhardt, M. (2022).
	\newblock The nested {Sinkhorn} divergence to learn the nested distance.
	\newblock {\em Computational Management Science}, 19(2):269--293.
	
	\bibitem[Postel-Vinay and Robin, 2002]{postel2002equilibrium}
	Postel-Vinay, F. and Robin, J.-M. (2002).
	\newblock Equilibrium wage dispersion with worker and employer heterogeneity.
	\newblock {\em Econometrica}, 70(6):2295--2350.
	
	\bibitem[Schrott et~al., 2023]{schrott2023denseness}
	Schrott, S., Beiglb{\"o}ck, M., and Pammer, G. (2023).
	\newblock Denseness of biadapted {Monge} mappings.
	\newblock {\em Annales de l’Institut Henri Poincar{\'e}-Probabilit{\'e}s et
		Statistiques}.
	
	\bibitem[Seguy et~al., 2018]{seguy2018large}
	Seguy, V., Damodaran, B.~B., Flamary, R., Courty, N., Rolet, A., and Blondel,
	M. (2018).
	\newblock Large-scale optimal transport and mapping estimation.
	\newblock In {\em International Conference on Learning Representations}, pages
	1--15.
	
	\bibitem[Shimer and Smith, 2000]{shimer2000assortative}
	Shimer, R. and Smith, L. (2000).
	\newblock Assortative matching and search.
	\newblock {\em Econometrica}, 68(2):343--369.
	
	\bibitem[Song et~al., 2019]{song2019firming}
	Song, J., Price, D.~J., Guvenen, F., Bloom, N., and Von~Wachter, T. (2019).
	\newblock Firming up inequality.
	\newblock {\em The Quarterly Journal of Economics}, 134(1):1--50.
	
	\bibitem[Steen and Seebach, 1978]{steen1978}
	Steen, L.~A. and Seebach, J.~A. (1978).
	\newblock {\em Counterexamples in Topology}.
	\newblock Springer.
	
	\bibitem[Strotz, 1955]{strotz1955myopia}
	Strotz, R. (1955).
	\newblock Myopia and inconsistency in dynamic utility maximization.
	\newblock {\em Review of Economic Studies}, 23(3):165--180.
	
	\bibitem[Stuart and Wolfram, 2020]{stuart2020inverse}
	Stuart, A.~M. and Wolfram, M.-T. (2020).
	\newblock Inverse optimal transport.
	\newblock {\em SIAM Journal on Applied Mathematics}, 80(1):599--619.
	
	\bibitem[Sundaram, 1996]{sundaram1996first}
	Sundaram, R.~K. (1996).
	\newblock {\em A First Course in Optimization Theory}.
	\newblock Cambridge University Press.
	
	\bibitem[Ta{\c{s}}kesen et~al., 2023]{tacskesen2023semi}
	Ta{\c{s}}kesen, B., Shafieezadeh-Abadeh, S., and Kuhn, D. (2023).
	\newblock Semi-discrete optimal transport: Hardness, regularization and
	numerical solution.
	\newblock {\em Mathematical Programming}, 199(1-2):1033--1106.
	
	\bibitem[Taylor, 2013]{taylor2013ceo}
	Taylor, L.~A. (2013).
	\newblock {CEO} wage dynamics: Estimates from a learning model.
	\newblock {\em Journal of Financial Economics}, 108(1):79--98.
	
	\bibitem[Torous et~al., 2021]{torous2021optimal}
	Torous, W., Gunsilius, F., and Rigollet, P. (2021).
	\newblock An optimal transport approach to causal inference.
	\newblock {\em arXiv preprint arXiv:2108.05858}.
	
	\bibitem[Villani, 2009]{villani2009optimal}
	Villani, C. (2009).
	\newblock {\em Optimal Transport: Old and New}, volume 338.
	\newblock Springer.
	
	\bibitem[Xu et~al., 2020]{xu2020cot}
	Xu, T., Li, W.~K., Munn, M., and Acciaio, B. (2020).
	\newblock {COT-GAN}: Generating sequential data via causal optimal transport.
	\newblock {\em Advances in Neural Information Processing Systems},
	33:8798--8809.
	
\end{thebibliography}

\appendix

\section{Proofs of results}\label{sec:proofs}

\subsection{The semi-discrete and Markovian case}\label{sec:proof_semi}
We need several auxiliary results to study the continuity. Lemmas \ref{lem:closed_corr} and \ref{lem:Cb_conti} do not rely on any particular choice of the metric.
\begin{lemma}\label{lem:closed_corr}
	Suppose
	\begin{itemize}
		\item[(1)] $(\cX_{t}, \cT_{\cX_t})$, $(\cX_{t+1}, \cT_{\cX_{t+1}})$,  $(\cY_{t}, \cT_{\cY_t}) $, $(\cY_{t+1}, \cT_{\cY_{t+1}})$ are Polish topological spaces and the product spaces between them are endowed with product topologies;
		\item[(2)] the stochastic kernels $\mu(dx_{t+1}|x_{t}): (\cX_{t}, \cT_{\cX_t}) \rightarrow (\cP(\cX_{t+1}), \cV[C_b(\cX_{t+1};  \cT_{\cX_{t+1}})])$ and $\nu(dy_{t+1} | y_{t}):  (\cY_{t}, \cT_{\cY_t}) \rightarrow (\cP(\cY_{t+1}), \cV[C_b(\cY_{t+1};  \cT_{\cY_{t+1}})])$ are continuous.
	\end{itemize}
	Denote a correspondence as
	\begin{align*}
		& D: (\cX_{t} \times \cY_{t}, \cT_{\cX_t} \times \cT_{\cY_t}) \twoheadrightarrow \big(\cP(\cX_{t+1} \times \cY_{t+1}), \cV[C_b(\cX_{t+1} \times \cY_{t+1}; \cT_{\cX_{t+1}} \times \cT_{\cY_{t+1}})]  \big) \\
		& \text{ that maps } (x_{t}, y_{t}) \mapsto \Pi(\mu(dx_{t+1}|x_{t}), \nu(dy_{t+1} | y_{t})).
	\end{align*}
	Then $D$ is upper hemicontinuous and $D(x_t, y_t)$ is non-empty, convex, and compact. Moreover, $D$ has a closed graph under the product topology $ \cT_{\cX_t} \times \cT_{\cY_t} \times \cV[C_b(\cX_{t+1} \times \cY_{t+1}; \cT_{\cX_{t+1}} \times \cT_{\cY_{t+1}})]$.
\end{lemma}
\begin{proof}
$\Pi(\mu(dx_{t+1}|x_{t}), \nu(dy_{t+1} | y_{t}))$ is non-empty since the independent coupling belongs to this set. This set is convex. Indeed, if $\gamma^1, \gamma^2 \in \Pi(\mu(dx_{t+1}|x_{t}), \nu(dy_{t+1} | y_{t}))$, then $\lambda \gamma^1 + (1 - \lambda) \gamma^2, \, \lambda \in [0, 1]$ is also a probability measure with marginals $\mu(dx_{t+1}|x_{t})$ and $\nu(dy_{t+1} | y_{t})$. 

Since Polish spaces are second countable, first countable, and metrizable, we can apply \citet[Theorem 17.20]{charalambos2013infinite}. We only need to show that if a sequence $\{(x^n_t, y^n_t, \gamma^n)\}$ is in the graph of $D$ and $(x^n_t, y^n_t) \rightarrow (x_t, y_t)$ under $\cT_{\cX_t} \times \cT_{\cY_t}$, then the sequence $\{\gamma^n\}$ has a limit point in $\Pi(\mu(dx_{t+1} |x_t), \nu(dy_{t+1}| y_t))$.

We have assumed $\mu(dx_{t+1} |x_t)$ and $\nu(dy_{t+1}| y_t)$ are continuous with the usual weak convergence. Thus, by Prokhorov's theorem, $\{ \mu(dx_{t+1} |x^n_t) \}^\infty_{n=1}$ and $\{ \nu(dy_{t+1}| y^n_t) \}^\infty_{n=1}$ are tight. The set of couplings $\Pi(\{ \mu(dx_{t+1} |x^n_t) \}^\infty_{n=1}, \{ \nu(dy_{t+1}| y^n_t) \}^\infty_{n=1})$ is also tight by \citet[Lemma 4.4]{villani2009optimal}, where we emphasize that the product topology $\cT_{\cX_{t+1}} \times \cT_{\cY_{t+1}}$ is imposed. Indeed, \citet[Lemma 4.4]{villani2009optimal} relies on the fact that, if $A \subset \cX_{t+1}$ and $B \subset \cY_{t+1}$ are compact, then $A \times B$ is compact under the product topology $\cT_{\cX_{t+1}} \times \cT_{\cY_{t+1}}$.

Since $\{\gamma^n\}^\infty_{n=1}$ is in $\Pi(\{ \mu(dx_{t+1} |x^n_t) \}^\infty_{n=1}, \{ \nu(dy_{t+1}| y^n_t) \}^\infty_{n=1})$, we can apply Prokhorov's theorem again. A subsequence $\{\gamma^{n_k}\}^\infty_{k=1}$ converges weakly in the usual sense to some $\gamma$. We show $\gamma \in \Pi(\mu(dx_{t+1} |x_t), \nu(dy_{t+1} |y_t))$. Denote the projection operator on the first and second component as $p_\cX : \cX_{t+1} \times \cY_{t+1} \rightarrow \cX_{t+1}$ and $p_\cY : \cX_{t+1} \times \cY_{t+1} \rightarrow \cY_{t+1}$.  With the product topology, these projection operators are continuous. Hence, if $f(x_{t+1})$ is a $\cT_{\cX_{t+1}}$-continuous function, then $f \circ p_\cX$ is $\cT_{\cX_{t+1}} \times \cT_{\cY_{t+1}}$-continuous. Therefore, we have
\begin{align*}
	\gamma^{n_k} \circ p^{-1}_\cX  & \rightarrow \gamma \circ p^{-1}_\cX, \quad \gamma^{n_k} \circ p^{-1}_\cY \rightarrow \gamma \circ p^{-1}_\cY  \text{ weakly for } k \rightarrow \infty.   
\end{align*}	
Since the marginals satisfy $\gamma^{n_k} \circ p^{-1}_\cX = \mu(dx_{t+1} |x^{n_k}_t)$ and $\gamma^{n_k} \circ p^{-1}_\cY = \nu(dy_{t+1} |y^{n_k}_t)$, together with $\mu(dx_{t+1} |x^{n_k}_t) \rightarrow \mu(dx_{t+1} |x_t)$, $\nu(dy_{t+1} |y^{n_k}_t) \rightarrow \nu(dy_{t+1} |y_t)$  weakly for $ k \rightarrow \infty$, we obtain $ \gamma \circ p^{-1}_\cX = \mu(dx_{t+1} |x_t)$ and $\gamma \circ p^{-1}_\cY = \nu(dy_{t+1} |y_t)$. Therefore, $\gamma \in \Pi(\mu(dx_{t+1} |x_t), \nu(dy_{t+1} |y_t))$.

By \citet[Theorem 17.20]{charalambos2013infinite}, $D$ is upper hemicontinuous and compact-valued. By \citet[Theorem 17.10]{charalambos2013infinite} and Polish spaces are Hausdorff, $D$ has a closed graph. 
\end{proof}

\begin{lemma}\label{lem:Cb_conti}
	Consider Polish topological spaces $(\cS_1, \cT_{\cS_1})$, $(\cS_2, \cT_{\cS_2})$, and $(\cS_3, \cT_{\cS_3})$. Suppose
	\begin{itemize}
		\item[(1)] the stochastic kernel $\gamma(ds_2| s_1): (\cS_1, \cT_{\cS_1}) \rightarrow (\cP(\cS_2), \cV[C_b(\cS_2; \cT_{\cS_2})])$ is continuous;
		\item[(2)] the function $h(s_3, s_2): (\cS_3 \times \cS_2, \cT_{\cS_3} \times \cT_{\cS_2}) \rightarrow (\R, \cT_\R)$ is continuous and bounded.
	\end{itemize}
	Then
	\begin{itemize}
		\item[(a)] $\int_{\cS_2} h(s_3, s_2) \gamma(ds_2 | s_1)$ is $\cT_{\cS_3} \times \cT_{\cS_1}$-continuous. If $(\cS_3, \cT_{\cS_3}) = (\cS_1, \cT_{\cS_1})$, it is understood as $\cT_{\cS_1}$-continuous;
		\item[(b)] with $\lambda \in \cP(\cS_2)$, $\int_{\cS_2} h(s_3, s_2) \lambda(ds_2)$ is $\cT_{\cS_3} \times \cV[C_b(\cS_2; \cT_{\cS_2})]$-continuous. 
	\end{itemize}
	
\end{lemma}

\begin{proof}
Consider a sequence $(s^n_3, s^n_1)$ converging to $(s_3, s_1)$ under the product topology $\cT_{\cS_3} \times \cT_{\cS_1}$. Equip $\cS_1$, $\cS_2$, and $\cS_3$ with some complete compatible metrics. $B := \{(s^n_3, s^n_1)\}^\infty_{n=1} \cup \{ (s_3, s_1) \}$ is a compact set. Moreover, the uniform continuity is defined under these metrics.

Fix an arbitrary $\varepsilon > 0$. Since the product topology is used, $s^n_1 \rightarrow s_1$ under $\cT_{\cS_1}$. With the assumption that $\gamma(ds_2 | s_1)$ is continuous in $s_1$, $\gamma(ds_2 | s^n_1)$ converges weakly to $\gamma(ds_2 | s_1)$. Thus $\{ \gamma(ds_2 | s^n_1) \}^\infty_{n=1}$ is tight by Prokhorov's theorem. We can find a compact subset $A \subset \cS_2$ such that  $\sup_n \gamma(\cS_2 \backslash A | s^n_1) \leq \varepsilon$.

Since $h$ is uniformly continuous on the compact set $B \times A$, there exists $N > 0$ such that
\begin{equation}
	\sup_{n > N, s_2 \in A}  | h(s_3, s_2) - h(s^n_3, s_2) | \leq \varepsilon. \label{eq:Cb_uniform}
\end{equation}
We have 
\begin{align*}
	& \Big| \int_{\cS_2} h(s^n_3, s_2) \gamma(ds_2 | s^n_1) - \int_{\cS_2} h(s_3, s_2) \gamma(ds_2 | s_1) \Big| \leq \rm{I} + \rm{II} + \rm{III},
\end{align*}
where
\begin{align*}
	\rm{I} = & \Big| \int_{\cS_2} h(s_3, s_2) [\gamma(ds_2 | s^n_1) - \gamma(ds_2 | s_1)] \Big|, \\
	\rm{II} = & \Big| \int_A [h(s^n_3, s_2) - h(s_3, s_2)] \gamma(ds_2 | s^n_1) \Big|, \\
	\rm{III} = & \Big| \int_{\cS_2 \backslash A} [h(s^n_3, s_2) - h(s_3, s_2)] \gamma(ds_2 | s^n_1) \Big|.
\end{align*}

Since $h$ is bounded and continuous, the term $\rm{I}$ converges to zero thanks to the weak convergence. As the set $A$ is compact, when $n > N$, we have $\rm{II} \leq \varepsilon$ by the uniform continuity. $\textrm{III} \leq C \varepsilon$ for a generic constant $C$, since $h$ is bounded and $\sup_n \gamma(\cS_2 \backslash A | s^n_1) \leq \varepsilon$. As $\varepsilon > 0$ is arbitrary, we obtain the continuity as desired. 

Claim (b) can be proved similarly.

\end{proof}

\begin{proof}[Proof of Lemma \ref{lem:semi_V}.]
Recalling the $\pi^{t, \gamma}$ in Definition \ref{def:equitrans} but with a generic bicausal $\pi$ after $t+1$, we can show a recursive relationship for the objective $J$ in \eqref{eq:semi_obj}:
\begin{align*}
	J(x_{t}, y_{t}; \pi^{t, \gamma}) = & \Big[ \int \Big(c_{t+1} (x_t, y_t, x_{t+1}, y_{t+1}) +  J (x_{t+1}, y_{t+1}; \pi) \Big) \gamma(dx_{t+1}, dy_{t+1}) \\ 
	& - \int G\Big(x_{t+1}, y_{t+1}, \int h(x_T, y_T) \pi(dx_T, dy_T | x_{t+1}, y_{t+1}) \Big) \gamma(dx_{t+1}, dy_{t+1})  \\
	& + G\Big(x_t, y_t, \int \int h(x_T, y_T) \pi(dx_T, dy_T | x_{t+1}, y_{t+1})  \gamma(dx_{t+1}, dy_{t+1}) \Big)  \\
	& - \int \int \sum^T_{k=t+2}  c_k (x_{t+1}, y_{t+1}, x_k, y_k) \pi(dx_{t+2:T}, dy_{t+2:T} | x_{t+1}, y_{t+1})  \gamma(dx_{t+1}, dy_{t+1})  \\
	& + \int \int  \sum^T_{k=t+2}  c_k (x_t, y_t, x_k, y_k) \pi(dx_{t+2:T}, dy_{t+2:T} | x_{t+1}, y_{t+1}) \gamma(dx_{t+1}, dy_{t+1})  \Big]. 
\end{align*}
The proof is similar to the discrete case. Hence, if an equilibrium transport $\pi^*$ exists and set $\pi = \pi^*$ for $t+1, ..., T-1$, we can derive the extended DP equation \eqref{eq:semi_dpp} in terms of $g_{t+1}$ and $b_k$. Next, we proceed by backward induction to show that the extended DP equation \eqref{eq:semi_dpp} is well-defined and $V_t$, $g_{t+1}$, and $b_k$ are continuous under a finer Polish topology.

Consider time $t = T-1$. Since $\nu(dy_T|y_{T-1})$ is Borel measurable, \citet[Theorem 13.11]{kechris2012classical} shows that there exists a finer Polish topology $\cT^{(1)}_{\cY_{T-1}} \supseteq \cT_{\cY_{T-1}}$ with the same Borel sets $\cB(\cT^{(1)}_{\cY_{T-1}}) = \cB(\cT_{\cY_{T-1}})$, such that $\nu(dy_T|y_{T-1}): (\cY_{T-1}, \cT^{(1)}_{\cY_{T-1}})  \rightarrow \big( \cP(\cY_T), \cV[C_b(\cY_T; \cT_{\cY_T})] \big)$ is continuous. Since previous open sets are still open, $\cT_{\cY_{T-1}}$-continuous functions are still $ \cT^{(1)}_{\cY_{T-1}}$-continuous. We check that the Borel measurable selection theorem \cite[Corollary 1]{brown1973} is applicable to our problem. By Lemma \ref{lem:Cb_conti} and Assumption \ref{a:semi_obj}, the objective at time $T-1$, given by
\begin{align*}
	f(x_{T-1}, y_{T-1}, \gamma) := & \int c_{T} (x_{T-1}, y_{T-1}, x_T, y_T)  \gamma(dx_T, dy_T) \\
	& + G\Big(x_{T-1}, y_{T-1}, \int h (x_T, y_T) \gamma(dx_T, dy_T) \Big),
\end{align*}
is $\cT_{\cX_{T-1}} \times \cT^{(1)}_{\cY_{T-1}} \times \cV[C_b(\cX_{T} \times \cY_{T}; \cT_{\cX_T} \times \cT_{\cY_T})]$-continuous.

By Lemma \ref{lem:closed_corr}, the graph of the correspondence
\begin{align*}
	& D: (\cX_{T-1} \times \cY_{T-1}, \cT_{\cX_{T-1}} \times \cT^{(1)}_{\cY_{T-1}}) \twoheadrightarrow \big(\cP(\cX_{T} \times \cY_{T}), \cV[C_b(\cX_{T} \times \cY_{T}; \cT_{\cX_{T}} \times \cT_{\cY_{T}})]  \big) \\
	&  \text{ that maps } (x_{T-1}, y_{T-1}) \mapsto \Pi(\mu(dx_{T}|x_{T-1}), \nu(dy_{T} | y_{T-1})),
\end{align*}
is closed and thus Borel. Therefore, the objective at time $T-1$ is defined on a Borel set. Moreover, for each $(x_{T-1}, y_{T-1})$, the section $\Pi(\mu(dx_{T}|x_{T-1}), \nu(dy_{T} | y_{T-1}))$ is compact. Hence, by \citet[Corollary 1]{brown1973}, there is a Borel measurable optimizer $\pi^*(dx_T, dy_T | x_{T-1}, y_{T-1})$ such that
\begin{equation*}
	f(x_{T-1}, y_{T-1}, \pi^*(dx_T, dy_T | x_{T-1}, y_{T-1})) = \inf_{\gamma \in \Pi(\mu(dx_{T}|x_{T-1}), \nu(dy_{T} | y_{T-1}))} f(x_{T-1}, y_{T-1}, \gamma),
\end{equation*} 
which is also the value function $V_{T-1}(x_{T-1}, y_{T-1})$. 

As a preparation for applying \citet[Corollary 1]{brown1973} at time $T-2$, we refine the topology on $\cY_{T-1}$ again. Since $\cX_{T-1}$ is finite, we suppose $\cX_{T-1} = \{1, ..., n\}$ without loss of generality. For a given $i \in \cX_{T-1}$, we apply \citet[Theorem 13.11]{kechris2012classical} recursively to
\begin{align*}
	\pi^*(dx_T, dy_T |i,\, y_{T-1}):  (\cY_{T-1}, \cT^{(i)}_{\cY_{T-1}}) \twoheadrightarrow \big( \cP(\cX_{T} \times \cY_{T}), \cV[C_b(\cX_{T} \times \cY_{T}; \cT_{\cX_{T}} \times \cT_{\cY_{T}})] \big). 
\end{align*}
There exists a stronger Polish topology $\cT^{(i+1)}_{\cY_{T-1}} \supseteq \cT^{(i)}_{\cY_{T-1}}$ with $\cB(\cT^{(i+1)}_{\cY_{T-1}}) = \cB(\cT^{(i)}_{\cY_{T-1}}) = \cB(\cT_{\cY_{T-1}})$, such that $\pi^*(dx_T, dy_T |i,\, y_{T-1})$ is $\cT^{(i+1)}_{\cY_{T-1}}$-continuous in $y_{T-1}$. We claim that $\pi^*(dx_T, dy_T |x_{T-1}, y_{T-1})$ is jointly continuous in $(x_{T-1}, y_{T-1})$ under the product topology $\cT_{\cX_{T-1}} \times \cT^{(n+1)}_{\cY_{T-1}}$. Consider a complete compatible metric $d_\cY$ for $\cT^{(n+1)}_{\cY_{T-1}}$ and the discrete metric $d_\cX$ on $\cX_{T-1}$, that is, $d_\cX(x, x') = 0$ if $x=x'$ and  $d_\cX(x, x') = 1$ if $x \neq x'$. $(\cX_{T-1} \times \cY_{T-1}, \cT_{\cX_{T-1}} \times \cT^{(n+1)}_{\cY_{T-1}})$ is a Polish space with the metric $d_\cX(x, x') + d_\cY(y, y')$, where we omitted the time subscript for simplicity. Denote $d_\cP$ as a complete compatible metric for $\big( \cP(\cX_{T} \times \cY_{T}), \cV[C_b(\cX_{T} \times \cY_{T}; \cT_{\cX_{T}} \times \cT_{\cY_{T}})] \big)$. Given $(x, y)$, for every $\varepsilon > 0$, we want to show that there exists $\delta > 0$, such that if $d_\cX(x, x') + d_\cY(y, y') < \delta$, then $d_\cP(\pi^*(\cdot|x, y), \pi^*(\cdot|x', y')) < \varepsilon$. In fact, since for each $i \in \cX_{T-1}$, we can find $\delta_i > 0$, such that if $d_\cY(y, y') < \delta_i$, then  $d_\cP(\pi^*(\cdot|i, y), \pi^*(\cdot|i, y')) < \varepsilon$. Therefore, we can take $\delta = \min\{ \min_i \{ \delta_i \}, 1 \}$, which guarantees $x = x'$ when $d_\cX(x, x') + d_\cY(y, y') < \delta$ and then apply the continuity on $y$ under a fixed $x$ to show $d_\cP(\pi^*(\cdot|x, y), \pi^*(\cdot|x', y')) < \varepsilon$.

Moreover, by Lemma \ref{lem:Cb_conti}, 
\begin{align*}
	g_{T-1} (x_{T-1}, y_{T-1}) & := \int h(x_{T}, y_{T}) \pi^*(dx_{T}, dy_{T} | x_{T-1}, y_{T-1}), \\
	b_{T} (x_{T-1}, y_{T-1}, x_{T-1}, y_{T-1}) & := \int c_T(x_{T-1}, y_{T-1}, x_T, y_T) \pi^*(dx_{T}, dy_{T} | x_{T-1}, y_{T-1}),
\end{align*}
and the value function $V_{T-1}(x_{T-1}, y_{T-1})$ are $\cT_{\cX_{T-1}} \times \cT^{(n+1)}_{\cY_{T-1}}$-continuous. For $i \in \{0, ..., T-2\}$, 
\begin{align*}
	b_{T} (x_i, y_i, x_{T-1}, y_{T-1}) & := \int c_T(x_i, y_i, x_T, y_T) \pi^*(dx_{T}, dy_{T} | x_{T-1}, y_{T-1})
\end{align*}
is $\cT_{\cX_i} \times \cT_{\cY_i} \times \cT_{\cX_{T-1}} \times \cT^{(n+1)}_{\cY_{T-1}}$-continuous.

At time $t=T-2$, $\cY_{T-1}$ is always endowed with $\cT^{(n+1)}_{\cY_{T-1}}$. By \citet[Proposition 7.25]{bertsekas1978stoch} or \citet[Theorem 17.24]{kechris2012classical}, $\big( \cP(\cY_{T-1}), \cV[C_b( \cY_{T-1}; \cT_{\cY_{T-1}})] \big)$ and $\big( \cP(\cY_{T-1}), $ $\cV[C_b( \cY_{T-1}; \cT^{(n+1)}_{\cY_{T-1}})] \big)$ have the same collection of Borel sets, since 
$\cB(\cT^{(n+1)}_{\cY_{T-1}}) = \cB(\cT_{\cY_{T-1}})$. Hence, $\nu(dy_{T-1} | y_{T-2}): (\cY_{T-2}, \cT_{\cY_{T-2}}) \rightarrow \big( \cP(\cY_{T-1}), \cV[C_b( \cY_{T-1}; \cT^{(n+1)}_{\cY_{T-1}})] \big)$ is still Borel measurable. Similarly, by \citet[Theorem 13.11]{kechris2012classical}, there exists a finer Polish topology $\cT^{(1)}_{\cY_{T-2}} \supseteq \cT_{\cY_{T-2}}$ with the same Borel sets $\cB(\cT^{(1)}_{\cY_{T-2}}) = \cB(\cT_{\cY_{T-2}})$, such that $\nu(dy_{T-1}|y_{T-2}): (\cY_{T-2}, \cT^{(1)}_{\cY_{T-2}})  \rightarrow \big( \cP(\cY_{T-1}), \cV[C_b(\cY_{T-1}; \cT^{(n+1)}_{\cY_{T-1}})] \big)$ is continuous.

By Lemma \ref{lem:Cb_conti} and continuity results of $V_{T-1}$, $g_{T-1}$, and $b_T$ above, the objective at time $T-2$, given by
\begin{align*}
	f(x_{T-2}, y_{T-2}, \gamma) := & \int \Big(c_{T-1} (x_{T-2}, y_{T-2}, x_{T-1}, y_{T-1}) +  V_{T-1} (x_{T-1}, y_{T-1}) \Big) \gamma(dx_{T-1}, dy_{T-1}) \\ 
	& - \int G(x_{T-1}, y_{T-1}, g_{T-1} (x_{T-1}, y_{T-1})) \gamma(dx_{T-1}, dy_{T-1})  \\
	& + G\Big(x_{T-2}, y_{T-2}, \int g_{T-1} (x_{T-1}, y_{T-1}) \gamma(dx_{T-1}, dy_{T-1}) \Big) \\
	& - \int b_T (x_{T-1}, y_{T-1}, x_{T-1}, y_{T-1})  \gamma(dx_{T-1}, dy_{T-1})  \\
	& +  \int b_T (x_{T-2}, y_{T-2}, x_{T-1}, y_{T-1})  \gamma(dx_{T-1}, dy_{T-1}), 
\end{align*}
is $\cT_{\cX_{T-2}} \times \cT^{(1)}_{\cY_{T-2}} \times \cV[C_b(\cX_{T-1} \times \cY_{T-1}; \cT_{\cX_{T-1}} \times \cT^{(n+1)}_{\cY_{T-1}})]$-continuous.

By Lemma \ref{lem:closed_corr}, the correspondence
\begin{align*}
	& D: (\cX_{T-2} \times \cY_{T-2}, \cT_{\cX_{T-2}} \times \cT^{(1)}_{\cY_{T-2}}) \twoheadrightarrow \big(\cP(\cX_{T-1} \times \cY_{T-1}), \cV[C_b(\cX_{T-1} \times \cY_{T-1}; \cT_{\cX_{T-1}} \times \cT^{(n+1)}_{\cY_{T-1}} )]  \big) \\
	&  \text{ that maps } (x_{T-2}, y_{T-2}) \mapsto \Pi(\mu(dx_{T-1}|x_{T-2}), \nu(dy_{T-1} | y_{T-2})),
\end{align*}
has a closed graph and is compact-valued. 

Again, \citet[Corollary 1]{brown1973} proves that there exists a Borel measurable optimizer $\pi^*(dx_{T-1}, dy_{T-1} |x_{T-2}, y_{T-2})$ such that
\begin{equation*}
	f(x_{T-2}, y_{T-2}, \pi^*(dx_{T-1}, dy_{T-1} | x_{T-2}, y_{T-2})) = \inf_{\gamma \in \Pi(\mu(dx_{T-1}|x_{T-2}), \nu(dy_{T-1} | y_{T-2}))} f(x_{T-2}, y_{T-2}, \gamma),
\end{equation*} 
which also gives the value function $V_{T-2}(x_{T-2}, y_{T-2})$.

Without loss of generality, suppose $\cX_{T-2} = \{1, ..., n\}$. There exists a finer topology  $\cT^{(n+1)}_{\cY_{T-2}} \supseteq \cT_{\cY_{T-2}}$ with $\cB(\cT^{(n+1)}_{\cY_{T-2}}) = \cB(\cT_{\cY_{T-2}})$, such that $\pi^*(dx_{T-1}, dy_{T-1} |x_{T-2}, y_{T-2})$ is jointly continuous in $(x_{T-2}, y_{T-2})$ under the product topology $\cT_{\cX_{T-2}} \times \cT^{(n+1)}_{\cY_{T-2}}$. Moreover,
\begin{align*}
	g_{T-2} (x_{T-2}, y_{T-2}) & := \int g_{T-1} (x_{T-1}, y_{T-1}) \pi^*(dx_{T-1}, dy_{T-1} | x_{T-2}, y_{T-2}), \\
	b_{T-1} (x_{T-2}, y_{T-2}, x_{T-2}, y_{T-2}) & := \int c_{T-1}(x_{T-2}, y_{T-2}, x_{T-1}, y_{T-1}) \pi^*(dx_{T-1}, dy_{T-1} | x_{T-2}, y_{T-2}), \\
	\int b_T (x_{T-1}, y_{T-1}, x_{T-1}, y_{T-1}) & \pi^*(dx_{T-1}, dy_{T-1} | x_{T-2}, y_{T-2}),  \\
	\int b_T (x_{T-2}, y_{T-2}, x_{T-1}, y_{T-1}) & \pi^*(dx_{T-1}, dy_{T-1} | x_{T-2}, y_{T-2}),
\end{align*}
and the value function $V_{T-2}(x_{T-2}, y_{T-2})$ are $\cT_{\cX_{T-2}} \times \cT^{(n+1)}_{\cY_{T-2}}$-continuous. For $i \in \{0, ..., T-3\}$, 
\begin{align*}
	b_{T-1} (x_i, y_i, x_{T-2}, y_{T-2}) & := \int c_{T-1}(x_i, y_i, x_{T-1}, y_{T-1}) \pi^*(dx_{T-1}, dy_{T-1} | x_{T-2}, y_{T-2}) \\
	\text{ and } \int b_T (x_{i}, y_{i}, x_{T-1}, y_{T-1}) & \pi^*(dx_{T-1}, dy_{T-1} | x_{T-2}, y_{T-2})
\end{align*}
are $\cT_{\cX_i} \times \cT_{\cY_i} \times \cT_{\cX_{T-2}} \times \cT^{(n+1)}_{\cY_{T-2}}$-continuous.

Therefore, we can prove the result by backward induction. $\pi^*(dx_{t+1}, dy_{t+1} | x_{t}, y_{t})$ and $V_t( x_{t}, y_{t})$ are Borel measurable. Indeed, there is a finer Polish topology such that they are continuous.
\end{proof}

The proof of Theorem \ref{thm:lsc} needs the following properties of l.s.c. functions.
\begin{lemma}\label{lem:int_lsc}
	Consider Polish topological spaces $(\cS_1, \cT_{\cS_1})$, $(\cS_2, \cT_{\cS_2})$, and $(\cS_3, \cT_{\cS_3})$. Suppose
	\begin{itemize}
		\item[(1)] the stochastic kernel $\gamma(ds_2| s_1): (\cS_1, \cT_{\cS_1}) \rightarrow (\cP(\cS_2), \cV[C_b(\cS_2; \cT_{\cS_2})])$ is continuous;
		\item[(2)] the function $h(s_3, s_2): (\cS_3 \times \cS_2, \cT_{\cS_3} \times \cT_{\cS_2}) \rightarrow (\R, \cT_\R)$ is  l.s.c. and bounded from below;
		\item[(3)] the function $g(s_2): (\cS_2, \cT_{\cS_2}) \rightarrow (\R, \cT_\R)$ is  l.s.c. The function $G(s_1, g): (\cS_1 \times \R, \cT_{\cS_1} \times \cT_{\R}) \rightarrow (\R, \cT_\R)$ is l.s.c. Moreover, $G(s_1, \cdot)$ is nondecreasing for each $s_1$.
	\end{itemize}
	Then
	\begin{itemize}
		\item[(a)] $\int_{\cS_2} h(s_3, s_2) \gamma(ds_2 | s_1)$ is $\cT_{\cS_3} \times \cT_{\cS_1}$-l.s.c. and bounded from below. If $(\cS_3, \cT_{\cS_3}) = (\cS_1, \cT_{\cS_1})$, it is understood as $\cT_{\cS_1}$-l.s.c. and bounded from below;
		\item[(b)] with $\lambda \in \cP(\cS_2)$, $\int_{\cS_2} h(s_3, s_2) \lambda(ds_2)$ is $\cT_{\cS_3} \times \cV[C_b(\cS_2; \cT_{\cS_2})]$-l.s.c. and bounded from below;
		\item[(c)] $G(s_1, g(s_2)): (\cS_1 \times \cS_2, \cT_{\cS_1} \times \cT_{\cS_2}) \rightarrow (\R, \cT_\R)$ is l.s.c. 
	\end{itemize}	
\end{lemma}

\begin{proof}
By \citet[Lemma 7.14]{bertsekas1978stoch}, there exists a sequence $h_k(s_3, s_2)$ of continuous and bounded functions such that $h_k$ converges increasingly to $h$. By monotone convergence theorem,
\begin{align*}
	\int_{\cS_2} h(s_3, s_2) \gamma(ds_2 | s_1) = \int_{\cS_2} \sup_k h_k (s_3, s_2) \gamma(ds_2 | s_1) = \sup_k \int_{\cS_2} h_k(s_3, s_2) \gamma(ds_2 | s_1).
\end{align*} 
Hence, claim (a) follows since $\int_{\cS_2} h(s_3, s_2) \gamma(ds_2 | s_1)$ is a supremum of continuous functions in $(s_3, s_1)$, by Lemma \ref{lem:Cb_conti}. Claim (b) can be proved similarly.

For claim (c), we consider a sequence $(s^n_1, s^n_2)$ converging to $(s_1, s_2)$ under the product topology $\cT_{\cS_1} \times \cT_{\cS_2}$. Then
\begin{align*}
	\liminf_{n \rightarrow \infty} G(s^n_1, g(s^n_2)) \geq \liminf_{n \rightarrow \infty} G(s^n_1, \inf_{k\geq n} g(s^k_2)) \geq G(s_1, \lim_{n\rightarrow \infty} \inf_{k \geq n} g(s^k_2)) \geq G(s_1, g(s_2)).
\end{align*}
The first inequality uses the fact that $g(s^n_2) \geq \inf_{k\geq n} g(s^k_2)$ and $G(s^n_1, \cdot)$ is nondecreasing. The second inequality holds since $G(\cdot, \cdot)$ is l.s.c. The last inequality is due to $\lim_{n\rightarrow \infty} \inf_{k \geq n} g(s^k_2) = \liminf_{n \rightarrow \infty} g(s^n_2) \geq g(s_2)$.  

\end{proof}

\begin{proof}[Proof of Theorem \ref{thm:lsc}.]
The proof is a modification for the case of Lemma \ref{lem:semi_V}.

Consider time $t = T-1$. $\nu(dy_T|y_{T-1}): (\cY_{T-1}, \cT^{(1)}_{\cY_{T-1}})  \rightarrow \big( \cP(\cY_T), \cV[C_b(\cY_T; \cT_{\cY_T})] \big)$ is continuous. By Lemma \ref{lem:int_lsc} and Assumption \ref{a:semi_lsc}, the objective at time $T-1$, given by
\begin{align*}
	f(x_{T-1}, y_{T-1}, \gamma) := & \int c_{T} (x_{T-1}, y_{T-1}, x_T, y_T)  \gamma(dx_T, dy_T) + G\Big(x_{T-1}, y_{T-1}, \int h (x_T, y_T) \gamma(dx_T, dy_T) \Big),
\end{align*}
is $\cT_{\cX_{T-1}} \times \cT^{(1)}_{\cY_{T-1}} \times \cV[C_b(\cX_{T} \times \cY_{T}; \cT_{\cX_T} \times \cT_{\cY_T})]$-l.s.c. Indeed, $\int c_{T} (x_{T-1}, y_{T-1}, x_T, y_T)  \gamma(dx_T, dy_T)$ is l.s.c. in $(x_{T-1}, y_{T-1}, \gamma)$ by Lemma \ref{lem:int_lsc} (b). Similarly, $\int h (x_T, y_T) \gamma(dx_T, dy_T)$ is l.s.c. in $\gamma$. $G(x_{T-1}, y_{T-1}, \int h (x_T, y_T) \gamma(dx_T, dy_T) )$ is l.s.c. in $(x_{T-1}, y_{T-1}, \gamma)$ by Lemma \ref{lem:int_lsc} (c).

Hence, by \citet[Corollary 1]{brown1973}, there is a Borel measurable optimizer $\pi^*(dx_T, $ $ dy_T | x_{T-1}, y_{T-1})$. It is jointly continuous in $(x_{T-1}, y_{T-1})$ under a finer product topology $\cT_{\cX_{T-1}} \times \cT^{(n+1)}_{\cY_{T-1}}$ by \citet[Theorem 13.11]{kechris2012classical}. We note that l.s.c. functions are still l.s.c. under a finer topology. Then $g_{T-1}(x_{T-1}, y_{T-1})$ and $	b_{T} (x_i, y_i, x_{T-1}, y_{T-1})$ are still l.s.c. when $\cX_{T-1} \times \cY_{T-1}$ is endowed with $\cT_{\cX_{T-1}} \times \cT^{(n+1)}_{\cY_{T-1}}$.

At time $T-2$, the objective is given by
\begin{align*}
	f(x_{T-2}, y_{T-2}, \gamma) := & \int c_{T-1} (x_{T-2}, y_{T-2}, x_{T-1}, y_{T-1})  \gamma(dx_{T-1}, dy_{T-1}) \\ 
	& + G\Big(x_{T-2}, y_{T-2}, \int g_{T-1} (x_{T-1}, y_{T-1}) \gamma(dx_{T-1}, dy_{T-1}) \Big) \\
	& +  \int b_T (x_{T-2}, y_{T-2}, x_{T-1}, y_{T-1})  \gamma(dx_{T-1}, dy_{T-1}).
\end{align*}
It is $\cT_{\cX_{T-2}} \times \cT^{(1)}_{\cY_{T-2}} \times \cV[C_b(\cX_{T-1} \times \cY_{T-1}; \cT_{\cX_{T-1}} \times \cT^{(n+1)}_{\cY_{T-1}})]$-l.s.c., by Lemma \ref{lem:int_lsc} and Assumption \ref{a:semi_lsc}. The remaining proof follows similarly as in Lemma \ref{lem:semi_V}.
\end{proof}

\subsection{The continuous and non-Markovian case}\label{sec:proof_params}

Lemma \ref{lem:metric_cont} extends Lemma \ref{lem:Cb_conti} to unbounded functions. We note that Lemma \ref{lem:metric_cont} relies on the metric used to define the growth rate.
\begin{lemma}\label{lem:metric_cont}
	Consider metric spaces $(\cS_1, d_1)$, $(\cS_2, d_2)$, and $(\cS_3, d_3)$ that are Polish. Moreover, $\cS_2$ is a finite dimensional vector space. With a given $p \in [1, \infty)$, suppose
	\begin{itemize}
		\item[(1)] the stochastic kernel $\gamma(ds_2| s_1): (\cS_1, d_1) \rightarrow \cP_p(\cS_2)$ is continuous;
		\item[(2)] the function $h(s_3, s_2) \in C_p(\cS_3 \times \cS_2)$.
	\end{itemize}
	Then
	\begin{itemize}
		\item[(a)] $\int_{\cS_2} h(s_3, s_2) \gamma(ds_2 | s_1)$ is continuous in $(s_3, s_1)$;
		\item[(b)] if $(\cS_3, d_3) = (\cS_1, d_1)$ and $\gamma(ds_2 | s_1)$ satisfies $\int_{\cS_2} d_2(s_2, \bar{s}_2)^p \gamma(ds_2|s_1) \leq C(1 + d_1(s_1, \bar{s}_1)^p)$ for all $s_1 \in \cS_1$, then $\int_{\cS_2} h(s_1, s_2) \gamma(ds_2 | s_1) \in C_p (\cS_1)$;
		\item[(c)] with $\lambda \in \cP_p(\cS_2)$, $\int_{\cS_2} h(s_3, s_2) \lambda(ds_2)$ is continuous in $(s_3, \lambda)$. 
	\end{itemize}
\end{lemma}
\begin{proof}
For the claim (a),  consider a sequence $(s^n_3, s^n_1)$ converging to $(s_3, s_1)$. $B := \{(s^n_3, s^n_1)\}^\infty_{n=1} \cup \{ (s_3, s_1) \}$ is a compact set. Fix an arbitrary $\varepsilon > 0$. Similar to Lemma \ref{lem:Cb_conti}, we can find a compact subset $A \subset \cS_2$ such that $\sup_n \gamma(  \cS_2 \backslash A | s^n_1) \leq \varepsilon$. Without loss of generality, we can assume that there is a sufficiently large radius $R$ such that the ball with radius $R$ is contained by $A$:
\begin{equation}\label{eq:ball}
	K_R : = \{ s_2 | d_2(s_2, \bar{s}_2) \leq R \} \subset A.
\end{equation}
Indeed, since $\cS_2$ is finite dimensional, the closed ball $K_R$ is compact \cite[Theorem 5.26]{charalambos2013infinite}. We can use $K_R \cup A$ to replace $A$ if needed.

Similarly, we have 
\begin{align*}
	& \Big| \int_{\cS_2} h(s^n_3, s_2) \gamma(ds_2 | s^n_1) - \int_{\cS_2} h(s_3, s_2) \gamma(ds_2 | s_1) \Big| \leq \rm{I} + \rm{II} + \rm{III},
\end{align*}
where
\begin{align*}
	\rm{I} = & \Big| \int_{\cS_2} h(s_3, s_2) [\gamma(ds_2 | s^n_1) - \gamma(ds_2 | s_1)] \Big|, \\
	\rm{II} = & \Big| \int_A [h(s^n_3, s_2) - h(s_3, s_2)] \gamma(ds_2 | s^n_1) \Big|, \\
	\rm{III} = & \Big| \int_{\cS_2 \backslash A} [h(s^n_3, s_2) - h(s_3, s_2)] \gamma(ds_2 | s^n_1) \Big|.
\end{align*}

Since $h(s_3, \cdot)$ satisfies the growth rate condition, the term $\rm{I}$ converges to zero when $n \rightarrow \infty$ by \citet[Definition 6.7 (iv)]{villani2009optimal}. As the set $B \times A$ is compact, when $n > N$, we have $\rm{II} \leq \varepsilon$ by the uniform continuity. For \rm{III}, since $A$ contains a large enough ball with radius $R$,
\begin{align*}
	\rm{III} \leq & \int_{\cS_2 \backslash A} \left( |h(s^n_3, s_2)| + |h(s_3, s_2)| \right) \gamma(ds_2 | s^n_1) \\		
	\leq & C \sup_n \left[ 1 + d_3(s^n_3, \bar{s}_3)^p \right] \varepsilon + C \left[ 1 + d_3(s_3, \bar{s}_3)^p \right] \varepsilon + 2 C \int_{\cS_2 \backslash K_R} d_2(s_2, \bar{s}_2)^p \gamma(ds_2 | s^n_1).
\end{align*} 
For the second inequality, we have used the growth rate condition of $h$, the inequality $\sup_n \gamma(  \cS_2 \backslash A | s^n_1) $ $\leq \varepsilon$, and the fact that the ball $K_R \subset A$. By \citet[Definition 6.7 (iii)]{villani2009optimal}, the last term is less than $2 C \varepsilon$ when $n \rightarrow \infty$ and $R$ is sufficiently large. In summary,
\begin{equation*}
	\limsup_{n \rightarrow \infty} {\rm I} + {\rm II} + {\rm III} \leq C \varepsilon,
\end{equation*}
with a generic constant that is independent of $\varepsilon$. As $\varepsilon > 0$ is arbitrary, we obtain the continuity as desired.

Claim (b) is direct. Claim (c) can be proved similarly to (a).
\end{proof}


Another useful result is the concatenation of continuous kernels is still continuous.
\begin{lemma}\label{lem:concat}
	Consider metric spaces $(\cS_1, d_1)$, $(\cS_2, d_2)$, and $(\cS_3, d_3)$ that are Polish. Moreover, $\cS_2$ and $\cS_3$ are finite dimensional vector spaces.  With a given $p \in [1, \infty)$, suppose
	\begin{itemize}
		\item[(1)] the stochastic kernel $\pi(ds_2| s_1): (\cS_1, d_1) \rightarrow \cP_p(\cS_2)$ is continuous and $\int_{\cS_2} d_2(s_2, \bar{s}_2)^p \pi(ds_2|s_1) \leq C(1 + d_1(s_1, \bar{s}_1)^p)$;
		\item[(2)]  the stochastic kernel $\pi(ds_3 |s_1, s_2): (\cS_1 \times \cS_2, (d^p_1 + d^p_2)^{1/p}) \rightarrow \cP_p(\cS_3)$ is continuous and $\int_{\cS_3} d_3(s_3, \bar{s}_3)^p \pi(ds_3|s_1, s_2) \leq C(1 + d_1(s_1, \bar{s}_1)^p + d_2(s_2, \bar{s}_2)^p)$.
	\end{itemize}
	Then there exists a unique kernel
	\begin{equation}\label{eq:concat}
		\pi(ds_{2:3} | s_1) =  \pi(ds_3 | s_1, s_2) \pi(ds_2 | s_1)
	\end{equation}
	from $(\cS_1, d_1)$ to $\cP_p(\cS_2 \times \cS_3)$. It is continuous in $s_1$.
\end{lemma}
\begin{proof}
By \citet[Proposition 7.28]{bertsekas1978stoch}, there exists a unique $\pi(ds_{2:3} | s_1)$ defined by \eqref{eq:concat}.  We only need to prove the continuity.

Consider a sequence $\{ s^n_1 \}$ converging to $s_1$. By \citet[Definition 6.7 and Theorem 6.8]{villani2009optimal}, we need to show that for all continuous functions $\varphi$ with a growth rate of $|\varphi(s_2, s_3)| \leq C(1 + d_2(s_2, \bar{s}_2)^p + d_3(s_3, \bar{s}_3)^p)$, one has
\begin{align}
	& \lim_{n \rightarrow \infty} \int_{\cS_2 \times \cS_3} \varphi(s_2, s_3) \pi(ds_{2:3} | s^n_1) = \int_{\cS_2 \times \cS_3} \varphi(s_2, s_3) \pi(ds_{2:3} | s_1). \label{eq:conver}
\end{align}
By Lemma \ref{lem:metric_cont} (b), $\Phi(s^n_1, s_2) := \int_{\cS_3} \varphi(s_2, s_3) \pi(ds_3| s^n_1, s_2)$ is in $C_p (\cS_1 \times  \cS_2)$. Thus, we can apply Lemma \ref{lem:metric_cont} (b) again to $\int_{\cS_2} \Phi(s^n_1, s_2) \pi(ds_2|s^n_1)$. It belongs to $C_p (\cS_1)$. Then the convergence \eqref{eq:conver} holds.

\end{proof}

\begin{proof}[Proof of Theorem \ref{thm:params}.]
The recursive relationship for $V_t(x_{1:t}, y_{1:t})$ in \eqref{eq:params_dpp} is direct. We prove that $g_{t+1}$ and $b_{t+1}$ are continuous and $\theta^*_{t+1}(x_{1:t}, y_{1:t})$ exists.

First, we have $V_T \in C_p(\cX \times \cY)$ by the boundary condition and Assumption \ref{a:params_obj}. 

Next, consider time $t = T-1$. Thanks to Assumption \ref{a:params}(1), Assumption \ref{a:params_obj} and Lemma \ref{lem:metric_cont}, the objective functional
\begin{align*}
	f(x_{1:T-1}, y_{1:T-1}, \theta_{T}) = & G\Big(x_{T-1}, y_{T-1}, \int h (x_{1:T}, y_{1:T}) \gamma(dx_T, dy_T| \theta_T) \Big) \\
	& + \int c (x_{T-1}, y_{T-1}, x_{1:T}, y_{1:T}) \gamma(dx_T, dy_T | \theta_T)
\end{align*}
is $\cT_{\cX_{1:T-1}} \times \cT_{\cY_{1:T-1}} \times \cT_{\Theta_{T}}$-continuous. Moreover, by Assumption \ref{a:params}(3), $f(x_{1:T-1}, y_{1:T-1}, \theta_{T})$ is strictly quasiconvex in $\theta_T$. Assumption \ref{a:params}(2) imposes the required properties of the correspondence $D_{T-1}$. By a version of Berge's maximum theorem with strict quasiconcavity, see \citet[Theorem 9.14 and Corollary 9.20]{sundaram1996first} or \citet[Theorem 17.31]{charalambos2013infinite}, $V_{T-1}$ given by
\begin{equation*}
	V_{T-1}(x_{1:T-1}, y_{1:T-1})  = \inf_{ \theta_{T} \in D_{T-1}(x_{1:T-1}, y_{1:T-1})}  	f(x_{1:T-1}, y_{1:T-1}, \theta_{T}) 
\end{equation*}
is $\cT_{\cX_{1:T-1}} \times \cT_{\cY_{1:T-1}}$-continuous and the infimum is attained. There exists a unique and continuous optimizer $\theta^*_T(x_{1:T-1}, y_{1:T-1}): (\cX_{1:T-1} \times \cY_{1:T-1}, d) \rightarrow (\Theta_{T}, d_{\Theta_{T}})$, which is a function instead of correspondence. As a composition, $\gamma(dx_{T}, dy_{T}| \theta^*_{T}(x_{1:T-1}, y_{1:T-1})): (\cX_{1:T-1} \times \cY_{1:T-1}, d) \rightarrow \cP_p(\cX_{T} \times \cY_T)$ is also continuous.

For the growth rate, we note the marginal constraint $$\gamma(dx_{T}, dy_{T}| \theta^*_{T}(x_{1:T-1}, y_{1:T-1})) \in \Pi(\mu(dx_T | x_{1:T-1}), \nu(dy_T | y_{1:T-1}) )$$ together with Assumption \ref{a:measure} (2). Then
\begin{align*}
	& \left| \int_{\cX_T \times \cY_T} c(x_{T-1}, y_{T-1}, x_{1:T}, y_{1:T}) \gamma(dx_{T}, dy_{T}| \theta^*_{T}(x_{1:T-1}, y_{1:T-1})) \right| \\
	&\qquad \leq \int_{\cX_T \times \cY_T} C \left[ 1 + \sum^T_{t=1} d_{\cX_t}(x_t, \bar{x}_t)^p + \sum^T_{t=1} d_{\cY_t}(y_t, \bar{y}_t)^p \right] \gamma(dx_{T}, dy_{T}| \theta^*_{T}(x_{1:T-1}, y_{1:T-1})) \\
	& \qquad = C + C \int_{\cX_T} \sum^T_{t=1} d_{\cX_t}(x_t, \bar{x}_t)^p \mu(dx_T | x_{1:T-1}) + C \int_{\cY_T} \sum^T_{t=1} d_{\cY_t}(y_t, \bar{y}_t)^p \nu(dy_T | y_{1:T-1}) \\
	& \qquad  \leq C \left[ 1 + \sum^{T-1}_{t=1} d_{\cX_t}(x_t, \bar{x}_t)^p + \sum^{T-1}_{t=1} d_{\cY_t}(y_t, \bar{y}_t)^p \right].
\end{align*}
Similarly, with Assumption \ref{a:params_obj}, one has
\begin{align}
	& \left| G \left( x_{T-1}, y_{T-1}, \int_{\cX_T \times \cY_T} h (x_{1:T}, y_{1:T}) \gamma(dx_{T}, dy_{T}| \theta^*_{T}(x_{1:T-1}, y_{1:T-1})) \right) \right| \label{eq:a}\\
	& \leq C \left(  1  + \left| \int_{\cX_T \times \cY_T} h (x_{1:T}, y_{1:T}) \gamma(dx_{T}, dy_{T}| \theta^*_{T}(x_{1:T-1}, y_{1:T-1})) \right|^{pr} + d((x_{T-1}, y_{T-1}), (\bar{x}_{T-1}, \bar{y}_{T-1}))^p \right) \nonumber \\
	& \leq C \left(  1  +  \int_{\cX_T \times \cY_T} \left| h (x_{1:T}, y_{1:T})  \right|^{pr} \gamma(dx_{T}, dy_{T}| \theta^*_{T}(x_{1:T-1}, y_{1:T-1})) + d((x_{T-1}, y_{T-1}), (\bar{x}_{T-1}, \bar{y}_{T-1}))^p \right)  \nonumber \\
	& \leq C[ 1 + d((x_{1:T-1}, y_{1:T-1}), (\bar{x}_{1:T-1}, \bar{y}_{1:T-1}))^{p}]. \nonumber
\end{align}
Therefore, $V_{T-1} \in C_p (\cX_{1:T-1} \times \cY_{1:T-1})$.

By Lemma \ref{lem:metric_cont} and $1/r \leq p$, the following functions are continuous in $(x_{1:T-1}, y_{1:T-1})$:
\begin{align*}
	g_{T-1} (x_{1:T-1}, y_{1:T-1}) &:= \int h(x_{1:T}, y_{1:T}) \gamma(dx_{T}, dy_{T}| \theta^*_{T}(x_{1:T-1}, y_{1:T-1})), \\
	b_{T-1} (x_{i}, y_{i}, x_{1:T-1}, y_{1:T-1}) & := \int c(x_{i}, y_{i}, x_{1:T}, y_{1:T}) \gamma(dx_{T}, dy_{T}| \theta^*_{T}(x_{1:T-1}, y_{1:T-1})), \; i \in \{0, ..., T-1\}.
\end{align*}

In each step of backward induction, we also use Lemma \ref{lem:concat} that the successive concatenation of continuous conditional kernels $ \gamma(dx_T, dy_T | \theta^*_{T}(x_{1:T-1}, y_{1:T-1})), \ldots, \gamma(dx_{t+1}, dy_{t+1} | \theta^*_{t+1}(x_{1:t}, y_{1:t}))$ induces a unique continuous conditional probability measure $\gamma(dx_{t+1:T}, dy_{t+1:T} | \theta^*_{t+1:T}(x_{1:t}, y_{1:t}))$. 
\end{proof}

\begin{proof}[Proof of Corollary \ref{cor:classical_nonMark}.]
With $\Theta_{t+1} = \cP_p(\cX_{t+1} \times \cY_{t+1})$, we only need to show the correspondence $D_t$ satisfies Assumption \ref{a:params} (2) in Lemma \ref{lem:corres}. Then Corollary \ref{cor:classical_nonMark} follows directly from Theorem \ref{thm:params}. Compared with \cite{neufeld2021stability}, Lemma \ref{lem:corres} considers the $W_p$ metric and does not have martingale constraints.
\end{proof} 

Denote the sum of Wasserstein distances between $\mu, \mu' \in \cP(\cX)$ and  $\nu, \nu' \in \cP(\cY)$ as
\begin{align*}
	W^\oplus_p ((\mu, \nu), (\mu', \nu')) := W_p (\mu, \mu') + W_p (\nu, \nu').  
\end{align*}

\begin{lemma}\label{lem:corres}
	Consider the metric spaces and Wasserstein distances in Section \ref{sec:metrics}. For a given $t$, denote a correspondence as
	\begin{align*}
		D: \cX_{1:t} \times \cY_{1:t} \twoheadrightarrow \cP_p(\cX_{t+1} \times \cY_{t+1}) \text{ that maps } (x_{1:t}, y_{1:t}) \mapsto \Pi(\mu(dx_{t+1}|x_{1:t}), \nu(dy_{t+1} | y_{1:t})).
	\end{align*}
	Suppose Assumption \ref{a:measure} Condition (1) holds. Then $D$ is a continuous correspondence and $D(x_{1:t}, y_{1:t})$ is non-empty, convex, and compact, under the product topology of $\cT_{\cX_{1:t}}$, $\cT_{\cY_{1:t}}$, and the topology induced by the Wasserstein distance on $\cP_p(\cX_{t+1} \times \cY_{t+1})$.
\end{lemma}

\begin{proof}
With the metric $d$ on $\cX_{1:t} \times \cY_{1:t}$ and the metric $W_p$ on $\cP_p(\cX_{t+1})$, $\cP_p(\cY_{t+1})$, and $\cP_p(\cX_{t+1} \times \cY_{t+1})$, we first show that the correspondence $(\alpha, \beta) \mapsto \Pi(\alpha, \beta)$ is continuous with $W^\oplus_p$ metric on the domain and $W_p$ on the range. 

{\bf Compactness and upper hemicontinuity}: The idea is similar to Lemma \ref{lem:closed_corr}. We show that if a sequence $\{(\alpha^n, \beta^n, \gamma^n)\}$ is in the graph of $D$ and
\begin{equation}\label{eq:conv-range}
	\lim_{n \rightarrow \infty} W^\oplus_p ((\alpha^n, \beta^n), (\alpha, \beta)) = 0,
\end{equation}
then the sequence $\{\gamma^n\}$ has a limit point in $\Pi(\alpha, \beta)$. 

By \citet[Definition 6.7]{villani2009optimal}, the convergence in \eqref{eq:conv-range} implies the usual weak convergence. With the same argument in Lemma \ref{lem:closed_corr}, we can prove that a subsequence $\{\gamma^{n_k}\}^\infty_{k=1}$ converges weakly in the usual sense to some $\gamma \in \Pi(\alpha, \beta)$. Furthermore,
\begin{align*}
	& \lim_{k \rightarrow \infty} \int_{\cX_{t+1} \times \cY_{t+1}} d((x_{t+1} , y_{t+1}), (\bar{x}_{t+1}, \bar{y}_{t+1}))^p \gamma^{n_k} (dx_{t+1}, dy_{t+1}) \\
	& \qquad = \lim_{k \rightarrow \infty} \int_{\cX_{t+1} \times \cY_{t+1}}  \left[ d_{\cX_{t+1}} (x_{t+1} , \bar{x}_{t+1})^p + d_{\cY_{t+1}}(y_{t+1}, \bar{y}_{t+1})^p \right] \gamma^{n_k} (dx_{t+1}, dy_{t+1}) \\
	& \qquad = \lim_{k \rightarrow \infty} \left( \int_{\cX_{t+1}} d_{\cX_{t+1}} (x_{t+1} , \bar{x}_{t+1})^p \alpha^{n_k} (dx_{t+1}) +  \int_{\cY_{t+1}} d_{\cY_{t+1}}(y_{t+1}, \bar{y}_{t+1})^p \beta^{n_k} (dy_{t+1}) \right) \\
	& \qquad = \int_{\cX_{t+1}} d_{\cX_{t+1}} (x_{t+1} , \bar{x}_{t+1})^p \alpha (dx_{t+1}) +  \int_{\cY_{t+1}} d_{\cY_{t+1}}(y_{t+1}, \bar{y}_{t+1})^p \beta (dy_{t+1}) \\
	& \qquad =  \int_{\cX_{t+1} \times \cY_{t+1}}  \left[ d_{\cX_{t+1}} (x_{t+1} , \bar{x}_{t+1})^p + d_{\cY_{t+1}}(y_{t+1}, \bar{y}_{t+1})^p \right] \gamma (dx_{t+1}, dy_{t+1}) \\
	& \qquad = \int_{\cX_{t+1} \times \cY_{t+1}} d((x_{t+1} , y_{t+1}), (\bar{x}_{t+1}, \bar{y}_{t+1}))^p \gamma (dx_{t+1}, dy_{t+1}).
\end{align*} 
Therefore, conditions in \citet[Definition 6.7 (i)]{villani2009optimal} are verified and we obtain the convergence $\lim_{k \rightarrow \infty} W_p(\gamma^{n_k}, \gamma) = 0$.

{\bf Lower hemicontinuity}: By \citet[Theorem 17.21]{charalambos2013infinite}, we have to prove the following claim: For any $(\alpha, \beta) \in \cP_p(\cX_{t+1}) \times \cP_p(\cY_{t+1})$, if $(\alpha^n, \beta^n)$ converges weakly to $(\alpha, \beta)$ in the metric $W^\oplus_p$, then for each $\gamma \in \Pi(\alpha, \beta)$, there exists a subsequence $(\alpha^{n_k}, \beta^{n_k})$ and $\gamma^k \in \Pi(\alpha^{n_k}, \beta^{n_k})$ for each $k$, such that $\gamma^k$ converges weakly to $\gamma$ in $\cP_p(\cX_{t+1} \times \cY_{t+1})$ \citep[Definition 6.7]{villani2009optimal}. Indeed, this result has been proved in \citet[Equations 2.1 and 2.2]{beiglbock2022approx}.

We have proved the continuity of the correspondence $(\alpha, \beta) \mapsto \Pi(\alpha, \beta)$ and the compactness of $\Pi(\alpha, \beta)$. Moreover, Assumption \ref{a:measure} Condition (1) guarantees that $$(x_{1:t}, y_{1:t}) \mapsto (\mu(dx_{t+1}|x_{1:t}), \nu(dy_{t+1} | y_{1:t}))$$ is continuous. By \citet[Theorem 17.23]{charalambos2013infinite} on the continuity of the composition of correspondence, $D$, viewed as the composition 
\begin{align*}
	(x_{1:t}, y_{1:t}) \mapsto (\mu(dx_{t+1}|x_{1:t}), \nu(dy_{t+1} | y_{1:t})) \mapsto \Pi(\mu(dx_{t+1}|x_{1:t}), \nu(dy_{t+1} | y_{1:t})),
\end{align*}
is also continuous.
\end{proof}


\section{Supplement to the executive job market}\label{sec:supp_exec}
\subsection{Data cleaning and summary statistics}\label{sec:clean}
We focus on a five-year time horizon 2017 -- 2021. This choice is motivated by the following factors. A longer time horizon incurs a heavy computational burden. Moreover, in our final dataset, the median tenure of a CEO is 4 years and the average tenure is 5.97 years. In the literature with different datasets, \citet[Table 1]{taylor2013ceo} reported the average CEO tenure as 7.9 years and the median as 6 years. Overall, a five-year time horizon is close to the CEO's tenure. Besides, the job market can change significantly over a longer period.

To obtain and clean the data, we adopt the following steps:
\begin{enumerate}[label={(\arabic*)}]
	\item We download the net sales data from Compustat and the executive compensation information from Execucomp. There are no restrictions on the net sales or market values of firms at this stage. The key variables are net sales (``sale" in Compustat) and total compensation (``tdc1" in Execucomp). In 2017 -- 2021, we have data on 1998 firms from Execucomp and 11911 firms from Compustat. Execucomp primarily includes salary information for S\&P 1500 components.
	\item We remove firms with missing net sales data for the five-year period. Compensation data are usually available for the CEO, chief financial officer, and the three other most highly compensated executive officers. We use their mean salary as the representative wage paid to top executives by the firm. After merging the sale and wage datasets, there are 1590 firms remaining.
	\item To further filter the firms, we restrict our analysis to companies with investment-grade credit ratings. We download the S\&P Domestic Long-Term Issuer Credit Rating (variable name: ``splticrm") and remove firms with a rating of CCC+ or lower. We also remove firms without ratings.
\end{enumerate}
Table \ref{tab:data} shows the descriptive statistics for the data. Our final sample consists of 790 firms, with a total market value of 61.7 trillion U.S. dollars in 2021. 396 firms are S\&P 500 components. 

\begin{table}[H]
	\centering
	\begin{tabular}{cccccc}
		\hline
		Industry & GICS Code & Number of firms & Market value &  Net sales & Compensation \\
		\hline
		Energy & 10  & 48   & 22693.93 & 19713.07  & 4.82  \\
		\hline
		Materials & 15 & 66  &	12889.25 &	 7301.68	& 3.79 \\
		\hline
		Industrials & 20  & 118 & 20781.59	& 12368.6	& 4.14  \\
		\hline
		Consumer   &  \multirow{2}{*}{25}  & \multirow{2}{*}{112} &	\multirow{2}{*}{32048.23}  & \multirow{2}{*}{16761.5}	    &  \multirow{2}{*}{5.80}	   \\
		Discretionary &    &     &    &      & \\
		\hline
		Consumer & \multirow{2}{*}{30}  & \multirow{2}{*}{48} &	\multirow{2}{*}{46617.8}	& \multirow{2}{*}{34468.41}	&  \multirow{2}{*}{5.55} \\
		Staples &       &      &        &             &       \\
		\hline
		Health Care & 35  & 63 & 51570.67 &	33801.81 & 6.07  \\
		\hline
		Financials & 40 & 110 &	24242.18 &	16041.93 & 5.38 \\
		\hline
		Information  & \multirow{2}{*}{45} & \multirow{2}{*}{79} &	\multirow{2}{*}{75539.2}	& \multirow{2}{*}{17750.21}	 &  \multirow{2}{*}{7.25} \\
		Technology & & & & & \\
		\hline
		Communication & \multirow{2}{*}{50} & \multirow{2}{*}{30} &	\multirow{2}{*}{78888.99}	& \multirow{2}{*}{26807.29}	& \multirow{2}{*}{9.02} \\
		Services & & & & & \\
		\hline
		Utilities & 55 & 50 & 18192.87 & 7446.87 & 3.21 \\
		\hline
		Real Estate &  60 & 66 & 11529.68 & 1949.91 & 3.88 \\
		\hline 
	\end{tabular}
	\caption{Summary statistics of firms over 2017 -- 2021. Values in the last three columns are measured in millions of USD.}\label{tab:data}
\end{table}

\subsection{Estimation methodology}\label{sec:method}
	\subsubsection{Number of clusters}\label{sec:num_group}
	Ranks are discrete data that are suitable for discrete OT. However, the smallest number of firms in a single industry, as shown in Table \ref{tab:data}, is 30. If we use the original ranks as the variables $x_t$ and $y_t$, a large amount of data would be needed to estimate the transition probability matrices $\mu(dx_{t+1} | x_t)$ and $\nu(dy_{t+1} | y_t)$. This is unrealistic, as the earliest wage data in Execucomp is primarily from 1992. In addition, transitions of original ranks can be noisy and sensitive to fluctuations in the data. There is also a practical concern about the computational burden of a large transport plan matrix.
	
	Therefore, we aggregate data into several clusters with orders. There are two questions to consider in this process: the number of clusters and whether to use even-sized or uneven-sized clusters. As for the first question, the number of clusters cannot be too small or too large. A small number of clusters fails to capture the variations in the data. If we use only one cluster, all industries would perfectly match firms with wages. On the other hand, using too many clusters leads to the problems discussed above. There is no theoretical result on the optimal number of clusters. \citet[Definition 1.2]{backhoff2020estimating} suggests using $N^{1/(T+1)}$ clusters, where $N$ is the number of time series. However, this result is only valid for $N \rightarrow \infty$, and would be too small when there are only a few data points. Given these difficulties, we suggest a practical and data-driven rule for determining the number of clusters.
	
	For simplicity, we will use the same number of clusters for each industry to ensure fairness when comparing cross-sectional differences. The main feature we want to preserve after clustering is the relationship between sales and wages. Clusters are ordered and sales/wages in the same cluster are assigned with the same rank. Denote the sale and wage cluster ranks for the firm $i$ at time $t$ as $x^i_t$ and $y^i_t$, respectively. We call the mean absolute value of the difference in wage and sale ranks, that is,
	$$\frac{1}{NT}\sum^N_{i=1}\sum^T_{t=1}|x^i_t - y^i_t|,$$
	the {\it sale-wage discrepancy} for this industry. If the number of clusters is chosen appropriately, the sale-wage discrepancy should be bigger for sectors with lower sale-wage efficiency reported in Table \ref{tab:corr}. In other words, if we calculate the correlation between the sale-wage discrepancy and the sale-wage efficiency, we should prefer the number of clusters that generates strong negative correlations.
	
\begin{figure}
		\centering
		\begin{minipage}{0.48\textwidth}
				\centering
				\includegraphics[width=0.95\textwidth]{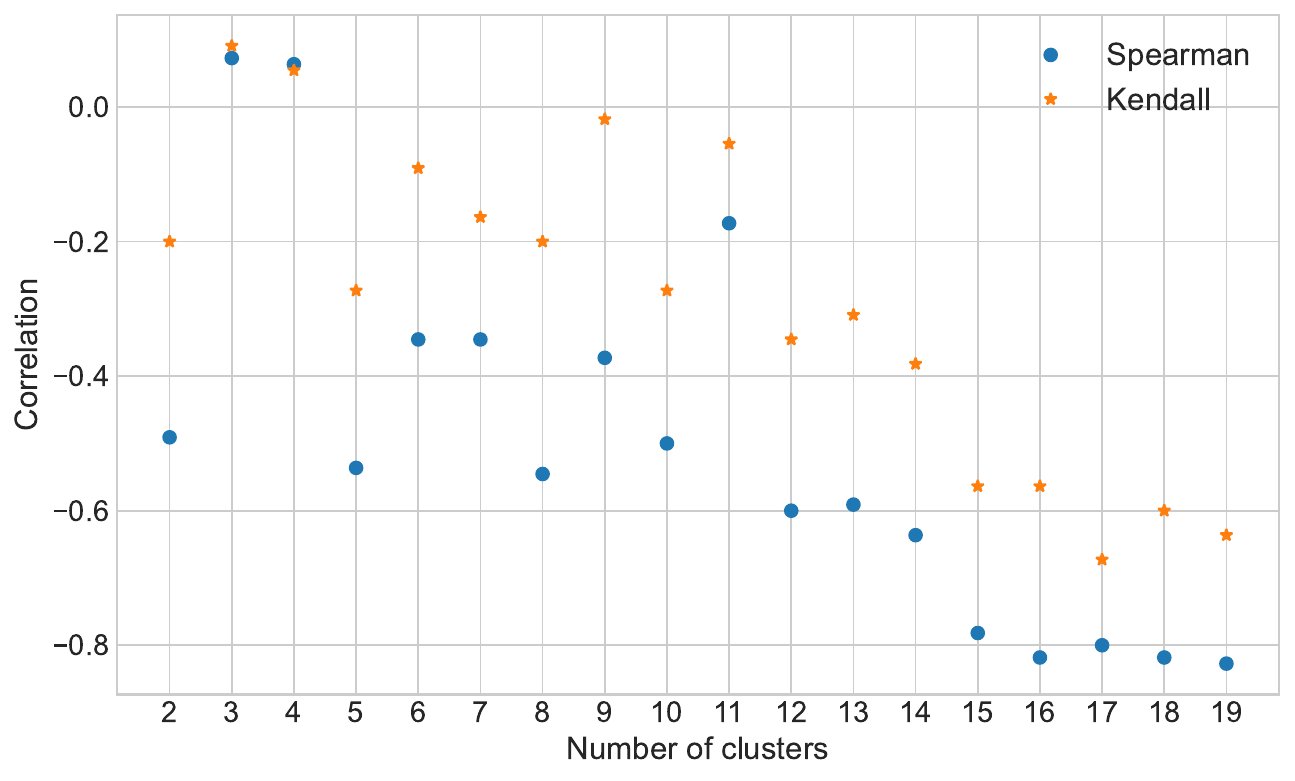}
				\subcaption{Uneven-sized clusters}\label{fig:uneven}
			\end{minipage}
			\begin{minipage}{0.48\textwidth}
				\centering
				\includegraphics[width=0.95\textwidth]{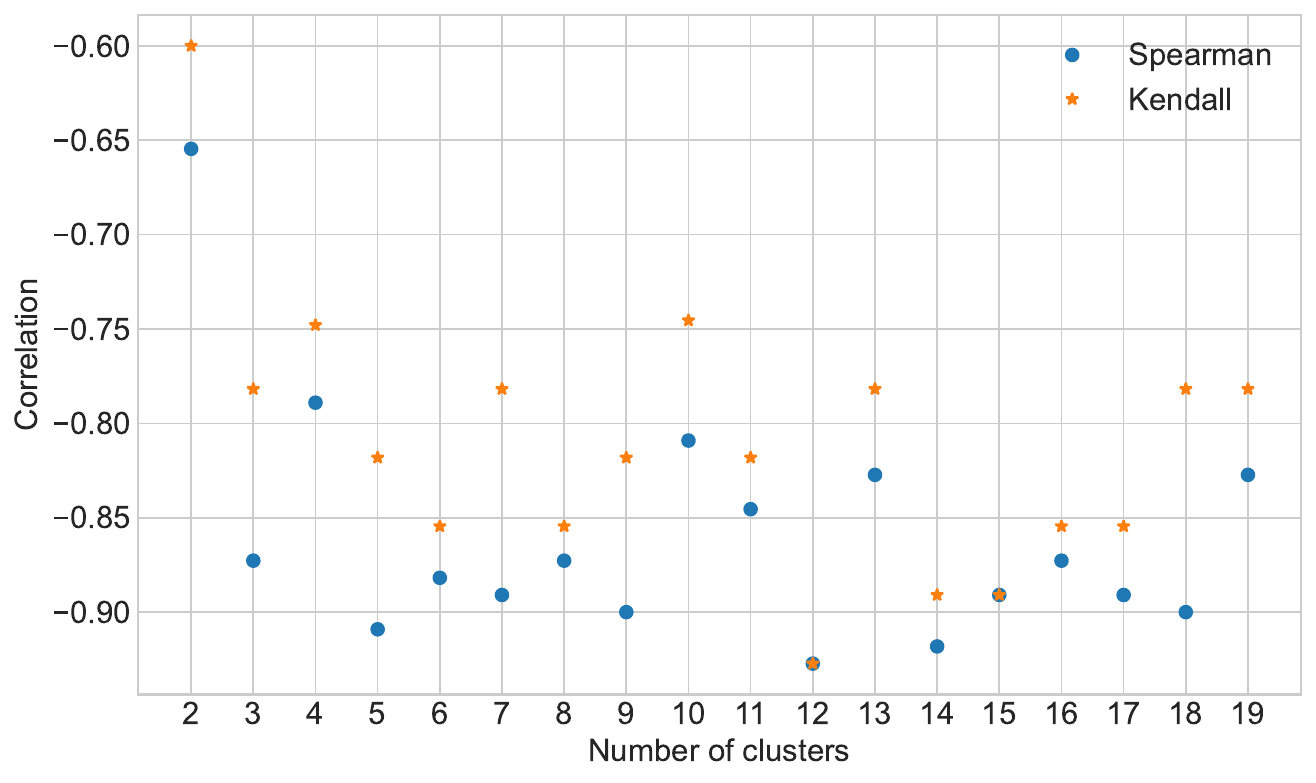}
				\subcaption{Even-sized clusters}\label{fig:even}
			\end{minipage}%
		\caption{Selection for the optimal number of clusters. \label{fig:cluster}}
	\end{figure}

	Figure \ref{fig:cluster} shows the correlations between the sale-wage discrepancy and efficiency using different numbers of clusters. The clusters are uneven in size in Figure \ref{fig:uneven} and even in size in Figure \ref{fig:even}. The uneven clusters are determined using the Jenks optimization method \citep{jenks} and the original values of wages and log values of sales. One advantage of uneven clusters is that the correlation is roughly monotonic with the number of clusters, as shown in Figure \ref{fig:uneven}. However, a disadvantage is that a large number of clusters is required to achieve a low, negative correlation, which is impractical due to the small amount of data. Even-sized clusters with ranks as the input, shown in Figure \ref{fig:even}, can achieve a low correlation with a relatively small number of clusters. Figure \ref{fig:cluster} reports both Kendall's and Spearman's rank correlations. We choose the first number of clusters that has a correlation less than $-0.8$. Based on Figure \ref{fig:even}, the suitable number of clusters is 5 or 6, depending on the type of the correlation used. As a robustness check, we also report the results using seven clusters in this e-companion.       
	
	\subsubsection{Estimation of transition matrices}
	
	After determining the number of clusters, we estimate the conditional probability kernels which are transition matrices of ranks. To improve the accuracy of these estimates, we use wage-sale data spanning from 2010 to 2021 for the 790 firms considered in Table \ref{tab:data}. We allow the planner to use all available data from 2010 to the final year of 2021 in the estimation process, in order to provide a larger pool of data. The transition matrices are calculated based on the frequencies of observed transitions between ranks over the years. 
	
	The transitions between ranks are crucial for our analysis. In the extreme case, $x_{t+1} = x_t$ and $y_{t+1} = y_t$ with probability one, that is, transition matrices are identity matrices. Then the state-dependent term will have no impact on the transport plans, because the marginal constraints have determined the matching already. There is no freedom to alter firm-wage pairs over time when rank transitions are forbidden. However, as long as transitions are not deterministic, we can capture $\alpha$ with enough data. 
	
\begin{figure}[H]
		\centering
			\begin{minipage}{0.45\textwidth}
				\centering
				\includegraphics[width=0.95\textwidth]{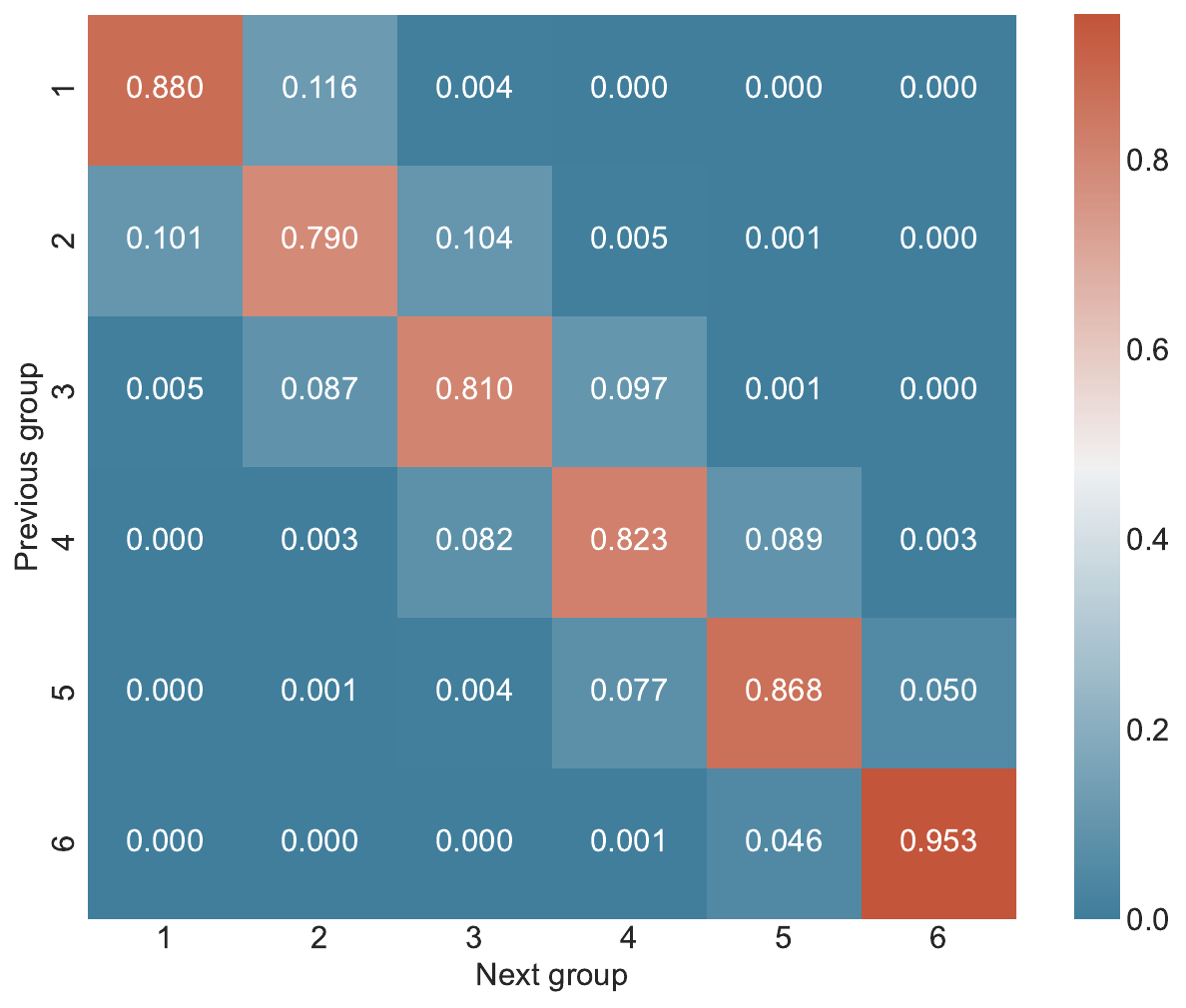}
				\subcaption{Sales}\label{fig:sale_trans}
			\end{minipage}
			\begin{minipage}{0.45\textwidth}
				\centering
				\includegraphics[width=0.95\textwidth]{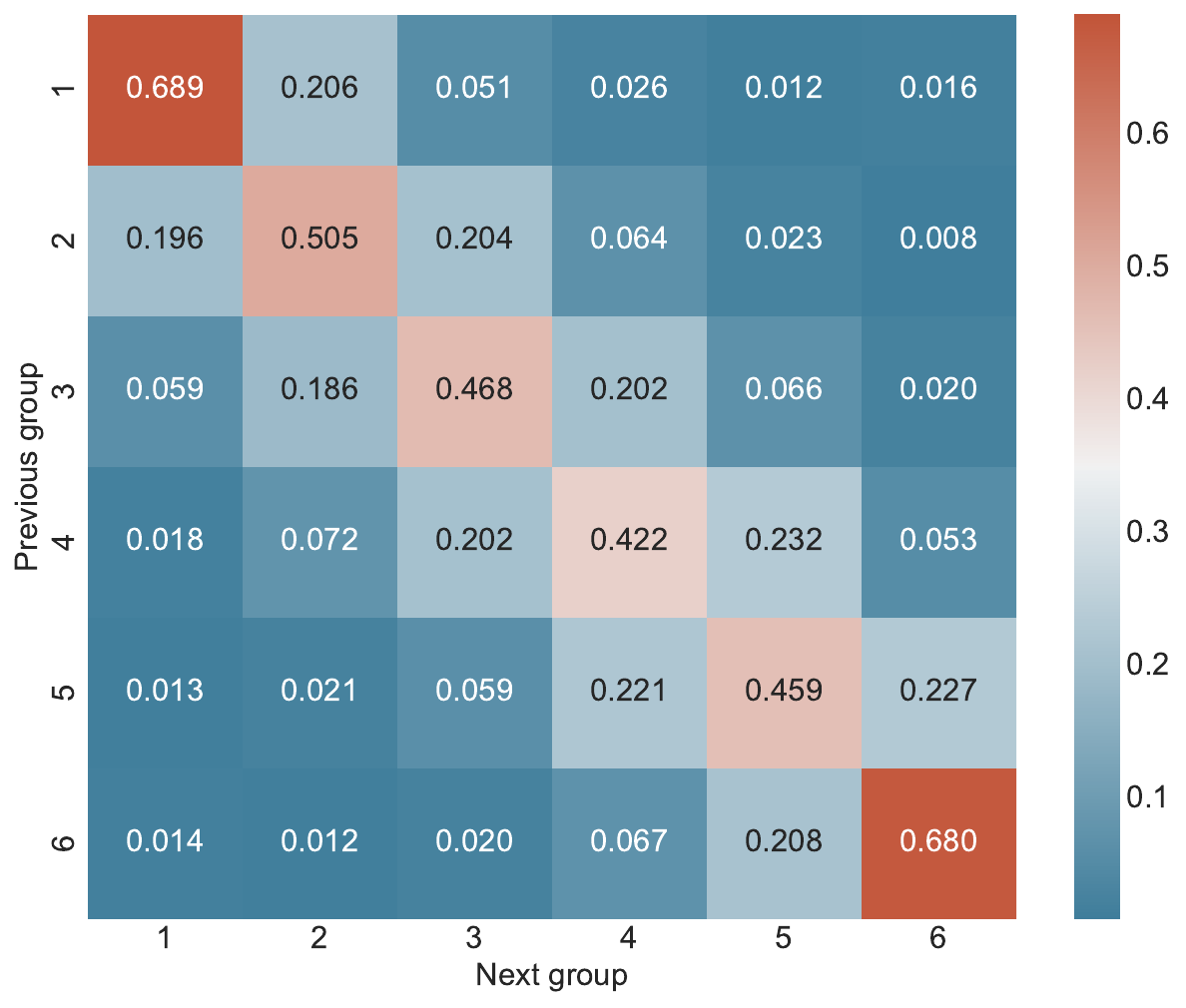}
				\subcaption{Wages}\label{fig:wage_trans}
			\end{minipage}%
		\caption{Transition matrices of wage and sale ranks. The matrices are averaged across 11 sectors. \label{fig:exec_trans}}
	\end{figure}

	Figure \ref{fig:exec_trans} gives the estimated transition matrices for sale and wage ranks. Group 1 is the group with the highest wages or sales. To save space, we average the transition matrices across 11 sectors and summarize them in two figures. A common pattern in all sectors is that sale ranks are more stable than wage ranks. Wages are more stable for the highest and the lowest group. For the wage groups in the middle, the probability of a rise or fall in wage ranks is approximately equal.
	
	\subsubsection{Model validation}\label{sec:valid}
	To stabilize the calibration, we adopt two additional methods. First, as our real dataset is of limited size, we use the re-sampling technique (bootstrapping) to generate a larger dataset by drawing with replacement from the real samples. Second, in practice, there are only a few non-zero entries in the one-period transport plans $\pi_t$. However, after multiplying the one-period transport plans as in equation \eqref{eq:decomp} to obtain $\pi(\alpha)$, there are many paths of $(x_{1:T}, y_{1:T})$ with very small probabilities. To speed up the calculation of the Wasserstein distance $\cW(\pi(\alpha), \pi_r)$ in \eqref{eq:wass_alpha} and amplify the impact of $\alpha$, we restrict the domain of $\pi(\alpha)$ to a set of paths with an equal size of the bootstrap samples and the highest probabilities.
	
	In our model, several approximations and hyperparameters may cause systematic bias, including the number of clusters, estimation errors in the transition matrices, the choice of the state-dependent preference function, the size of the bootstrap samples, and the restriction of $\pi(\alpha)$ on paths with the highest probabilities. To address concerns about potential bias in our model, we propose to first apply it to synthetic data consisting only of perfectly matched sale-wage pairs. That is, firms ranked number $n$ always pay wages ranked number $n$. In this case, $x_t = y_t$ for all time points $1 \leq t \leq T$. Ideally, we would expect to observe the optimal value of $\alpha$ to be close to zero in this scenario. There are several advantages to this validation procedure. First, it allows us to test whether the framework is functioning correctly in the simplest possible case. Second, it provides a method for selecting appropriate hyperparameters such that the benchmark outcome is close to zero. Finally, we can use the outcome of this analysis to correct any systematic bias in the model.

	\subsubsection{Parameter settings}
	In Section \ref{sec:exec}, the number of clusters is 6 and the preference function is the indicator function. We divided the ranks by the number of clusters $n$, such that $x_t$ and $y_t$ are in $[0, 1-1/n]$, with $0$ representing the highest rank. The candidate values for $\alpha$ are discrete and evenly spaced as $\{ -1.5 , -1.44, ..., 0, ..., 1.44, 1.5\}$, with a step size of $0.06$. The discount factor is set to $\beta = 0.9$. We ran the calibration with ten simulation instances. The size of both the synthetic and bootstrap samples is set to 500. We restricted the domain of $\pi(\alpha)$ to 500 paths with the highest probabilities.

\subsection{Robustness check with different number of clusters}\label{sec:robust}
As a robustness check, we examine the results in Section \ref{sec:exec} with different numbers of clusters. The screening method in Figure \ref{fig:even} suggests that the choice of five clusters is also reasonable. All other specifications are kept the same. The first row of Table \ref{tab:5} shows that the systematic bias moves further toward the negative side as the number of clusters decreases, indicating that the outcomes may be more similar to the one-cluster case. When resampled real data are used, we obtain the same conclusion as in the six-cluster situation. Table \ref{tab:5_test} finds that the correlations are negative and the $p$-values are smaller than 5\%.
\begin{table}[!ht]
	\centering
	{\scriptsize
		\begin{tabular}{cccccccccccc}
			\hline
			Sector & 10 & 15 & 20 & 25 & 30 & 35 & 40 & 45 & 50 & 55 & 60 \\
			\hline
			Benchmark $\alpha$ & $-0.42$ & $-0.12$ & $-0.318$ & $-0.324$ & $-0.282$ & $-0.366$ & $-0.42$ & $-0.21$ & $-0.252$ & $-0.222$ & $-0.294$ \\
			Raw $\alpha$ & 0.012 & $-0.072$ & $-0.006$ &  0.45 & $-0.132$ & $-0.51$ & $-0.444$ & $-0.066$ & $-0.492$ & $-0.462$ &  0.054 \\
			Adjusted $\alpha$ & 0.432 &  0.048 &  0.312 &  0.774 &  0.15 & $-0.144$ & $-0.024$ &  0.144 & $-0.24$ & $-0.24$ & 0.348 \\
			\hline 
		\end{tabular}
	}
	\caption{Mean values of the optimal $\alpha$ using five clusters. Other settings are the same as in the six-cluster counterpart.}\label{tab:5}
\end{table}

\begin{table}[!ht]
	\centering
	\begin{tabular}{ccc}
		\hline
		Correlation & Spearman ($p$-value) & Kendall ($p$-value) \\
		\hline 
		Raw $\alpha$ & $-0.709$ $(0.015)$ &  $-0.636$ $(0.006)$ \\
		\hline
		Adjusted $\alpha$ & $-0.700$ $(0.016)$ & $-0.636$ $(0.006)$ \\
		\hline
	\end{tabular}
	\caption{Testing the association between job market efficiency and inertia when five clusters are adopted. }\label{tab:5_test}
\end{table}

When seven clusters are used, the systematic bias shifts to the positive side, as shown in the first row of Table \ref{tab:7}. The original and adjusted optimal $\alpha$ tend to be larger than those obtained using five or six clusters. Although the correlations in Table \ref{tab:7_test} are still negative, the power of the test is much lower. One possible explanation is that Figure \ref{fig:even} shows that the choice of seven clusters results in a slightly worse correlation than the five and six clusters options. In addition, estimation errors in the transition matrices can be larger. Unfortunately, using a longer period of empirical data to estimate the transition matrices may not be a viable solution, as the job market may not be stationary over a long time horizon. Our parameter settings were fine-tuned for six clusters, so it is expected that the performance can be worse with a different number of clusters.
\begin{table}
	\centering
	{\small
		\begin{tabular}{cccccccccccc}
			\hline
			Sector & 10 & 15 & 20 & 25 & 30 & 35 & 40 & 45 & 50 & 55 & 60 \\
			\hline
			Benchmark $\alpha$ & 0.234 & $-0.06$ &  0.33 &  0.324 & $-0.06$ &  0.582 & $-0.096$ &  0.198 & $-0.102$ &  0.066 &  0.078  \\
			Raw $\alpha$ &  0.372 &  $-0.06$ &  0.756 &  0.42 &  0.222 &  0.636 &  0.054 &  0.438 & 0.204 &  0.264 &  0.666 \\
			Adjusted $\alpha$ &0.138 & 0.0  &  0.426 & 0.096 & 0.282 & 0.054 & 0.15 & 0.24 & 0.306 & 0.198 & 0.588 \\
			\hline 
		\end{tabular}
	}
	\caption{Mean values of the optimal $\alpha$ using seven clusters. Other settings are the same as in the six-cluster counterpart.}\label{tab:7} 
\end{table}

\begin{table}
	\centering
	\begin{tabular}{ccc}
		\hline
		Correlation & Spearman ($p$-value) & Kendall ($p$-value) \\
		\hline 
		Raw $\alpha$ & $-0.345$ $(0.298)$ &  $-0.382$ $(0.121)$ \\
		\hline
		Adjusted $\alpha$ & $-0.245$ $(0.467)$ & $-0.309$ $(0.218)$ \\
		\hline
	\end{tabular}
	\caption{Testing the association between job market efficiency and inertia when seven clusters are adopted. }\label{tab:7_test}
\end{table}

\section{Supplement to the academic job market}\label{sec:supp_acad}
\subsection{Summary statistics}
 The University of California (UC) compensation data are from the Government Compensation in California website [\footnote{\url{https://publicpay.ca.gov/Reports/RawExport.aspx}}]. The dataset includes total annual wages for various positions. In our analysis, we focus on four specific full-time positions: ``Prof-Ay-B/E/E", ``Assoc Prof-Ay-B/E/E", ``Asst Prof-Ay-B/E/E", and ``Postdoc-Employee". ``AY" is short for Academic Year and ``B/E/E" means Business/Economics/Engineering. In our current dataset, these three areas are reported jointly. The employees under the ``Postdoc-Employee" title are from all departments. These full-time positions may provide a better representation of compensation levels than part-time positions. Tables \ref{tab:prof} -- \ref{tab:postdoc} present the summary statistics of salaries in 2017 -- 2021 for nine universities in the UC system, except UC San Francisco since it focuses on medical research and does not have faculty in the B/E/E departments. These tables also include the number of employees under each job title for 2021 only.
	
	In assessing university quality, we utilize the U.S. News rankings, based on historical data compiled by Andrew G. Reiter [\footnote{\url{https://andyreiter.com/datasets/}}]. While acknowledging the methodological limitations of university rankings, particularly the subjectivity of reputational assessments, our focus is specifically on the rankings of nine universities within the UC system. It is important to highlight that faculty salaries contribute to these rankings but carry a weight of only 7\% [\footnote{ More details in the article at \url{https://www.usnews.com/education/best-colleges/articles/how-us-news-calculated-the-rankings}}]. Therefore, the ranks of faculty wages and university rankings do not align automatically. To determine the annual pay levels at each university, we calculate a single value by first ranking all wage observations among a given year and job title, and then selecting the median rank as a measure of compensation levels. This method mitigates the impact of outliers in the data.

\begin{table}
	\centering
	\begin{tabular}{ccccc}
		\hline
		University & Number of employees & Min wage & Median wage & Max wage\\
		\hline
		UC Berkeley    & 222 & 3000.0 & 250141.5 &  632257.0\\
		UC  Davis         &  164 & 1250.0 & 211538.0 & 472937.0 \\
		UC Irvine          & 172 & 2117.0 & 225990.0 & 609000.0 \\
		UC Los Angeles  & 195 & 1906.0 & 330058.0 & 778102.0 \\
		UC  Merced        & 29 & 19401.0 & 207291.5 & 350160.0 \\
		UC Riverside       &  72 & 500.0 &  230418.5 & 412552.0 \\
		UC San Diego        &  181 & 12400.0 & 246321.0 & 586884.0 \\ 
		UC Santa Barbara  &  114 & 21667.0 &  258774.0 & 580967.0 \\
		UC Santa Cruz        &  66 & 6936.0 & 182750.0 & 379110.0 \\
		\hline 
	\end{tabular}
	\caption{Wage statistics of B/E/E professors in 2017 -- 2021. The numbers of employees are in 2021 only. Min, median, and max wages (in U.S. dollars) are for five years from 2017 to 2021.}\label{tab:prof}
\end{table}

\begin{table}[!ht]
	\centering
	\begin{tabular}{ccccc}
		\hline
		University & Number of employees & Min wage & Median wage & Max wage\\
		\hline
		UC Berkeley    & 73 & 11086.0 &  209901.0 & 582942.0 \\
		UC  Davis         &   49 & 54768.0 & 160465.0 & 313901.0 \\
		UC Irvine          &  64 & 45856.0 & 171821.0 &  387084.0 \\
		UC Los Angeles  & 40 & 4057.0 &  261654.0 & 494989.0 \\
		UC  Merced        &  27 & 62208.0 & 155763.0 & 230762.0 \\
		UC Riverside       &   44 & 10608.0 &  175183.0 &  301166.0 \\
		UC San Diego        &   60 & 11258.0  & 199728.0 & 436367.0 \\ 
		UC Santa Barbara  &  13 & 36633.0 &  176921.0 & 430169.0 \\
		UC Santa Cruz        &  22 & 10667.0 & 163482.0 & 235600.0 \\
		\hline 
	\end{tabular}
	\caption{Wage statistics of B/E/E associate professors. Numbers are reported in the same way as in Table \ref{tab:prof}.}\label{tab:assocprof}
\end{table}

\begin{table}[!ht]
	\centering
	\begin{tabular}{ccccc}
		\hline
		University & Number of employees & Min wage & Median wage & Max wage\\
		\hline
		UC Berkeley    & 72 & 9992.0 &   172602.5 & 356679.0\\
		UC  Davis         &  51 & 1000.0 & 141866.0 & 262833.0 \\
		UC Irvine          &  74 & 10833.0 &   143116.0 & 331450.0 \\
		UC Los Angeles  & 82 & 8917.0 &  175874.5 & 408104.0 \\
		UC  Merced        &  31 & 8450.0 & 129069.0 & 203730.0\\
		UC Riverside       &   58 & 27775.0 &   132105.5 &  266866.0\\
		UC San Diego        &  97 & 2386.0 &  155516.0 &  429625.0\\ 
		UC Santa Barbara  &  48 & 1295.0 & 162483.0 &  302019.0\\
		UC Santa Cruz        &  45 & 9017.0 &  138514.0 &   247281.0 \\
		\hline 
	\end{tabular}
	\caption{Wage statistics of B/E/E assistant professors. Numbers are reported in the same way as in Table \ref{tab:prof}.}\label{tab:asstprof}
\end{table}

\begin{table}[!ht]
	\centering
	\begin{tabular}{ccccc}
		\hline
		University & Number of  employees & Min wage & Median wage & Max wage \\
		\hline
		UC Berkeley    & 1246 & 4.0 &  44151.0 & 197419.0 \\
		UC  Davis         &   933 & 2.0 & 45978.0 &  137724.0 \\
		UC Irvine          & 524 & 40.0 &  43032.0 &  102079.0 \\
		UC Los Angeles  & 1056 & 10.0 & 48566.0 & 181982.0 \\
		UC  Merced        &  89 & 1694.0 & 39271.5 & 100772.0\\
		UC Riverside       &    299 & 42.0 &  44651.0 & 87372.0 \\
		UC San Diego        &  1394 & 14.0 &  46350.0 &  139683.0 \\ 
		UC Santa Barbara  &     430 & 54.0 &  47104.0 & 148142.0 \\
		UC Santa Cruz        &  184 & 106.0 &  45629.5 & 80133.0 \\
		\hline 
	\end{tabular}
	\caption{Wage statistics of postdocs. Numbers are reported in the same way as in Table \ref{tab:prof}.}\label{tab:postdoc}
\end{table}

\subsection{Number of clusters and estimation of transition matrices}
	To determine the number of clusters for our analysis, we consider the U.S. News rankings and domain knowledge, as the small number of universities in our dataset (nine) makes it difficult to use the criterion applied to the executive data. We find that there are usually two universities with closely ranked positions, with the top two universities in the UC system ranking around the 20th place among national universities (and sometimes tied). The other universities tend to rank around the 30th, 40th, 80th, and 90th places. Based on these observations, we set the number of clusters to five, with two universities in each of the top four clusters and one university in the last cluster.

\begin{figure}[H]
	\centering
	\includegraphics[width=0.4\textwidth]{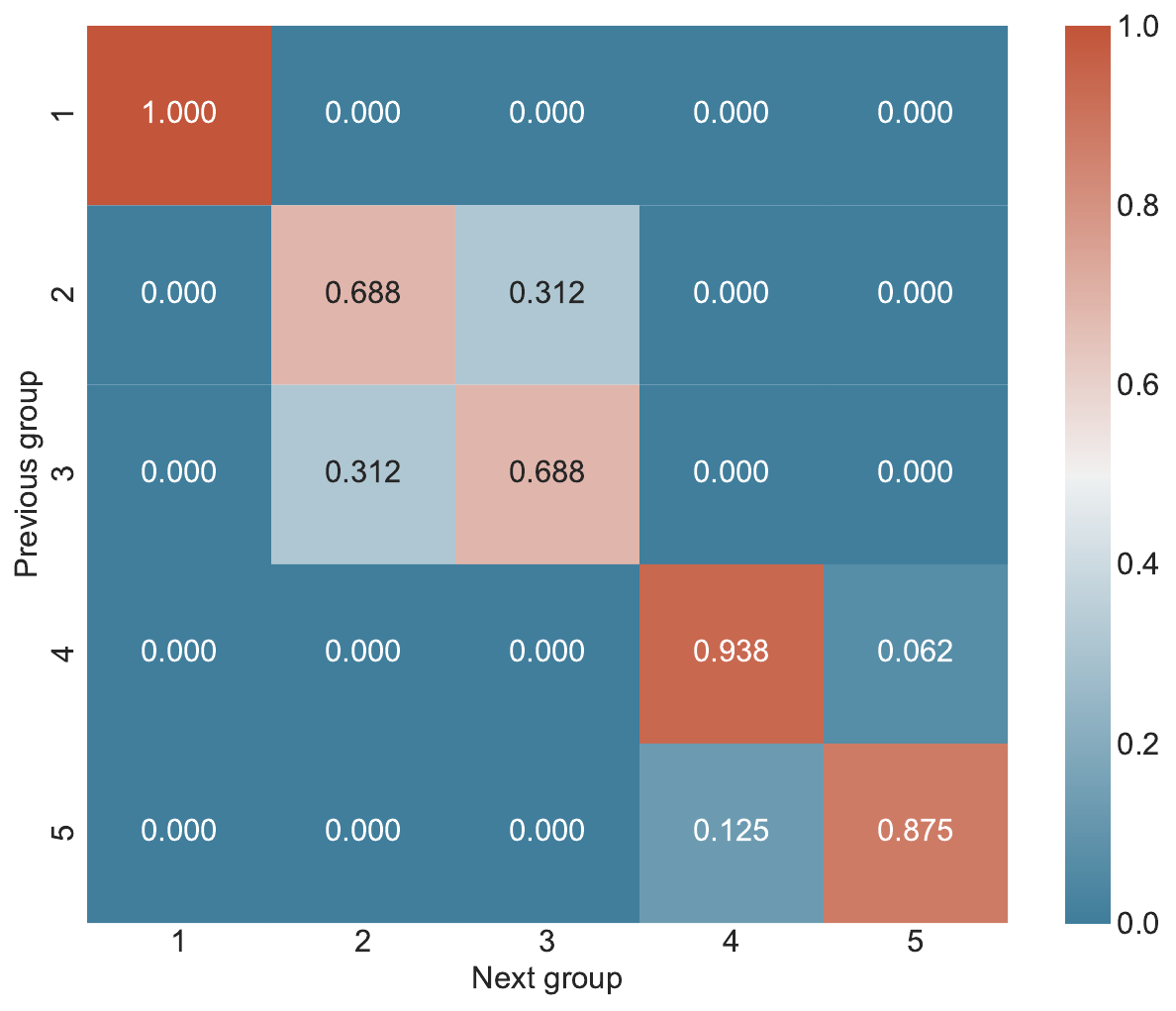}
	\caption{Transition matrices of university rankings. \label{fig:uni_trans}}
\end{figure}

\begin{figure}[!ht]
	\caption{Transition matrices of wages.} \label{fig:prof_trans}
	\centering
	\begin{minipage}{0.45\textwidth}
		\centering
		\includegraphics[width=0.9\textwidth]{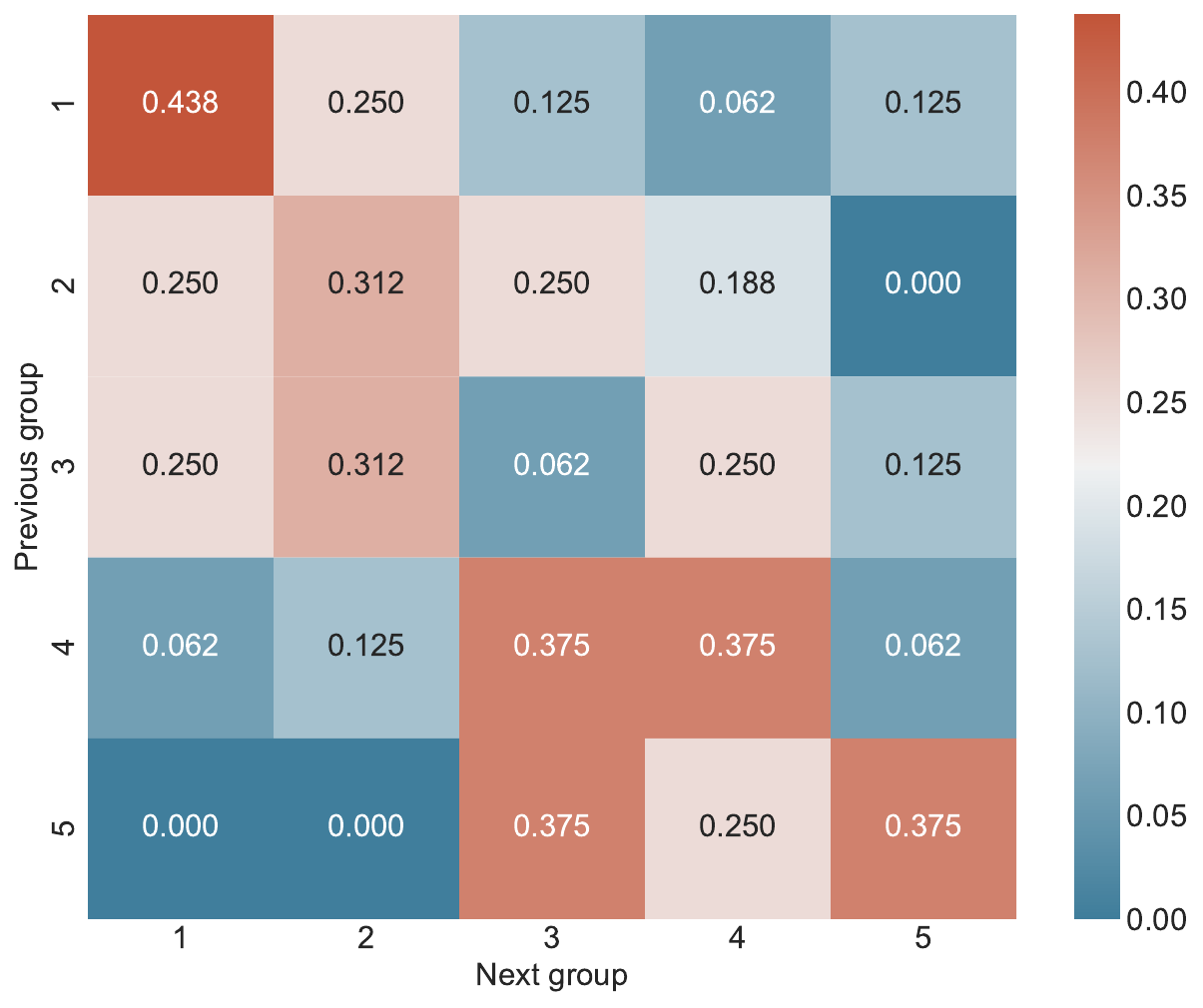}
		\subcaption{Postdoc}\label{fig:postdoc_trans}
	\end{minipage}
	\begin{minipage}{0.45\textwidth}
		\centering
		\includegraphics[width=0.9\textwidth]{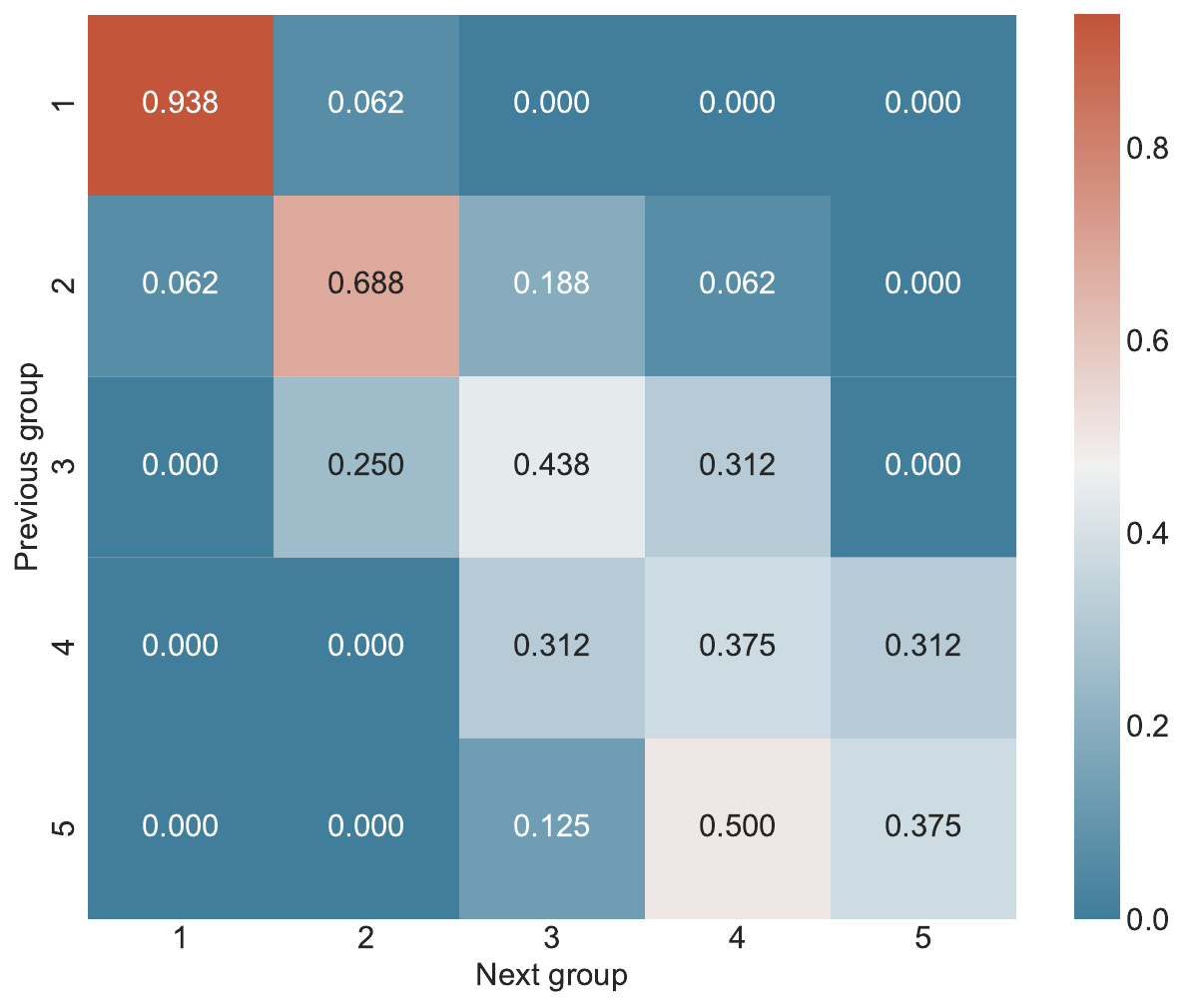}
		\subcaption{Assistant professor}\label{fig:assitprof_trans}
	\end{minipage}
	\begin{minipage}{0.45\textwidth}
		\centering
		\includegraphics[width=0.9\textwidth]{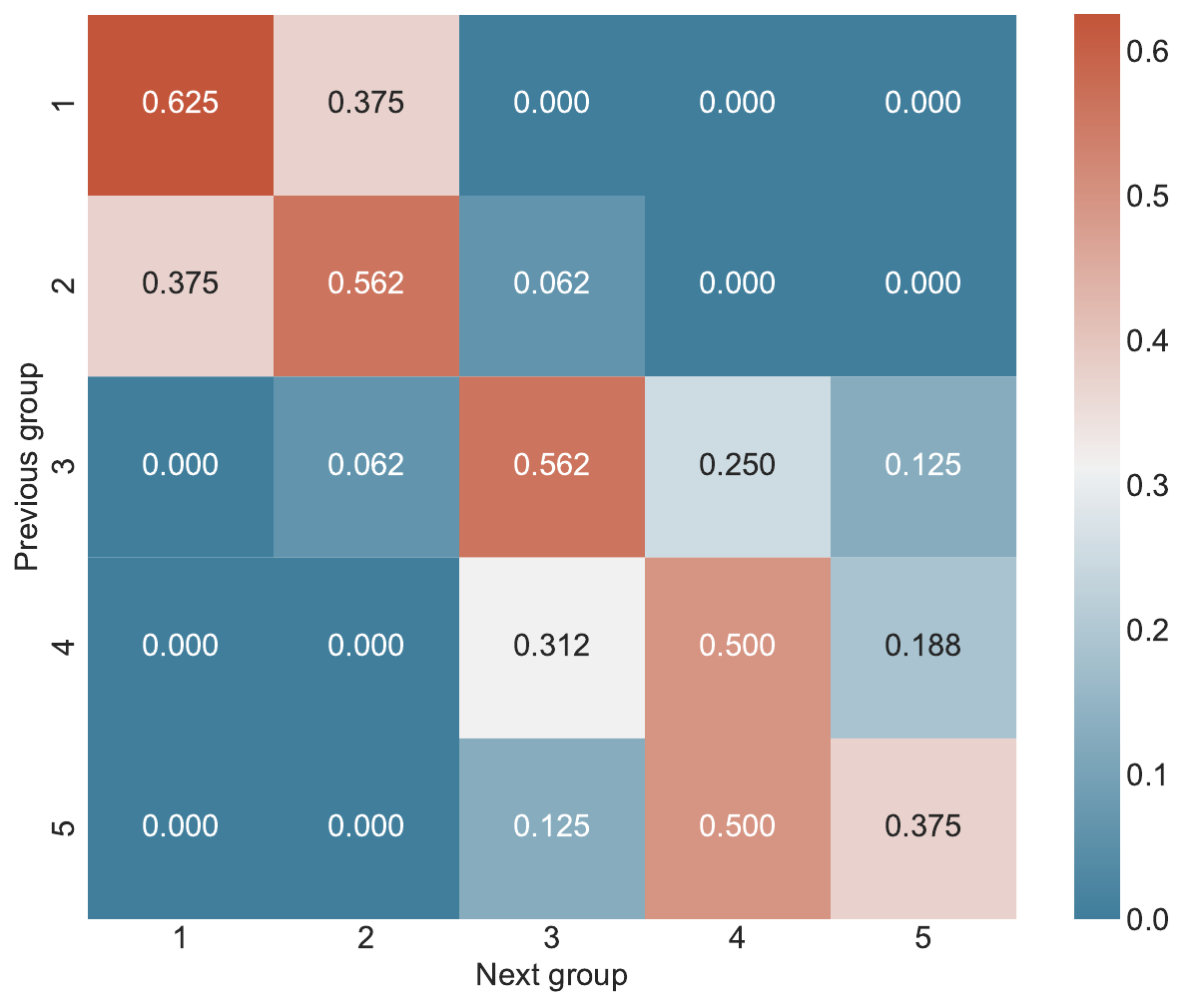}
		\subcaption{Associate professor}\label{fig:assocprof_trans}
	\end{minipage}
	\begin{minipage}{0.45\textwidth}
		\centering
		\includegraphics[width=0.9\textwidth]{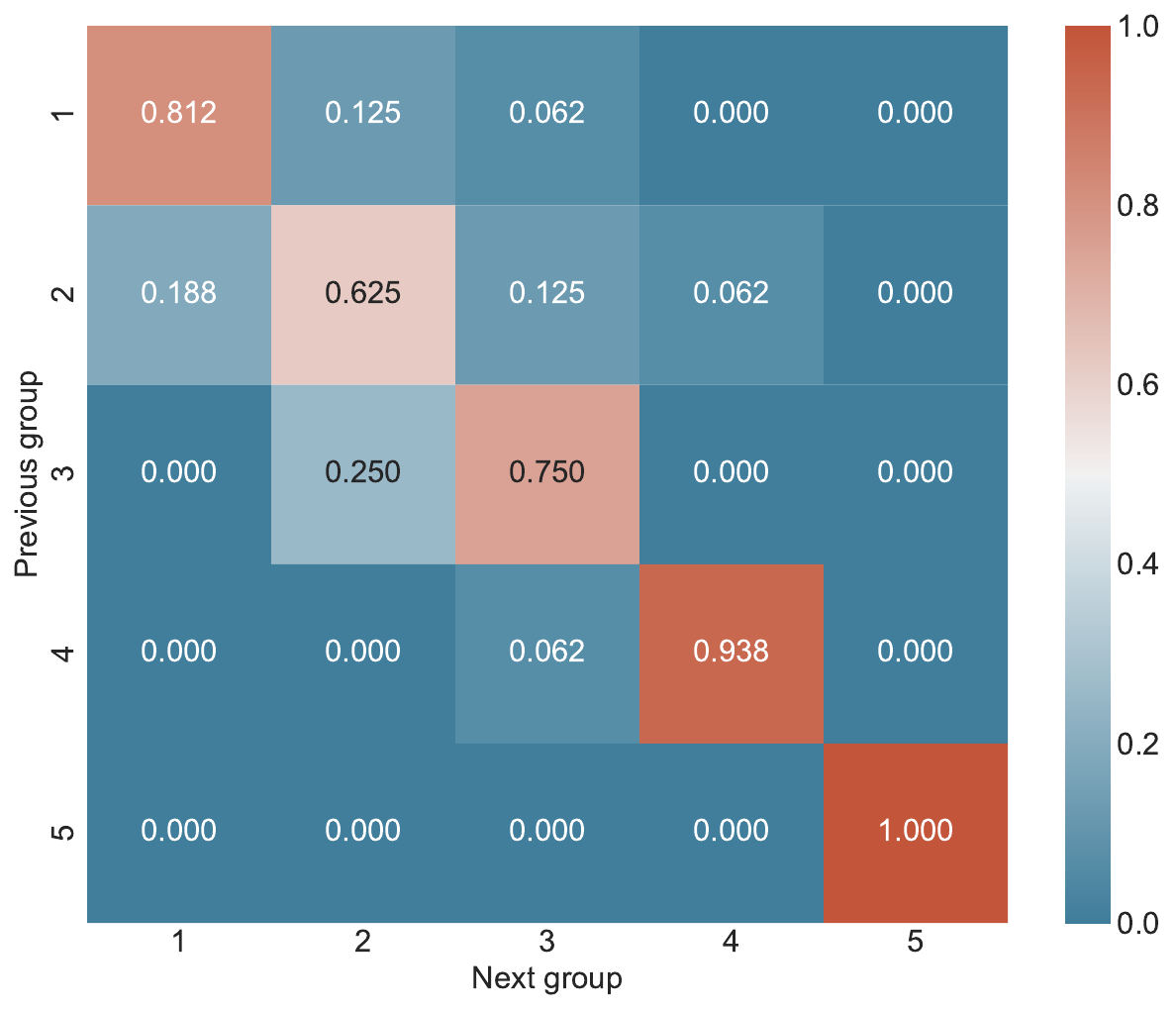}
		\subcaption{Professor}\label{fig:full_trans}
	\end{minipage}%
\end{figure}

To estimate the transition matrices of wages and university rankings, we use data from 2013 to 2021, the longest period available for UC employee salaries on the Government Compensation in California website. Figure \ref{fig:uni_trans} shows that university rankings are stable in this time horizon. Most university ranking transitions happen between the second and the third group. For wages, Figure \ref{fig:prof_trans} demonstrates that postdoc wage ranks are relatively unstable over time. The transition matrices become more and more stable when the academic job ranks become higher. The first group with the highest wages is stable among tenure-track and tenured faculty, while it is unstable for postdoc employees.

\subsection{Parameter settings}
Since the number of universities in our dataset is much smaller than the number of firms, we set the size of our bootstrap samples to 200 for each simulation. We use the state-dependent preference function with $\tau = 1$.

	\subsection{Calibration curves for professor-level positions}
	
\begin{figure}
	\centering
		\begin{minipage}{0.48\textwidth}
			\centering
			\includegraphics[width=0.9\textwidth]{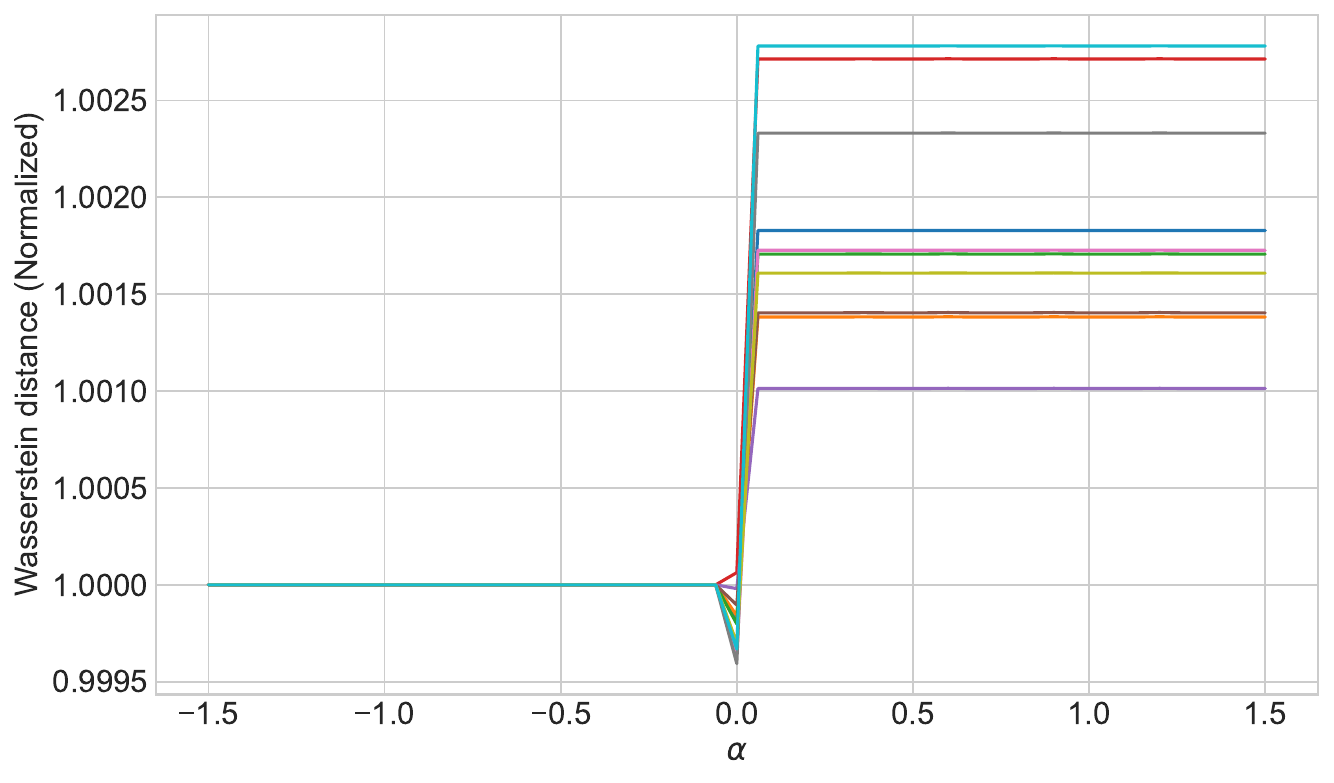}
			\subcaption{Synthetic data with perfect matching}\label{fig:prof_zero}
		\end{minipage}
		\begin{minipage}{0.48\textwidth}
			\centering
			\includegraphics[width=0.9\textwidth]{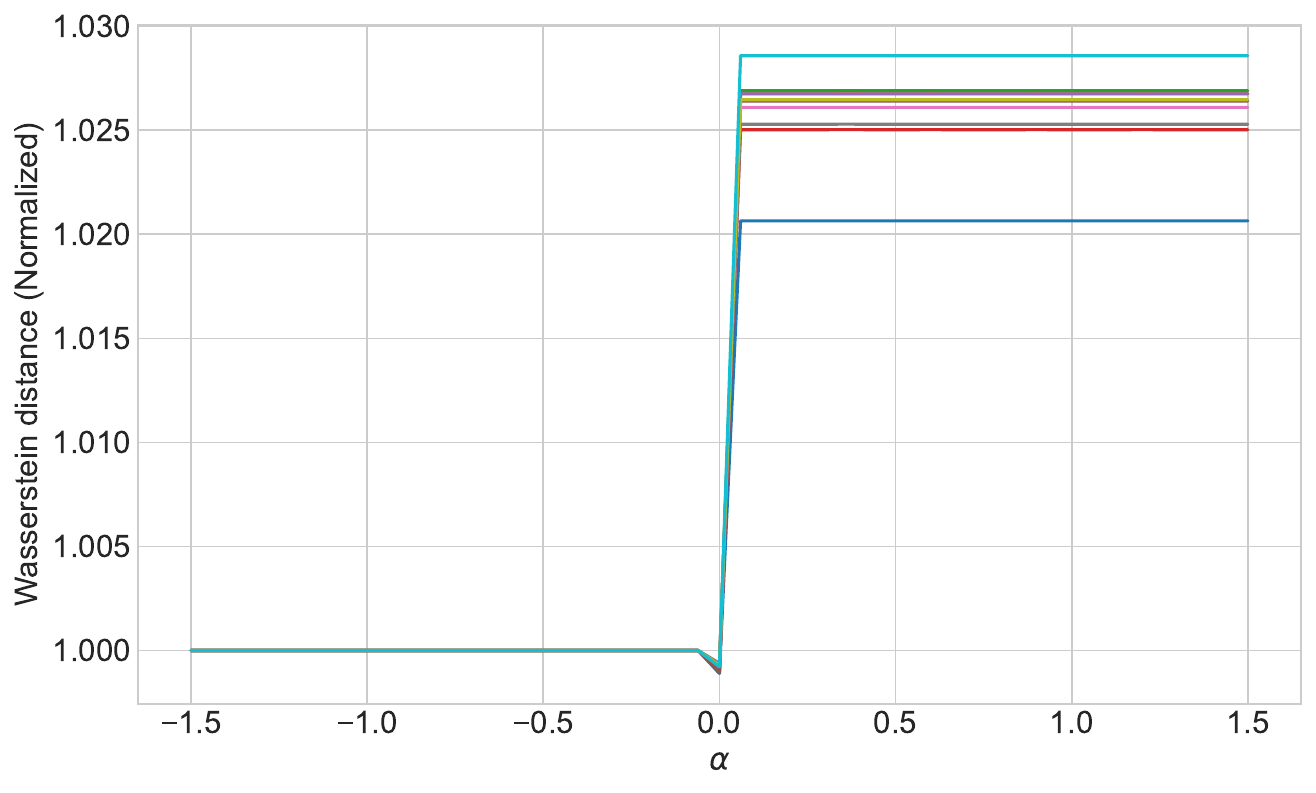}
			\subcaption{Resampled real data}\label{fig:prof_real}
		\end{minipage}%
	\caption{Calibration curves for B/E/E professors. In each subplot, ten curves represent individual calibrations with ten independent sampled data. \label{fig:prof}}
\end{figure}

	Using the B/E/E Professors data, Figure \ref{fig:prof} plots the Wasserstein distance between the model-implied transport plan $\pi(\alpha)$ and the empirical transport plan $\pi_r$. The Wasserstein distance is divided by $\cW(\pi(-1.5), \pi_r)$ such that all curves start from 1 when $\alpha = -1.5$. As a validation test and benchmark, Figure \ref{fig:prof_zero} calibrates $\alpha$ to synthetic data with perfect matching and plots curves from ten simulation instances. It shows that the optimal $\alpha$ is zero and there is no evidence of inertia when synthetic data with perfect matching are used. In contrast, Figure \ref{fig:prof_real} shows the results using bootstrap samples. The optimal $\alpha$ is still zero across different simulation runs, indicating that the real data for the professor category are very similar to the perfectly matched synthetic data. Table \ref{tab:uc_alph} reports the mean values of optimal $\alpha$ across these ten simulations. 
	
\begin{figure}
	\centering
	\begin{minipage}{0.48\textwidth}
		\centering
		\includegraphics[width=0.9\textwidth]{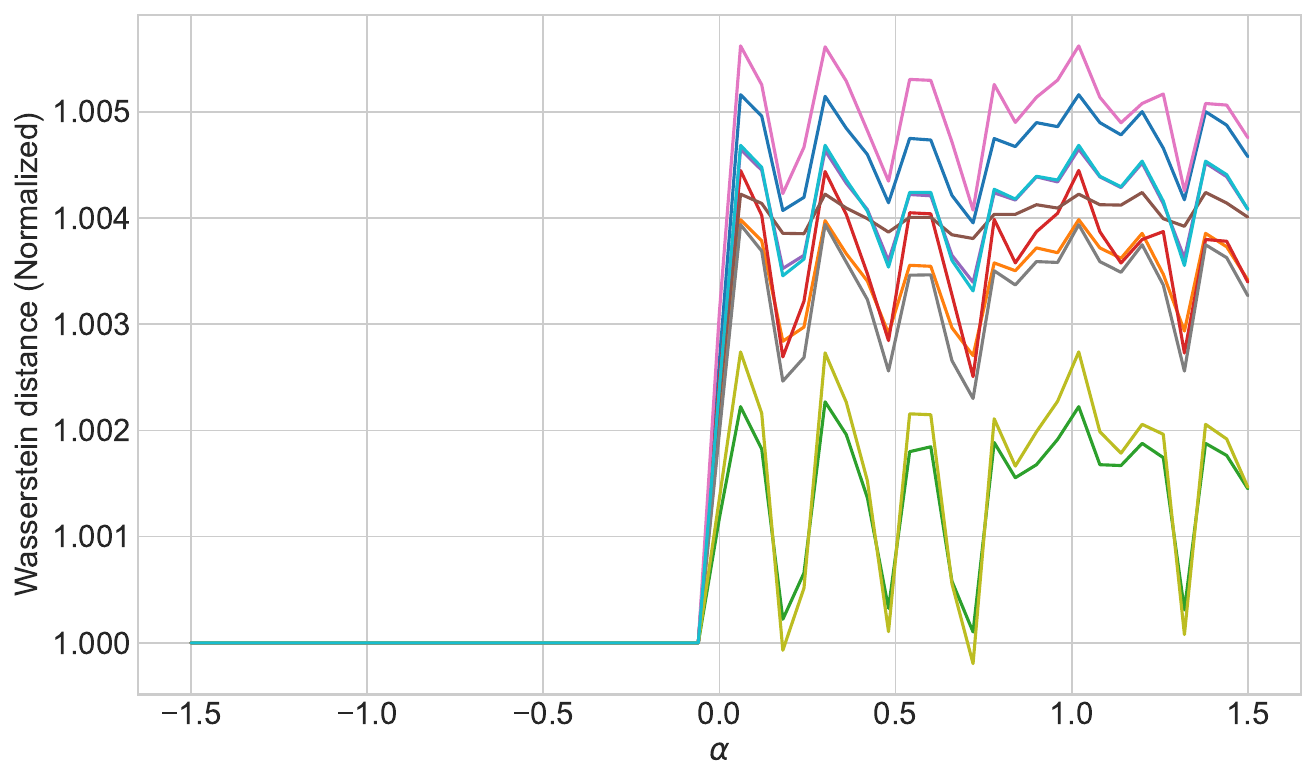}
		\subcaption{Synthetic data with perfect matching}\label{fig:asst_zero}
	\end{minipage}
	\begin{minipage}{0.48\textwidth}
		\centering
		\includegraphics[width=0.9\textwidth]{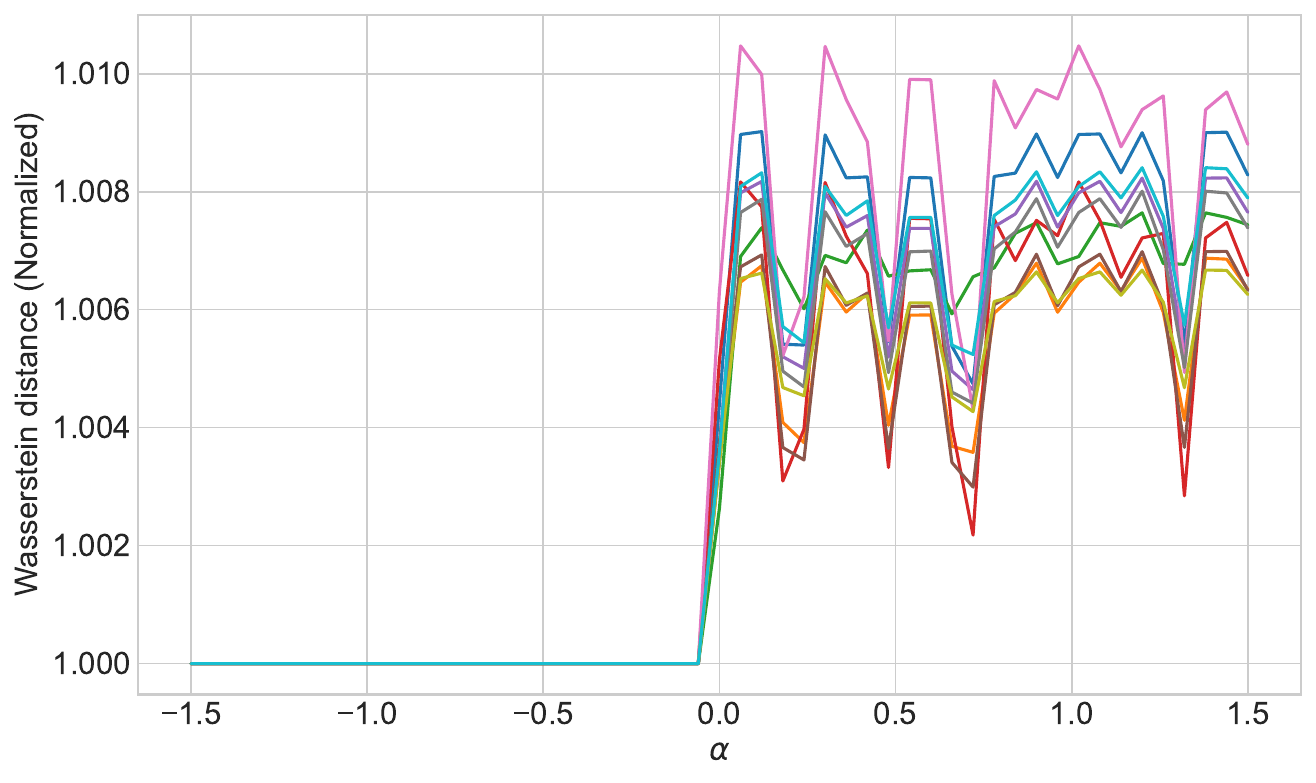}
		\subcaption{Resampled real data}\label{fig:asst_real}
	\end{minipage}%
	\caption{Calibration curves for B/E/E assistant professors. }\label{fig:asst}
\end{figure}

\begin{figure}
	\centering
	\begin{minipage}{0.48\textwidth}
		\centering
		\includegraphics[width=0.9\textwidth]{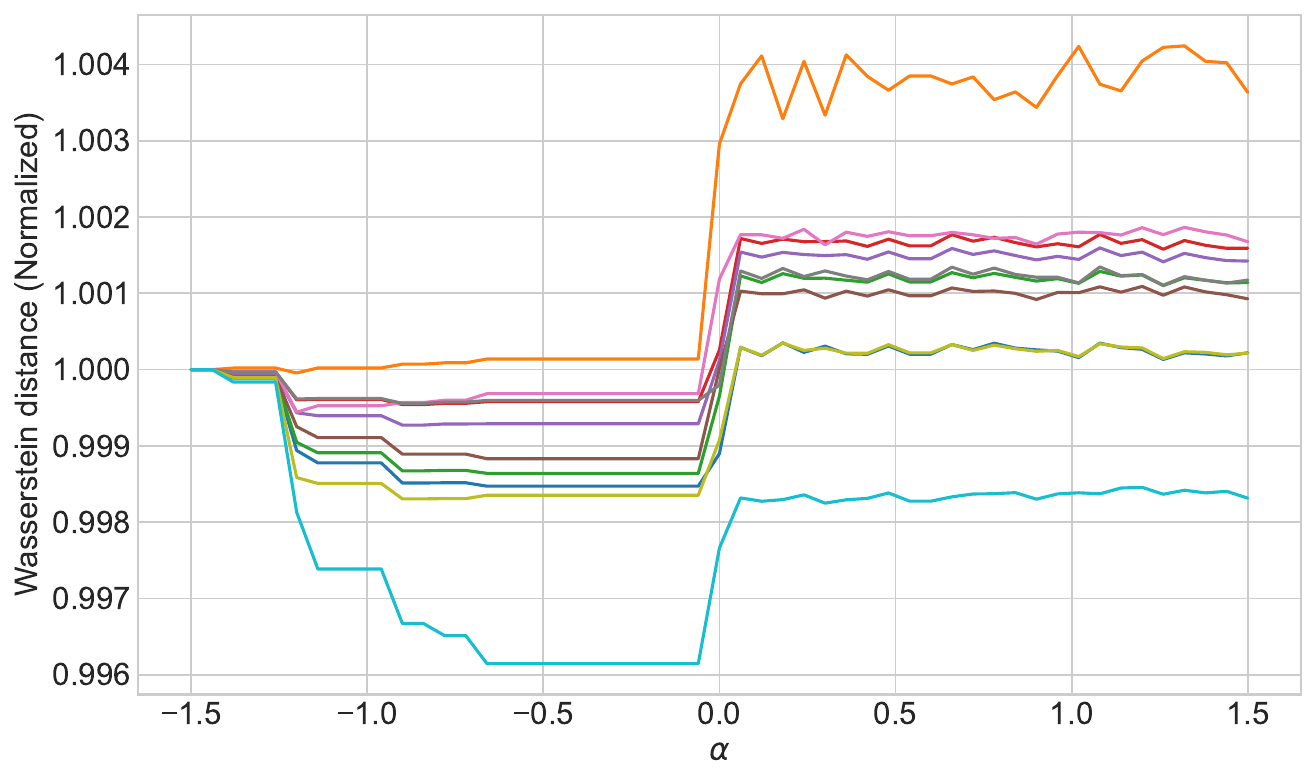}
		\subcaption{Synthetic data with perfect matching}\label{fig:assoc_zero}
	\end{minipage}
	\begin{minipage}{0.48\textwidth}
		\centering
		\includegraphics[width=0.9\textwidth]{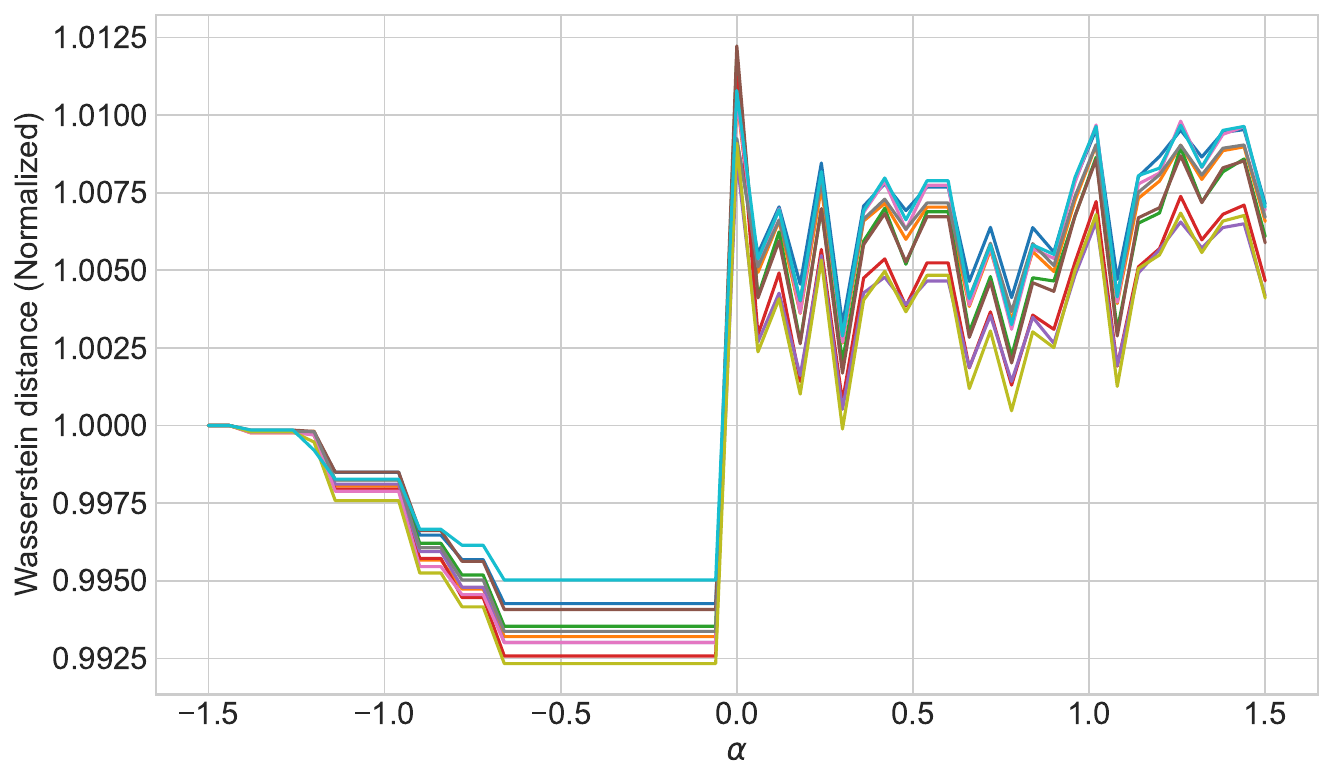}
		\subcaption{Resampled real data}\label{fig:assoc_real}
	\end{minipage}%
	\caption{Calibration curves for B/E/E associate professors.} \label{fig:assoc}
\end{figure}

	We observe similar curves for associate professors and assistant professors. Their wages are also highly matched with university rankings, with no significant evidence of inertia. We report the mean values of the optimal $\alpha$ in Table \ref{tab:uc_alph} and show the calibration curves in Figures \ref{fig:assoc} -- \ref{fig:asst}.

Figure \ref{fig:postdoc} also shows that the Wasserstein distance curves can be noisy when $\alpha$ varies. To estimate the optimal $\alpha$ robustly, we filter the original curves using the one-dimensional total variation denoising algorithm \citep{condat2013direct}. This algorithm is suitable for signals with piecewise constant behaviors. Table \ref{tab:uc_alph} reports the mean optimal $\alpha$ calculated from the filtered data. There is no significant evidence of inertia in wages for faculty. In contrast, the effect of inertia is very strong in the postdoc sector. 
	\begin{table}[H]
		\centering
		\begin{tabular}{ccccc}
			\hline
			Position & Professor & Associate Professor & Assistant Professor & Postdoc \\
			\hline
			Benchmark $\alpha$  & 0.0 & $-0.06$ & $-0.06$ & 0.0 \\
			Raw $\alpha$ & 0.0 &  $-0.06$ & $-0.06$ & 1.194 \\
			Adjusted $\alpha$ & 0.0 & 0.0 & 0.0 & 1.194 \\
			\hline 
		\end{tabular}
		\caption{Mean values of optimal $\alpha$ after total variation filtering.}\label{tab:uc_alph}
	\end{table}
	We test the connection between the inertia and university ranking-wage correlations, summarized in Table \ref{tab:uc_finaltest}. Correlations in Table \ref{tab:uc_finaltest} are negative, showing that when the inertia is stronger, wages are less matched with university rankings. However, the power of our tests is low since only four types of jobs are considered.  
	\begin{table}[H]
		\centering
		\begin{tabular}{ccc}
			\hline
			Correlation & Spearman ($p$-value) & Kendall ($p$-value) \\
			\hline 
			Raw $\alpha$ & $-0.632$ (0.368) &  $-0.548$ (0.279) \\
			\hline
			Adjusted $\alpha$ & $-0.775$ (0.225) & $-0.707$ (0.180) \\
			\hline
		\end{tabular}
		\caption{Correlations between optimal $\alpha$ and the ranking-wage dependence. }\label{tab:uc_finaltest}
	\end{table}

\end{document}